\documentclass[12pt,openany]{amsbook}
\usepackage{epsfig, amscd}

\numberwithin{equation}{chapter}
\begingroup
\newtheorem{propo}[equation]{Proposition}
\newtheorem{corol}[equation]{Corollary}
\newtheorem{theor}[equation]{Theorem}
\newtheorem{lemma}[equation]{Lemma}
\newtheorem{genpr}[equation]{General Problem}
\endgroup
\theoremstyle{definition}
\newtheorem{defin}[equation]{Definition}
\newtheorem{notat}[equation]{Notation}
\newtheorem{examp}[equation]{Example}
\newtheorem{exerc}[equation]{Exercise}
\newtheorem{remar}[equation]{Remark}
\newtheorem{claim}[equation]{Claim}
\newtheorem{quest}[equation]{Question}
\newtheorem{prope}[equation]{Properties}

\newenvironment{probl}{\ \\ \indent\textsc{Problem:}\ }{\ \\}
\newcommand{\titulo}[1]{\hfill\newline\textbf{#1}\ }
\newenvironment{reftc}{\ \\ \noindent\textbf{References for this chapter:}\ }{\ \\}

\newcommand{\al}{\alpha}
\newcommand{\Gm}{\Gamma}
\newcommand{\gm}{\gamma}
\newcommand{\ka}{\kappa}
\newcommand{\lm}{\lambda}
\newcommand{\ep}{\epsilon}
\newcommand{\om}{\omega}
\newcommand{\Om}{\Omega}
\newcommand{\sg}{\sigma}
\newcommand{\Dt}{\Delta}
\newcommand{\vlm}{\vec\lambda}

\newcommand{\RR}{\mathbb{R}}
\newcommand{\CC}{\mathbb{C}}
\newcommand{\ZZ}{\mathbb{Z}}
\newcommand{\Ss}{\mathbb{S}}
\newcommand{\Oc}{\mathcal{O}}
\newcommand{\Cc}{\mathcal{C}}
\newcommand{\hh}{\mathfrak{h}}
\newcommand{\ptl}{\partial}
\newcommand{\tr}{\operatorname{tr}}
\newcommand{\const}{\operatorname{const}}
\newcommand{\vol}{\operatorname{vol}}
\newcommand{\sol}{\operatorname{Sol}}
\newcommand{\GL}{\operatorname{GL}}
\newcommand{\SL}{\operatorname{SL}}
\newcommand{\Aut}{\operatorname{Aut}}
\newcommand{\oco}{\otimes\cdots\otimes}
\newcommand{\oij}{\Omega^{(ij)}}
\newcommand{\I}[1]{\operatorname{I}^{\gm}_{#1}}
\newcommand{\smd}{\operatorname{Sing} M^{\otimes m}[|m|-2]}
\newcommand{\smk}{\operatorname{Sing} M^{\otimes m}[|m|-2k]}
\newcommand{\smkq}{\operatorname{Sing} M^{\otimes m,q}[|m|-2k]}

\newcommand{\sld}{\operatorname{Sing} L^{\otimes m}[|m|-2]}
\newcommand{\slk}{\operatorname{Sing} L^{\otimes m}[|m|-2k]}
\newcommand{\pg}{\varphi^{(\gm)}}
\newcommand{\nup}[2]{(#1_1,\ldots,#1_{#2})}
\newcommand{\nupp}[3]{(#1_{#2},\ldots,#1_{#3})}
\newcommand{\dkf}[2]{d#1_1\wedge\cdots\wedge d#1_{#2}}
\newcommand{\lcomb}[2]{\left(\begin{smallmatrix}#1 \\ #2\end{smallmatrix}\right)}
\newcommand{\pk}{\Phi^{1/\ka}}
\newcommand{\pkkn}{\Phi^{1/\ka}_{k,n}}
\newcommand{\pkun}{\Phi^{1/\ka}_{1,n}}
\newcommand{\pzzsk}{\prod_{i,j}(z_i-z_j)^{\frac{m_im_j}2}}
\newcommand{\pttsk}{\prod_{i,j}(t_i-t_j)^2}
\newcommand{\ptzsk}{\prod_i\prod_l(t_i-z_l)^{-m_l}}
\newcommand{\pzz}{\prod_{i,j}(z_i-z_j)^{\frac{m_im_j}{2\ka}}}
\newcommand{\ptt}{\prod_{i,j}(t_i-t_j)^{\frac 2\ka}}
\newcommand{\ptz}{\prod_i\prod_l(t_i-z_l)^{\frac{-m_l}\ka}}
\newcommand{\slt}{\mathfrak{sl}_2}
\newcommand{\nb}{\nabla}
\newcommand{\ddx}[1]{\frac{\partial}{\partial #1}}
\newcommand{\dfx}[2]{\frac{\partial #1}{\partial #2}}
\newcommand{\ckn}{\Cc_{k,n}(z)}
\newcommand{\mdpd}[4]{\left(\begin{array}{rr}#1 & #2 \\ #3 & #4\end{array}\right)}

\newcommand{\Sing}{\operatorname{Sing}}
\newcommand{\Sym}{\operatorname{Sym}}
\newcommand{\im}{\operatorname{Im}}

\newcommand{\HH}{\mathbb{H}}
\newcommand{\GG}{\mathbb{G}}

\begin{document}
\begin{titlepage}
\vspace*{\stretch{1}}
\begin{center}
\Large\bfseries
Special functions, KZ type equations \\
and Representation theory \\
\vspace{2cm}
\normalsize\mdseries
Notes of a course given at MIT during the spring of 2002 by\\
\vspace{1cm}
\Large Alexander Varchenko \\
\normalsize
\vspace*{\stretch{1}}
\begin{minipage}{10cm}
Notes taken by: \\
Josh Scott: Lectures 1,2,3 \\
Mat\'\i as Gra\~na: Lectures 4,5,6 \\
Igor Mencattini: Lectures 7,8,9,10,11,12
\end{minipage}
\end{center}
\end{titlepage}

\vspace*{\stretch{1}}
This paper is a set of lecture notes of my course
``Special functions, KZ type equations, and representation theory''
given at MIT during the spring semester of 2002.
The notes were prepared by J.Scott (Lectures 1-3), M.Gra\~na
(Lectures 4-6), and I.Mencattini (Lectures 7-12). 

\bigskip
The notes do not contain new results, and are an exposition
(mostly without proofs) of various published 
results in this area, illustrated by the simplest
nontrivial examples. Some references are given at the end of each
lecture, but their list is not complete; the reader is referred
to the original articles for proofs and for complete references.  
\vspace*{\stretch{2}}

\setcounter{tocdepth}{1}
\tableofcontents

\markboth{A. Varchenko. Special functions, KZ type equations and Representation theory}{\textsc{Lecture \thechapter}}
\chapter{}
\section{Setting of Knizhnik--Zamolodchikov equation}
Given $m = (m_1,\dots,m_n)$ and $z = (z_1,\dots,z_n)$ fixed vectors in $\mathbb{C}^n$ and given
fixed $\kappa \in \mathbb{C}$ we define the following multi-valued function
and associated set of 1-forms:

\[ \Phi \  = \  \prod_{1 \leq i < j \leq n} 
          (z_i - z_j)^{ {m_i \frac{m_j}{2 \kappa} }} \
          \prod_{l=1}^{n} (t - z_l)^{\frac {-m_l}{\kappa} } \]  
    
\[ \eta_j \  = \ \Phi \  \frac{dt}{t - z_j} \ \ j=1, \dots,n \]

\bigskip
\noindent
The 1-forms $\eta_j$ are closed and, moreover, are cohomologically
dependent as they satisfy the relation

\[ m_1\eta_1 + \cdots + m_n\eta_n \ = \ -\kappa d\phi. \]

\bigskip
\noindent
Let $G = \pi_1(\mathbb{C}-\{ z_1,\dots,z_n \},*)$ and choose 
any contour $\gamma$ emanating from the base point * with
the property that $[\gamma]$ resides in the
commutator subgroup $[G,G]$. Since $[\gamma] \in [G,G]$
one can pull back the multi-valued 1-form $\eta_j$ to
a single valued 1-form on $S^1$ - after choosing
a branch of $\eta_j$ about the base point *.
In this manner the contour integral $\int_\gamma \eta_j$ becomes
well defined and its value will only depend on $\gamma$'s
homotopy class. Define $I^\gamma$ to be the following
holomorphic vector valued function:

\[ I^\gamma \ = \ (I_1, \dots, I_n) \ = \
                  \Big( \int_\gamma \eta_1, \dots, \int_\gamma \eta_n \Big) \]

\bigskip
\noindent
\begin{theor} $I^\gamma$ satisfies the following system of
            differential equations:

\[ \frac{\partial I}{\partial z_i} \ = \ 
   {\frac 1\kappa} \sum_{j \ne i} 
   \frac{\Omega_{ij}}{z_i - z_j} I \ \ \text{for} \ \ i=1, \dots, n
   \ \text{where} \]

\[ \Omega_{ij} \ = \ \begin{pmatrix} 

                   & \vdots^i &  & \vdots^j &  \\
                 
        {\scriptstyle i} \cdots & { (m_i -2)\frac{m_j}2 } & \cdots &
                 m_j & \cdots \\

                   & \vdots &  & \vdots &   \\

        {\scriptstyle j} \cdots & m_i & \cdots & \frac{m_i(m_j-2)}2 &
                 \cdots \\

                   & \vdots &  & \vdots &   

                   \end{pmatrix} \]

\bigskip
\noindent
All other diagonal entries are $\frac{m_im_j}2$ and the
remaining off-diagonal entries are all zero.

\end{theor}
\begin{proof}
For example, in the case of $i=1$ and the component
$\frac{\partial I^\gamma_2}{\partial z_1}$ of 
$\frac{\partial I^\gamma}{\partial z_1}$, we have

\[ \frac{\partial I^\gamma_2}{\partial z_1} = 
\int_\gamma \frac{\partial \eta_2}{\partial z_1} dt 
= \int_\gamma \Big( 
              \sum_{j \ne 1}^n \frac{m_i m_j}{2 \kappa}
              \frac{1}{z_1 - z_j} + \frac{m_1}\kappa
              \frac{1}{t-z_1} 
              \Big) \Phi \frac{dt}{t-z_2}. \]

\bigskip
\noindent
Using the identity $\frac{1}{t-z_1}\frac{1}{t - z_2}
= \frac{1}{z_1-z_2} \Big( \frac{1}{t-z_1} -
\frac{1}{t-z_2} \Big) $ we obtain the desired result

\[ \frac{\partial I^\gamma_2}{\partial z_1}  =
   \int_\gamma \Bigg( \frac{m_1(m_2-2)}{2 \kappa}
                      \frac{1}{z_1-z_2} \eta_2 +
                      \sum_{j \notin \{ 1,2 \} }
                      \frac{m_1m_j}{2 \kappa}
                      \frac{1}{z_1 - z_j} \eta_2 +
                      \frac{m_1}{\kappa} \frac{1}{z_1 -z_2} \eta_1
               \Bigg). \]

\end{proof}

\bigskip
\noindent
We will now describe a generalization of the above system
of differential equations which arises in the context
$\mathfrak{sl}_2(\mathbb{C})$-representations.

\bigskip
\noindent 
Recall that $ \mathfrak{sl}_2(\mathbb{C})$ is
the Lie algebra with generators 

\[ \Bigg< e = \begin{pmatrix} 0 & 1 \\ 0 & 0 \end{pmatrix}, \
       f = \begin{pmatrix} 0 & 0 \\ 1 & 0 \end{pmatrix}, \               
       h = \begin{pmatrix} -1 & 0 \\ 0 & -1 \end{pmatrix} 
 \Bigg> \]

\bigskip
\noindent
satisfying the familiar commutation relations 

\[ [h,e] = 2e \qquad \qquad  [e,f] = h \qquad \qquad [h,f] = -2f. \] 

\bigskip
\noindent
The element $\Omega =  e \otimes f + f \otimes e +
                       \frac{1}{2} h \otimes h  \in \mathfrak{g} 
                       \otimes \mathfrak{g} $ is called
the Casimir element and
for all $x \in \mathfrak{sl}_2(\mathbb{C})$ it satisfies the commutation 
relation
$[ x \otimes 1 + 1 \otimes , \Omega] = 0$ inside 
$\mathcal{U}(\mathfrak{sl}_2(\mathbb{C})) \otimes 
\mathcal{U}(\mathfrak{sl}_2(\mathbb{C}))$.

\bigskip
\noindent
Let $V_1, \dots, V_n$ be ${\mathfrak{sl}_2(\mathbb{C})}$-modules
and set $V = V_1 \otimes \cdots \otimes V_n$. An
element $x \in \mathfrak{sl}_2(\mathbb{C})$ acts on $V$
by $x\otimes 1 \otimes \cdots \otimes 1 \ + \ \dots \ + \ 1 
\otimes 1 \otimes \cdots \otimes x$. For indices $1 \leq i<j \leq n$ let
$\Omega^{ij}: V \longrightarrow V$ be the operator which
acts by $\Omega$ on the $i$th and $j$th positions and 
as the identity on all others. For example:

\[ \Omega^{12} (v_1 \otimes \cdots \otimes v_n)\  = \
   \Omega (v_1 \otimes v_2) \otimes v_3 \cdots \otimes v_n \]                

\bigskip
\noindent
\begin{defin} the KZ equation for a $V$-valued function 
$\phi(z_1, \dots, z_n)$ is 

\[ \frac{\partial \phi}{\partial z_i} =
   \frac{1}{\kappa} \sum_{j \ne i} 
   \frac{\Omega^{ij}}{z_i-z_j} \ \phi 
   \ \ \text{for} \ i=1 \dots n. \]

\bigskip
\noindent
This equation is defined on $U = \Big\{ \ (z_1,\dots,z_n) 
\in \mathbb{C}^n \ \Big| \ z_i \ne z_j \ \Big\}$.
\end{defin}       

\bigskip
\noindent
Mention connections with CFT and decorated Riemann surfaces.

\bigskip
\noindent
A solution $\phi$ of the KZ equation one can be
translated by an element $x \in \mathfrak{sl}_2(\mathbb{C})$
to obtain a new $V$-valued function $x \cdot \phi$. It
follows from the commutation relation
$[ \Omega, x \otimes 1 + 1 \otimes x ] =0$ that 
$x \cdot \phi$ is also a solution of the KZ equation.
From this we can conclude that

\begin{corol} For any $x \in \mathfrak{sl}_2(\mathbb{C})$ the
KZ equation preserves the eigenspaces of the action
of $x$ on $V$. In other words, if $\phi$ is a KZ solution
such that $x \cdot \phi(z_0) = \lambda \phi( z_0)$
for some $z_0 \in \mathbb{C}^n$ then $x \cdot
\phi(z) = \lambda \phi(z)$ for all $z \in
\mathbb{C}^n$.

\end{corol}   

\bigskip
\noindent
For $m \in \mathbb{C}$, the Verma module $M_m$ is the
infinite dimensional $\mathfrak{sl}_2(\mathbb{C})$-representation
with highest weight $m$ generated by a single vector
$v_m$ where $hv_m = m v_m$ and $ev_m =0$. The 
vectors $f^kv_m$ for $k=0,1,\dots$ form a basis for $M_m$.
The generators $e$, $f$, and $h$ act on this basis as
indicated below:

\[ f \cdot f^kv_m \ = \ f^{k+1}v_m \qquad
   h \cdot f^kv_m \ = \ (m-2k)f^kv_m \qquad
   e \cdot f^kv_m \ = \ k(m-k+1)f^{k-1}v_m \]

\bigskip
\noindent
If $m \notin \mathbb{Z}_{\geq 0}$ then $M_m$ is irreducible,
otherwise $ <f^{m+1}v_m> =
\mathfrak{sl}_2(\mathbb{C}) \cdot f^{m+1}v_m$ will
be an invariant subspace. Let $L_m = M_m / <f^{m+1}v_m>$. 
This quotient is irreducible with basis 
$v, fv, \dots, f^mv$.

\bigskip
\noindent
It is a familiar exercise in theory of $\mathfrak{sl}_2(\mathbb{C})$-
representations to check the following tensor decomposition:

\[ L_m \otimes L_l \ = \ L_{m+l} \oplus L_{m+l-2} \oplus 
   \cdots \oplus L_{m-l} \] 

\bigskip
\noindent
Let $m = (m_1,\dots,m_n) \in \mathbb{C}^n$ and set
$|m| = m_1 + \cdots + m_n$. Let $M^{\otimes m}
= M_{m_1} \otimes \cdots \otimes M_{m_n}$. For
$J=(j_1,\dots,j_n) \in \mathbb{Z}_{\geq 0}$ let
$f_Jv = f^{j_1}v_{m_1}\otimes \cdots \otimes f^{j_n}v_{m_n}$
where $v_{m_i}$ is the primitive vector in $M_{m_i}$.
The vectors $f_Jv$ form a basis for the module
$M^{\otimes m}$.

\bigskip
\noindent
Now $h \cdot f_Jv = \Big( |m|-2|J| \Big) f_Jv$. Define
$M^{\otimes m}[\lambda] = \Big\{ \ v \in M^{\otimes m} \ 
\Big| \ h\cdot v = \lambda v \ \Big\}$ and $\Sing M^{\otimes m}
[\lambda] = \Big\{ \ v \in M^{\otimes m}[\lambda] \ \Big| \ 
e \cdot v = 0 \ \Big\}$. The singular eigenspaces
$\Sing M^{\otimes m}[\lambda]$ generate the entire
module $M^{\otimes m}$ and, in view of Corollary 1, 
it follows that it is enough to produce KZ solutions
for the singular eigenspaces only.

\section{Solutions}
Case 1: $ \Sing M^{\otimes m}[|m|]$

\bigskip
\noindent
$\Sing M^{\otimes m}[|m|]$ is spanned by the single
vector $v_{m_1} \otimes \cdots \otimes v_{m_n}$. If $\phi$ is
a KZ solution residing in the $|m|$-eigenspace
then it is of the form $\phi = I_0 (z_1, \dots ,z_n) 
v_{m_1} \otimes \cdots \otimes v_{m_n}$ where
$I_0$ is a scalar valued function. In an effort to
compute the right hand side of the KZ equation we
obtain:

\[ \Omega^{ij} \Big( v_{m_1}\otimes \cdots \otimes v_{m_n} \Big) \ = \
   \frac{m_i m_j}{2} \  v_{m_1}\otimes \cdots \otimes v_{m_n} \]

\bigskip
\noindent
In this case the KZ equation reads as:

\[ \frac{\partial I_0}{\partial z_i} \ = \ 
   \sum_{j \ne i} \frac{m_i m_j}{2 \kappa}
   \frac{I_0}{z_i - z_j} \ \ \text{for} \ i=1 \dots n \]

\noindent
with solution given by

\[ I_0 \ = \ \prod_{1 \leq i<j \leq n} (z_i-z_j)^{\frac{m_im_j}{2 \kappa}}  \]

\bigskip
\noindent
In view of Corollary 1 we can conclude that
$f^{k} \cdot I_0 v_{m_1} \otimes \cdots \otimes v_{m_n} =
\newline I_0 f^{k} \cdot v_{m_1} \otimes \cdots \otimes v_{m_n}$ is
also a KZ solution for any $k \in \mathbb{Z}_{\geq 0}$.

\bigskip
\noindent
Case 2: $\Sing M^{\otimes m}[|m| -2]$

\bigskip
\noindent
This eigenspace is 
\[ \Sing M^{\otimes m}[|m| -2] =
\Big\{ \ w = \sum_{l =1}^{m} I_l v_{m_1} 
             \otimes \cdots \otimes f \cdot v_{m_l} \otimes 
             \cdots \otimes v_{m_n} \ \Big| \ 
       e \cdot w = 0 \ \Big\} \]

\bigskip
\noindent
where the $I_l$'s are scalars. In this case the condition $e \cdot w = 0$ 
can be reformulated as $m_1I_1 + \cdots + m_n I_n =0$.
In order to compute the right hand side of the KZ equation
we compute the action of $\Omega^{ij}$ on the
basis vectors $v_{m_1} \otimes \cdots \otimes f \cdot v_{m_l}
\otimes \cdots \otimes v_{m_n}$. Assume $i<j$. If the basis
vector is of the form $ {\bf v} = \cdots \otimes v_{m_i} \otimes
\cdots \otimes v_{m_j} \otimes \cdots$ then
$\Omega^{ij}({\bf v}) = \frac{m_im_j}{2} {\bf v}$. If
${\bf v}$ is of the form $ \cdots \otimes f \cdot
v_{m_i} \otimes \cdots \otimes v_{m_j} \otimes \cdots $ then

\[ \Omega^{ij}{\bf v} \ = \ m_i \Big( \cdots \otimes v_{m_i}
                            \otimes \cdots \otimes f \cdot
                            v_{m_j} \otimes \cdots \Big) \ + \
                            \frac{(m_i-2)m_j}{2} \Big( \cdots
                            \otimes f \cdot v_{m_i} \otimes
                            \cdots \otimes v_{m_j} \otimes
                            \cdots \Big). \]

\bigskip
\noindent
If the basis vector ${\bf v}$ is of the form
$ \cdots \otimes v_{m_i} \otimes \cdots \otimes f \cdot
v_{m_j} \otimes \cdots $ then 

\[ \Omega^{ij}{\bf v} \ = \ m_j \Big( \cdots \otimes f \cdot v_{m_i}
                            \otimes \cdots \otimes 
                            v_{m_j} \otimes \cdots \Big) \ + \
                            \frac{(m_j-2)m_i}{2} \Big( \cdots
                            \otimes v_{m_i} \otimes
                            \cdots \otimes f \cdot v_{m_j} \otimes
                            \cdots \Big). \]

\bigskip
\noindent
From these computations we may surmise that the matrix
representing $\Omega^{ij}$ is of the form

\[ \Omega_{ij} \ = \ \begin{pmatrix} 

                   & \vdots^i &  & \vdots^j &  \\
                 
        {\scriptstyle i} \cdots & \frac{ (m_i -2)m_j}{2 } & \cdots &
                 m_j & \cdots \\

                   & \vdots &  & \vdots &   \\

        {\scriptstyle j} \cdots & m_i & \cdots & \frac{m_i(m_j-2)}{2} &
                 \cdots \\

                   & \vdots &  & \vdots &   

                   \end{pmatrix} \]

\bigskip
\noindent
where all other diagonal entries are $\frac{ m_im_j}{2}$ and the
remaining off-diagonal entries are all zero. Of course
this is the same matrix we encountered in the beginning
of the lecture. We know explicit solutions namely,

\[ \phi^\gamma (z_1, \dots, z_n) \ = \
   \sum_{j=1}^{n} \int_\gamma \Phi
   \frac{dt}{{t-z_j}} v_{m_1} \otimes \cdots \otimes
   f \cdot v_{m_j} \otimes \cdots \otimes v_{m_n} 
   \ \ \text{where} \]

\[ \Phi \ = \ \prod_{ i<j } \ (z_i - z_j)^\frac{m_i m_j}{{2 \kappa}} \
              \prod_{l=1}^{n} \ (t-z_l)^\frac{- m_l}{\kappa}. \]

\bigskip
\noindent
Applying our translation argument again we may conclude
that $f^k \cdot \phi^\gamma$ is a KZ solution with values
in $M^{\otimes m}\big[ |m| - 2 - 2k \big]$.

\bigskip
\noindent
Case 3: $\Sing M^{\otimes m} \big[ |m| - 2k \big]$ where
$k \in \mathbb{Z}_{\geq 0}$

\bigskip
\noindent
Define $\Phi_{k,n}(t,z,m) = 
\Phi_{k,n} (t_1, \dots,t_k,z_1, \dots,z_n,m_1, \dots,m_n)$ to be
the following expression:

\[ \prod_{i<j} \  (z_i-z_j)^{\frac{m_im_j}{2}} \ 
   \prod_{1 \leq i \leq j \leq k}  (t_i-t_j)^2 \
   \prod_{l=1}^{n} \ \prod_{i=1}^{k} \ (t_i-z_l)^{-m_l}. \]

\bigskip
\noindent
Comments about discriminantal arrangement of hyperplanes
associated to the singularities $t_i=z_l$ and
vanishings $t_i=t_j$. 

\bigskip
\noindent
Recall that $M^{\otimes m} \big[ |m|-2k \big]$ has a 
basis of vectors $f_J v$ where $J=(j_1,\dots,j_n)$ and $|J|=k$. Associate
to each basis vector $f_Jv$ the following rational
function:

\[ A_J(t,z) \ = \ \frac{1}{j_1! \cdots j_n!} \
                            \Sym_t \Bigg[ \ \prod_{l=1}^{n} 
                            \prod_{i=1}^{j_l}
                            \ \frac{1}{t_{ j_1 + \cdots + j_{l-1}+i - z_l}} 
                            \Bigg] \ \ \text{where}\]

\[ \Sym_t f(t_1,\dots,t_k) := \sum_{\sigma \in S_k}
f(t_{\sigma(1)}, \dots ,t_{\sigma(k)}). \]

\bigskip
\noindent
For example: $A(1,0,\dots,0) = \frac{1}{t_1 - z_1}$ and
$A(2,0,\dots,0) = \frac{1}{t_1 - z_1} \frac{1}{t_2 - z_1}$ and
\newline \indent \indent \indent \indent \indent  
$A(1,1,0,\dots,0) = \frac{1}{t_1-z_1} \frac{1}{t_2-z_2}
+ \frac{1}{t_2-z_1} \frac{1}{t_1-z_2} $.

\bigskip
\noindent
We define the following $M^{\otimes m}[ |m| -2k]$-valued
functions:

\[ \phi^\gamma(z_1,\dots,z_n) \ = \ 
   \sum_{|J|=k} \Bigg( \int_\gamma \Phi_{k,n}(t,z, m))^
   \frac{1}{\kappa} A_J(t,z) \ dt_1 \wedge \cdots
   \wedge dt_k \Bigg) {\bf \cdot} f_Jv.  \]

\bigskip
\noindent
In the next lecture we shall prove the following theorem:

\begin{theor}

$\phi^\gamma (z) \in \Sing M^{\otimes m} \Big[ |m|-2k \Big]$
and it is a KZ solution.

\end{theor}

\bigskip
\noindent
\begin{defin}

The functions $I_J^\gamma := \int_\gamma \Phi_{k,n}^{\frac{1}{\kappa}}
A_J \ dt_1 \wedge \cdots \wedge dt_k $ are called 
the hypergeometric functions associated with
$\Sing M^{\otimes m} \Big[ |m| -2k \Big]$.

\end{defin}

\begin{claim}
Let $J + {\bf 1}_l = (j_1, \dots, j_l +1, \dots, j_n)$.
If $|J|=k-1$ then 

\[ \sum_{l=1}^{n}I_{J+ {\bf 1}_l}(z) (j_l+1)(m_l-j_l) \ = \ 0. \]

\bigskip
\noindent
This statement is equivalent to the claim
that $\phi^\gamma(z) \in \Sing M^{\otimes m} \Big[ |m| -2k \Big]$.
\end{claim}

\begin{reftc}
\cite{V2}, \cite{SV}.
\end{reftc}


\chapter{}
We begin by sketching the proof of Theorem 2 from
Lecture 1. To do so, we modify the expressions for $\phi^\gamma$
slightly. For any differential form $f(t_1,\dots,t_k)$ let

\[ \operatorname{Ant}_t f(t_1,\dots,t_k) \ = \
   \sum_{\sigma \in S_k} (-1)^{l(\sigma)} f(t_{\sigma(1)}, 
   \dots, t_{\sigma(k)}) \]

\bigskip
\noindent
where $l(\sigma)$ is the length of the permutation $\sigma$. For
$J=(j_1,\dots,j_n) \in \mathbb{Z}_{\geq 0}^n$ with $|J| = j_1 + 
\cdots j_n = k$ set

\[ \eta_J \ = \ \frac{1}{j_1 ! \cdots j_n !} 
              \operatorname{Ant}_t \Bigg[
              \bigwedge_{l=1}^{n} \bigwedge_{i=1}^{j_l}
              \frac{d(t_{j_1 + \cdots + j_{l-1} + i} - z_l)}{t_{ j_1 + \cdots + j_{l-1} + i} -z_l}
              \Bigg]. \]

\bigskip
\noindent
When all $z_i$ specialized to a fixed constant then $\eta_J$
is $A_J dt_1 \wedge \cdots \wedge dt_k$ as defined
in Lecture 1. For $\kappa \in \mathbb{C}^*$ let

\[ \phi^\gamma \ = \ \sum_{|J|=k} \int_{\gamma} 
                     \Phi_{k,n}^{\frac{1}{\kappa}} \eta_J f_Jv. \]

\begin{claim} The fact that $\phi^\gamma(z) \in
              \Sing M^{\otimes m} \Big[ |m|-2k \Big]$ 
              follows from the identities

\[ \Bigg[ \kappa \ d \Big( \Phi_{k,n}^{\frac{1}{\kappa}} \eta_J \Big)
   \ + \ \sum_{l=1}^n \Phi_{k,n}^{\frac{1}{\kappa}} \eta_{J + {\bf 1}_l}
   (j_l +1)(m_l-j_l) \Bigg] \Bigg|_{ \text{all} \ z_i= \ \text{constant} }
   \ = \ 0   \]

\bigskip
\noindent
where $J=(j_1,\dots,j_n)$ is any multi-index with $|J|=k-1$. 
\end{claim}

\bigskip
\noindent
The validity of this identity is easily checked
by expanding $d (\Phi_{k,n}^{\frac{1}{\kappa}} \eta_J )$. 
To see that this identity implies $\phi^\gamma \in
\Sing M^{\otimes m} \Big[ |m|-2k \Big]$ integrate
the above expression over $\gamma$. Since $\gamma$ is
a cycle the integral $\int_\gamma \kappa \ d( \Phi_{k,n}^{\frac{1}{\kappa}}
\eta_J)$ vanishes. The remaining part of the integral is 
exactly the coefficient of $f_Jv$ in $e \cdot \phi^\gamma$. 

\begin{claim}
The fact that $\phi^\gamma$ is a KZ solution follows from the
identity 

\[ d \Bigg( \sum_{|J|=k} \Phi_{k,n}^{\frac{1}{\gamma}} \eta_J f_Jv \Bigg) 
   \  = \ \frac{1}{\kappa} \sum_{i<j} \Omega^{ij}
   \Bigg( \frac{d(z_i-z_j)}{z_i-z_j} \wedge 
   \sum_{|J|=k} \Phi_{k,n}^{\frac{1}{\kappa}} \eta_J f_Jv \Bigg). \]

\end{claim}

\bigskip
\noindent
The identity is easily checked by expanding the
left hand side with $d= d_z + d_t$.   
Integrate the above expression over a cycle $\gamma$ in
the coordinate $t$. From the identity we obtain

\[   \int_\gamma d_t \Big[ \sum_{|J|=k} \Phi_{k,n}^{\frac{1}{\kappa}}
            \eta_J f_J v \Big] 
     \ = \ -\int_\gamma d_z \Big[ \sum_{|J|=k} 
            \Phi_{k,n}^{\frac{1}{\kappa}} \eta_J f_J v \Big] \ + \
           \int_\gamma \ \frac{1}{\kappa} \sum_{i<j} \Omega^{ij}
   \Bigg( \frac{d(z_i-z_j)}{z_i-z_j} \wedge 
   \sum_{|J|=k} \Phi_{k,n}^{\frac{1}{\kappa}} \eta_J f_Jv \Bigg). \]

\bigskip
\noindent
By Stokes' Theorem we know
that $\int_\gamma d_t \Big[ \sum_{|J|=k} \Phi_{k,n}^{\frac{1}{\kappa}}
\eta_J f_J v \Big]$ vanishes so we may conclude that

\[ \int_\gamma d_z \Big[ \sum_{|J|=k} 
   \Phi_{k,n}^{\frac{1}{\kappa}} \eta_J f_J v \Big] \ = \
   \int_\gamma \ \frac{1}{\kappa} \sum_{i<j} \Omega^{ij}
   \Bigg( \frac{d(z_i-z_j)}{z_i-z_j} \wedge 
   \sum_{|J|=k} \Phi_{k,n}^{\frac{1}{\kappa}} \eta_J f_Jv \Bigg) 
   \ \text{or} \]

\[ d_z \Bigg[ \int_\gamma \sum_{|J|=k} 
   \Phi_{k,n}^{\frac{1}{\kappa}} \eta_J f_J v \Bigg] \ = \
   \frac{1}{\kappa} \sum_{i<j} \Omega^{ij} \Bigg(
   \frac{d(z_i-z_j)}{z_i-z_j} \Bigg) \cdot \int_\gamma 
   \sum_{|J|=k} \Phi_{k,n}^{\frac{1}{\kappa}} \eta_J f_Jv. \]

\bigskip
\noindent
By isolating the $\frac{\partial}{\partial z_j}$ component
in the right hand side of
the above expression we obtain the KZ equation as required.

\section{Hyperplane Arrangements}

The modern theory of arrangements has some origins in Hilbert's 13th problem. 
An algebraic function $z = z(x_1,\dots,x_k)$ is
a multi-valued function defined by an
equation of the form 

\[ z^n + P_1(x_1,\dots,x_k)z^{k-1} + \cdots + P_n(x_1,\dots,x_k) = 0 \]

\bigskip
\noindent
where the $P_i$'s are polynomials. 

\bigskip
\noindent
{\bf Hilbert's 13th Problem:}
Show that an algebraic function $z$ defined by

\[ z^7 + az^2 + bz + c = 0 \]

\bigskip
\noindent
can not be represented as a composition of continuous functions
in two variables.

\bigskip
\noindent
Kolmogorov and Arnold showed that in fact one
can find such a two variable decomposition. 
Despite the Kolmogorov/Arnold result it is
believed that the composition is impossible
if one restricts to algebraic functions.
Arnold wanted to understand why algebraic
functions of many variables become more
complicated. The idea was to invent invariants
of algebraic functions which would detect
when compositions are possible. The simplest
characteristic of an algebraic function $z(x_1,\dots,x_k)$ 
is its discriminant

\[ \Delta_P \ = \ \Big\{ \ {\bf x} \in \mathbb{C}^k \ \Big| \ 
                z^n + P_1({\bf x}) + \cdots
                P_n({\bf x}) = 0 \ \text{has multiple roots} \ 
                \Big\}. \]

\bigskip
\noindent
Any algebraic function 
is induced from the universal algebraic
function $z^n + a_1z^{n-1} + \cdots + a_n =0$ by the map
$P:\mathbb{C}^k \longrightarrow \mathbb{C}^n$ given 
by ${\bf x}=(x_1,\dots,x_k) \mapsto \big(P_1({\bf x}),\dots, P_n({\bf x})
\big)$. The discriminant can then be expressed
as $\Delta_P = P^{-1}(\Delta)$ where 

\[ \Delta \ = \ \Big\{ \ (a_1,\dots,a_n) \in \mathbb{C}^n \ \Big| \
   z^n + a_1z^{n-1} + \cdots + a_n = 0 \ \text{has multiple roots} 
   \ \Big\}. \]

\bigskip
\noindent
By passing to cohomology we obtain the map
 $P^*:H^*(\mathbb{C}^n - \Delta) \longrightarrow H^*(\mathbb{C}^k - \Delta_P)$.
For any $\alpha \in H^*(\mathbb{C}^n - \Delta)$ one has a
characteristic class $P^*\alpha$ in $H^*(\mathbb{C}^k - \Delta_P)$
for the algebraic function $z(x_1,\dots,x_k)$.
One might hope that the cohomology of the compliment of 
the discriminant might give obstructions to
the representability of an algebraic function as a composition.
As a first step one might try to describe $H^*(\mathbb{C}^n - \Delta)$.
There is, however, a simplification.
The Vieta map $V:\mathbb{C}^n \longrightarrow \mathbb{C}^n$ is defined
by $z=(z_1,\dots,z_n) \mapsto \big(a_1(z),\dots,a_n(z)
\big)$ where $a_k$ is $(-1)^k$ times the $k$th elementary
symmetric function in the variables $z_1,\dots,z_n$.
Under the Vieta map the set $D = \{ \ (z_1,\dots,z_n) \ | \
z_i = z_j \ \exists \ i,j \ \}$ maps to $\Delta$.
The problem of computing $H^*(\mathbb{C}^n - \Delta)$ thus reduces
to the study of $H^*(\mathbb{C}^n - D)$ since
$H^*(\mathbb{C}^n-D)$ is the symmetric part of the $S_n$ action
on $H^*(\mathbb{C}^n - \Delta)$.

\bigskip
\noindent
Arnold described $H^*(\mathbb{C}^n - D)$.
Consider the 1-forms $\omega_{ij} = \frac{1}{2\pi i}
\frac{d(z_i-z_j)}{z_i-z_j}$. Clearly $\omega_{ij}
= \omega_{ji}$. An easy exercise of Arnold reveals that
for distinct $i$, $j$, and $k$ 

\[ \omega_{ij} \wedge \omega_{jk} \ + \ \omega_{jk}\wedge \omega_{ki}
   \ + \ \omega_{ki} \wedge \omega_{ij} \ = \ 0. \]

\begin{theor}[Arnold] Let $\mathcal{A} = \mathbb{C}[\omega_{ij}]$ be 
the exterior algebra
generated by the differential forms $\omega_{ij}$.
The map $\alpha \mapsto [ \alpha ]$ of $\mathcal{A}$
to $H^*(\mathbb{C}^n - D)$ is an isomorphism.
\end{theor}

\begin{theor}[Arnold] Let $P_\Delta(t) = \sum_{k=0}^n \dim H^k
(\mathbb{C}^n - D) t^k$ be the Poincar\'e polynomial, then

\[ P_\Delta (t) \ = \ \Big( 1+t \Big) \Big( 1+2t \Big) \cdots 
   \Big( 1 + (n-1)t \Big). \]
\end{theor}

\bigskip
\noindent
Brieskorn proved Theorem 1 for any hyperplane arrangement and
Orlik and Solomon combinatorially described the algebra
of logarithmic differential forms.
Realize each hyperplane $H$ in the arrangement $\mathcal{H}$ as the
vanishing of an equation $f_H$. Associate to $f_H$ the
closed 1-form $\omega_H = \frac{1}{2 \pi i} \frac{df_H}{f_H}$;
note that $\omega_H$ does not depend on the choice of $f_H$.
Consider the exterior algebra $\mathcal{A} = \mathbb{C}[\omega_H]_{H \in 
\mathcal{H}}$. The relations in $\mathcal{A}$
are given as follows. A collection of hyperplanes
$H_1,\dots,H_j$ is said to intersect transversally if
$\operatorname{codim}\big( H_1 \cap \cdots \cap H_j \big) = j.$
To obtain the relations for $\mathcal{A}$, take
any collection of hyperplanes $H_1,\dots,H_{j+1}$ such
that $H_1 \cap \cdots \cap H_{j+1} \ne \phi$ and
such that $H_1,\dots,H_{j+1}$ do {\bf not} intersect transversally.
Associate to this collection the relation

\[ \sum_{i=1}^{j+1} \ (-1)^j \ \omega_{H_1}\wedge \cdots
   \wedge \hat{\omega_{H_i}} \wedge \cdots 
   \wedge \omega_{H_{j+1}} \ = \ 0. \]

\begin{theor}[Orlik \& Solomon]
These are the defining relations for $\mathcal{A} =
\mathbb{C}[\omega_H]$.
\end{theor}

\bigskip
\noindent
Theorem 2 is generalized in the following manner. Let
$G$ be a Coxeter group and let $\mathcal{H}$ be
the collection of reflecting hyperplanes for $G$. Set
$\mathcal{H}_\mathbb{C}$ to be $\mathcal{H}$'s complexification.

\begin{theor}[Terao] Let $P(t)$ be the
Poincar\'e polynomial of $\mathbb{C}^k - \bigcup_{H \in \mathcal{H}}H$
then,

\[ P(t) \ = \ \prod_j \Big( 1 + (d_j-1)t \Big) \] 

\bigskip
\noindent
where the $d_j$'s are the degrees of the basis set
of homogeneous polynomials invariant under the action
of $G$.
\end{theor}

\section{Classical Hypergeometric Series}

\bigskip
\noindent
For a complex number $z$ define the
$n$th rising power $(z)_n = z(z+1)\cdots(z+n-1)$ with
the provision that $(z)_0 =1$.
The classical hypergeometric series is defined as:

\[ F(a,b,c;z) \ = \ \sum_{n=0}^{\infty} \frac{(a)_n(b)_n}{(c)_n}
                    \ \frac{z^n}{n!} \]

\bigskip
\noindent
The rising power can be expressed via the Gamma function as
$(z)_n = \frac{ \Gamma(a+n)}{\Gamma(a)}$ where 

\[ \Gamma(z) \ = \ \int_0^\infty e^{-t}t^{z-1} dt. \]

\bigskip
\noindent
Note that the Gamma function is defined and holomorphic for $z \in
\mathbb{C}-\mathbb{Z}_{<0}$. 

\bigskip
\noindent
\section{Properties of the Hypergeometric series (Gauss 1812)}

\bigskip
\indent 1. The series $F(a,b,c;z)$ is convergent for $|z|<1$
and it converges for $|z|=1$ \indent \indent \ provided $\Re(c-a-b)>0.$

\bigskip

\indent 2. Contiguous Recurrence Relations: 

\[ (c-a-b)F \ + \ a(1-z)F(a+) \ - \ (c-b)F(b-) \ = \ 0 \]

\bigskip
\indent \indent where $F(a \pm)=F(a \pm 1,b,c;z)$, 
$F(b \pm 1) = F(a,b \pm 1,c;z)$, and
$F(c \pm 1)=$ \indent \indent $F(a,b,c \pm 1;z)$ are the
associated continuous functions. Gauss showed there \indent
\indent is a linear relation between $F$ and any two continuous 
functions with \newline \indent
\indent coefficients which are linear in $z$.

\bigskip

\indent 3. $F(a,b,c;1) \ = \ \frac{\Gamma(c)\Gamma(c-a-b)}{\Gamma(c-a)\Gamma(c-b)}.$

\begin{theor}[Euler] $F$ satisfies the Gauss Hypergeometric Equation

\[ z(1-z) \frac{d^2u}{dz^2} \ + \ \Big( c- (a+b+1)z \Big)\frac{du}{dz} \ - \ abu \ = \ 0. \]

\end{theor} 

\begin{theor}[Euler] 

\[ F(a,b,c,z) \ = \ \frac{\Gamma(c)}{ \Gamma(b) \Gamma(c-b)}
\int_{0}^{1} t^{b-1}(1-t)^{c-b-1} (1-zt)^{-a} dt. \]

\end{theor}

\bigskip
\noindent
If the specialization $z=1$ is made in the above expression
formula \#3 is recovered.

\bigskip
\noindent
\[ \text{Set} \ J \ = \ \int_a^b (s-a)^\alpha (b-s)^\beta (z-s)^\gamma ds. \]

\bigskip
\noindent
If we make the change of variables $t=\frac{s-a}{b-a}$ in the integral $J$ we
obtain

\[ J \ = \ \frac{\Gamma(\alpha + 1)\Gamma( \beta +1)}{\Gamma( \alpha + \beta + 2)} (b-a)^{\alpha + \beta +1} (z-a)^\gamma
F\Bigg(-\gamma,\alpha +1, \alpha + \beta + 2, \frac{b-a}{z-a}\Bigg). \]

\begin{reftc}
\cite{V2}, \cite{SV}, \cite{R}, \cite{OS}
\end{reftc}


\chapter{}
In this lecture we address which cycles $\gamma$ can
be chosen in order to make the expressions

\[ \phi^\gamma(z) \ = \ \sum_{|J|=k} \
\int_{\gamma} \Phi_{k,n}^{\frac{1}{\kappa}}(t,z,m)
    A_J(t,z) dt_1 \wedge \cdots \wedge dt_k \]

\bigskip
\noindent
single valued non trivial KZ solutions, where $\kappa \in \mathbb{C}^*$ and $m=(m_1,\dots,m_n) \in \mathbb{C}^n$ are fixed and 

\[ \Phi_{k,n}(t,z,m)
   \ = \ \prod_{1 \leq i<j \leq n} (z_i-z_j)^{\frac{m_i m_j}{2}}
                    \prod_{1 \leq i \leq j \leq k}
                    (t_i-t_j)^2 \
                    \prod_{l=1}^n \prod_{i=1}^k
                    \ (t_i-z_l)^{-m_l}   \]

\[ A_J(t,z) \ = \ \frac{1}{j_1! \cdots j_n!} \
                            \Sym_t \Bigg[ \ \prod_{l=1}^{n} 
                            \prod_{i=1}^{j_l}
                            \ \frac{1}{t_{ j_1 + \cdots + j_{l-1}+i - z_l}} 
                            \Bigg]. \] 

\bigskip
\noindent
When $k=1$ the cycle $\gamma$ should be
chosen so that its homotopy class $[\gamma]$
resides in the commutator subgroup of
$\pi_1\Big(\mathbb{C}-\{z_1,\dots,z_n\}\Big)$. For example,
a Pochhammer cycle about any points $z_i$ and $z_j$ will
serve for $\gamma$. In general the appropriate $\gamma$ will
reside in homology groups twisted by coefficients
associated with $\Phi_{k,n}^{\frac{1}{\kappa}}$.
We start by defining and computing this homology for $k=1$.

\begin{defin}
A twisted $k$-cell is a pair $(\Delta^k,s)$ where
$\Delta^k \subset \mathbb{C}- \{z_1,\dots,z_n\}$ is a 
singular $k$-cell and $s$ is a univalent
branch of $\Phi_{1,n}^{\frac{1}{\kappa}}$ over $\Delta^k$.
A twisted $k$-chain is a formal linear (over $\mathbb{C}$)
combination of $k$-cells.

\bigskip
\noindent
Define the boundary operator $d_k$ on any
$k$-cell by the formula $d_k(\Delta^k,s) = (\partial
\Delta^k,s|_{\Delta^k})$. Extend it by linearity
to all $k$-chains.

\bigskip
\noindent
The $k$th $\Phi_{1,n}^{\frac{1}{\kappa}}$-twisted
homology is defined as 

\[ H_k\Big( \mathbb{C}-\{z_1,\dots,z_n\},
\Phi_{1,n}^{\frac{1}{\kappa}} \Big) \ := \
 \frac{ \ker d_k}{ \im d_{k+1}  }. \]

\end{defin}

\bigskip
\noindent
Note that $\phi^\gamma$ is a well-defined
non-trivial KZ solution provided $\gamma \in
H_1$. 

\bigskip
\noindent
For what follows let $q = e^{\frac{2\pi i}{\kappa}}$ and
for any $a \in \mathbb{C}$ set

\[ [a] \ = \ \frac{q^{\frac{a}{2}} - q^{\frac{-a}{2}}}{ q^{\frac{1}{2}} - q^{\frac{-1}{2}} } \]

\bigskip
\noindent
where $q^a = e^{\frac{2a \pi i}{\kappa}}$. Note that if $q \rightarrow
1$ (or $\kappa \rightarrow \infty$) then $[a] \rightarrow a$. We
will begin by computing the twisted homology in a toy example;
namely $H_*(\mathbb{C}-\{0\},t^{\frac{-m}{\kappa}})$. Take $p \in
\mathbb{C}-\{0\}$ and choose a loop $l_1$ about $0$ passing
through $p$. Fix a value $s_0$ of $t^{\frac{-m}{\kappa}}$ at $p$ 
and a branch $s$ of $t^{\frac{-m}{\kappa}}$ over $l$ whose value
at $l$'s endpoint is $s_0$. Computing the boundary map $d_1$ we
obtain:

\[ d_1(l,s) \ = \ ( p,s|_{\text{end point}} ) -
                  ( p,s|_{\text{start point}} ) \ = \
                  ( p,s_0) - (p,q^{m} s_0)  \] 

\[ = \ (1-q^{m})(p,s_0) \]

\bigskip
\begin{lemma}
The chain complex $0 \longrightarrow \mathbb{C}\cdot(l,s) 
\overset{d_1}{\longrightarrow}
\mathbb{C}\cdot(p,s_0) \longrightarrow 0$ computes the homology groups
$H_*(\mathbb{C}-\{0\},t^{\frac{-m}{\kappa}})$.

\end{lemma}

\begin{corol}
$H_0$ and $H_1$ are both $0$ if $[m] \ne 0$ and
they are both 1 dimensional if $[m]=0$.
\end{corol}

\bigskip
\noindent
Consider $H_*\Big( \mathbb{C}-\{z_1,\dots,z_n\},\Phi_{1,n}^{\frac{1}{\kappa}}
\Big)$ and assume that $z_1<\cdots<z_n$ are all real.
Fix a point $p_0$ in the upper-half plane and choose
loops $l_i$ emanating from $p_0$ and circling only around 
$z_i$ for $i=1,\dots,n$. Choose branches $s_0,s_1,\dots,s_n$
of $\Phi_{1,n}^{\frac{1}{\kappa}}$ over $p_0,l_1,\dots,l_n$
respectively. Note that if we analytically continue
$s_i$ along the loop $l_j$ (for $1 \leq,i,j \leq n$) then
$s_i$ at $l_j$'s endpoint is $q^{-m_j}$ times its
value at the start point. 

\begin{lemma}
The complex $ 0 \longrightarrow \bigoplus_{j=1}^n \mathbb{C}\cdot (l_j,s_j)
\overset{d_1}{\longrightarrow}
\mathbb{C}\cdot (p_0,s_0) \longrightarrow 0$
computes $H_*\Big( \mathbb{C}- \{z_1,\dots,z_n\}, \Phi_{k,n}^{\frac{1}{\kappa}} \Big)$. 
\end{lemma}
 
\begin{lemma}
One can choose $c_j \in \mathbb{C}$ for $j=0,\dots,n$ such
that 

\[ d_1(\omega_j) \ = \ [m_j] \Big( q^{\frac{ -m_1 - \cdots - m_{j-1} -m_{j+1} + \cdots m_n}{4}}
\Big) \omega_0  \ \ j=1, \dots,n \]

\bigskip
\noindent
where $\omega_o = c_0(p_0,s_0), \cdots, \omega_n = c_n(l_n,s_n)$.
\end{lemma}

\begin{corol}
\[ H_1 = \Big\{ \ I_1 \omega_1 + \cdots + I_n \omega_n \ \Big|
\ \sum_{j=1}^{n} \ [m_j] \ q^{{\frac{1}{4}}\big( -\sum_{i<j}m_i + 
\sum_{i>j}m_i \big)} I_j = 0 \ \Big\}. \]
\end{corol}

\bigskip
\noindent
Note that as $q \rightarrow 1$ the formula above becomes
$\sum_{j=1}^n m_jI_j = 0$ which is precisely the singular
vector condition; i.e. the condition to be in 
$\Sing M^{\otimes m} \Big[ |m| -2 \Big]$.

\begin{corol}
If $[m_1] = \cdots = [m_n] =0$ then $\dim H_0 = 1$ and $\dim H_1 = n$.
Otherwise $\dim H_0 = 0$ and $\dim H_1 = n-1$.
\end{corol}

\bigskip
\section{Quantum Groups}

\bigskip
\noindent
In what follows set $q^x = e^{{\frac{2 \pi i}{\kappa}}\cdot x}$ 
for any linear operator $x$.

\begin{defin} The algebra $U_q(\mathfrak{sl}_2)$ 
is the algebra, with unit $1$, generated by $e$,
$f$, and $q^xh$ (for $x \in \mathbb{C}$) subject to the
following relations:

\begin{alignat*}{2}
&q^{xh}e = eq^{x(h+2)} && \quad q^{xh}f = fq^{x(h-2)} \\
& [e,f] = \frac{q^{\frac{h}{2}} - q^{\frac{-h}{2}}}{{ q^{\frac{1}{2}} - q^{\frac{-1}{2}}}}
	&& \quad q^{xh}q^{x'h}=q^{(x+x')h}
\end{alignat*}

\bigskip
\noindent
We also require $q^{0h}=1$. The algebra $U_q(\mathfrak{sl}_2)$
is called the quantum group of $\mathfrak{sl}_2(\mathbb{C})$.
It is equipped with a Hopf algebra structure with
comultiplication $\Delta$ given by 

\[ \Delta(e) = e\otimes q^{\frac{h}{4}}
+ q^{\frac{-h}{4}}\otimes e \] \[ \Delta(f)= f \otimes q^{\frac{h}{4}}
+ q^{\frac{-h}{4}}\otimes f \] \[\Delta(q^{xh}) =
q^{xh} \otimes q^{xh}. \]
\end{defin}

\bigskip
\noindent
Note that when $q \rightarrow 1$ these relations
degenerate into the defining relations for
the classical enveloping algebra $U(\mathfrak{sl}_2)$.
By analogy with $\mathfrak{sl}_2(\mathbb{C})$, we
define the quantum Verma modules. For $\lambda \in \mathbb{C}$,
the Verma module $M^q_{\lambda}$ is the infinite dimensional
$U_q(\mathfrak{sl}_2)$-module generated by one vector $v_\lambda$ 
satisfying $e \cdot v_\lambda = 0$ and $q^{xh}\cdot v_\lambda =
q^x\lambda v_\lambda$. It has a basis given by 
$v_\lambda, fv_\lambda, f^2v_\lambda, \dots$ subject to

\[ f \cdot f^kv_\lambda = f^{k+1}v_\lambda \qquad
   q^{xH} \cdot f^kv_\lambda = q^{x(\lambda - 2k)}f^kv_\lambda \qquad
   e \cdot f^kv_\lambda = [k][\lambda-k+1]f^{k-1}v_\lambda. \]

\bigskip
\noindent
If $\lambda \in \mathbb{Z}_{\geq 0}$ then 
$M^q_\lambda$ possesses a non-trivial submodule spanned
by the vectors $f^{\lambda +1}v_\lambda,f^{\lambda +2}v_\lambda,
\dots$. The quotient $L^q_\lambda$ is irreducible. Again,
as $q \rightarrow 1$ the representations $L^q_\lambda$ degenerate
to the classical finite dimensional irreducible representations
$L_\lambda$ of $\mathfrak{sl}_2(\mathbb{C})$. If no confusion arises
we shall omit the superscript $q$ in the notation and
write only $M_\lambda$ and $L_\lambda$.

\bigskip
\noindent
Unlike the classical theory, the linear isomorphism $v \otimes w
\mapsto w \otimes v$ between $V \otimes W$ and $W \otimes V$
of two $U_q(\mathfrak{sl}_2)$-modules
$V$ and $W$ is {\bf not} an isomorphism of representations.
Indeed, the $U_q(\mathfrak{sl}_2)$-module isomorphism
is constructed using the $R$ matrix:

\[ R \ = \ q^{\frac{h \otimes h}{4}}
	\sum_{k \geq 0} q^{\frac{-k(k+1)}{4}}
		\frac{(q^{\frac 12} - q^{\frac{-1}2})^k}{[k]!} \ q^{\frac{kh}{4}}e^k
	\otimes q^{\frac{-kh}{4}} f^k \ \in U_q(\mathfrak{sl}_2)
	\otimes U_q(\mathfrak{sl}_2)
\]

\bigskip
\noindent
The above expression is in fact a finite sum when applied
to any particular vector in $V \otimes W$ since $e$ is
always locally nilpotent. The $U_q(\mathfrak{sl}_2)$-isomorphism
between $V \otimes W$ and $W \otimes V$ is $P \cdot R$ where
$P(v \otimes w) = w \otimes v$.

\bigskip
\noindent
As an example, consider the case of $L_m \otimes L_l$.
Its basis is given by vectors of the form $v_m \otimes v_l,
fv_m \otimes v_l, v_m \otimes fv_l, \dots$.
Assume $q$ is not a root of unity (i.e. $\kappa$ is not
rational). We will evaluate $P \cdot R$ on the first three
basis vectors. Since the highest power of $f$ occurring
in any of these vectors is 1 it is enough to truncate $R$ as

\[ q^{\frac{h \otimes h}{4}} \ + \ q^{\frac{h \otimes h}{4}} q^{\frac{-1}{2}}
   q^{\frac{h}{4}} e \otimes q^{\frac{-h}{4}}f. \] 

\bigskip
\noindent
The isomorphism $P \cdot R:L_m \otimes L_l \longrightarrow 
L_l \otimes L_m$ evaluates on $v_m \otimes v_l$, 
$fv_m \otimes v_l$, and $v_m \otimes fv_l$ as 

\begin{alignat*}{2} & v_m \otimes v_l && \longmapsto \quad
   q^{\frac{ml}{4}}v_l \otimes v_m  \\ 
   & v_m \otimes fv_l && \longmapsto \quad q^{\frac{m(l-2)}{4}}
   v_l \otimes fv_m \\ & fv_m \otimes v_l && \longmapsto \quad
   q^{\frac{l(m-2)}{4}} v_l \otimes fv_m + 
   [m] \ q^{\frac{ml-m-l}{4}} fv_l \otimes v_m
\end{alignat*}

\bigskip
\noindent
\begin{theor} Let $V_1$, $V_2$, and $V_3$ be 
$U_q(\mathfrak{sl}_2)$-modules then the following diagram is commutative

\begin{alignat*}{5}
 & V_1 \otimes V_2 \otimes V_3
   && \overset{(23)}{\longrightarrow} && V_1 \otimes V_3 \otimes V_2 
   && \overset{(12)}{\longrightarrow} &&  V_3 \otimes V_1 \otimes V_2 \\
   & \quad {\scriptstyle (12)} \Big\downarrow  &&  &&  && 
   && \quad \quad \ \Big\downarrow {\scriptstyle (23)} \\
   & V_2 \otimes  V_1 \otimes V_3
   && \overset{(23)}{\longrightarrow} && V_2 \otimes V_3 \otimes V_1
   && \overset{(12)}{\longrightarrow} && V_3 \otimes V_2 \otimes V_2
\end{alignat*}

\bigskip
\noindent
where $(12)= PR \otimes 1$ and $(23)=1 \otimes PR$.
\end{theor}

\bigskip
\noindent
Commutativity of the above diagram implies that
the $R$-matrix satisfies the Yang-Baxter equation.

\section{Braid Groups}

\begin{defin}
The Braid group $B_n$ on $n$-strings is the group with
generators $b_1,\dots,b_{n-1}$ with defining relations 

\bigskip
\indent \indent \indent $\bullet$ $b_ib_j = b_jb_i$ if \ $|i-j|>1$

\bigskip
\indent \indent \indent $\bullet$ $b_ib_{i+1}b_i = b_{i+1}b_ib_{i+1}$ 
for $1 \leq i \leq n-2$.

\bigskip
\noindent
The Pure Braid group $PB_n$ is defined as the kernel
of the homomorphism $\sigma: B_n \longrightarrow S_n$ 
given by $b_i \mapsto (i,i+1)$. 

\end{defin}

\bigskip
\noindent
An element of the Braid
group is visualized as a family of $n$ non-intersecting
paths in space joining $n$ points in a plane to $n$ point
in a parallel plane. Multiplication in $B_n$ is 
interpreted as concatenation of paths. The Pure Braid group
$PB_n$ is visualized as families of paths which start
and end at the same point - where the second plane
and the points are identified with the first. 

\bigskip
\noindent
Recall that $\Delta = \Big\{ \ (a_1,\dots,a_n) \in \mathbb{C}^n \ \Big| \ 
z^n + a_1z^{n-1} + \cdots + a_n \ \text{has multiple roots} \
\Big\}$. Every element in the compliment $\mathbb{C}^n - \Delta$ can
be identified with a unique monic $n$th order polynomial with
distinct roots; namely associate the point $(a_1,\dots,a_n)$
with the polynomial $(z-a_1)\cdots(z-a_n)$. Given a braid -
which we interpret as a family of non-intersecting paths
joining the distinct points $(a_1,0) \dots,(a_n,0) \in \mathbb{C} \times \{0\}
\subset \mathbb{R}^3$ and the points $(a_1,1),\dots,(a_n,1) \in
\mathbb{C} \times \{1\} \subset \mathbb{R}^3$- we may associate a 
loop in $\mathbb{C}^n - \Delta$ emanating from $(a_1,\dots,a_n)$.
The position of the loop at time $t \in [0,1]$ is the
point $(b_1,\dots ,b_n)$ where each $(b_i,t)$ is the intersection
point of a strand of the braid with the plane
$\mathbb{C} \times \{t\}$. Using this construction one obtains
the following theorem.

\begin{theor} The classifying space for the Braid
group $B_n$ is $\mathbb{C}^n - \Delta$. In other words,
$\pi_1\big(\mathbb{C}^n-\Delta \big) = B_n$ and 
$\pi_k\big(\mathbb{C}^n-\Delta \big) = 0$ for $k>1$.
\end{theor}

\bigskip
\noindent
Let $D = \Big\{ \ (z_1,\dots,z_n) \in \mathbb{C}^n \ \Big| \ 
z_i=z_j \ \exists \ i \ne j \ \Big\}$. By a similar
construction using pure braids one obtains the following
theorem.

\begin{theor}
The classifying space for the Pure Braid group $PB_n$ is
$\mathbb{C}^n - D$.
\end{theor}

\bigskip
\noindent
Let $V_1,\dots,V_n$ be $U_q(\mathfrak{sl}_2)$-modules and 
let $R_i^{\vee}: V_1 \otimes \cdots \otimes V_i \otimes
V_{i+1} \otimes \cdots \otimes V_n \longrightarrow
V_1 \otimes \cdots \otimes V_{i+1} \otimes V_i \otimes
\cdots \otimes V_n$ be the map $P\cdot R$ in the $i,i+1$
positions and identity in all others. 

\begin{theor}
For $1 \leq i \leq n-2$ we have the Yang-Baxter identity
$R_i^{\vee}R_{i+1}^{\vee}R_i^{\vee} = R_{i+1}^{\vee}
R_i^{\vee}R_{i+1}^{\vee}$.
\end{theor}

\bigskip
\noindent
Let $V$ be a $U_q(\mathfrak{sl}_2)$-module. In view
of Theorem 4 we can surmise that the
map $B_n \longrightarrow \operatorname{GL}(V^{\otimes n})$
given by $b_i \mapsto R_i^{\vee}$ is a representation of
$B_n$. Moreover, since $\big(R_i^{\vee}\big)^2$ maps
$V_1 \otimes \cdots \otimes V_n$ to itself it follows
that the map $PB_n \longrightarrow \operatorname{GL}(
V_1 \otimes \cdots \otimes V_n)$ given by $b_i \mapsto R_i^{\vee}$ is
a representation of the Pure Braid group $PB_n$.

\section{Quantum Singular Vectors}

\bigskip
\noindent
For $m = (m_1, \dots,m_n) \in \mathbb{C}^n$ let $M^{\otimes m}$ be
the $U_q(\mathfrak{sl}_2)$ Verma module 
$M_{m_1} \otimes \cdots \otimes M_{m_n}$. 
Let $\Sing M^{\otimes m} \Big[ |m|-2k \Big] =
\Big\{ \ w \in M^{\otimes m} \ \Big| \ e \cdot w = 0 \ \text{and}
\ q^{\lambda h} \cdot w = q^{\lambda(|m|-2k)}w \ \Big\}$. 
Note that $e$ acts as $e\otimes q^{\frac{h}{4}} \otimes
\cdots \otimes q^{\frac{h}{4}} \ + \ \cdots \ + \
q^{\frac{-h}{4}} \otimes \cdots \otimes q^{\frac{-h}{4}}
\otimes e$. As in the classical case,
$\Sing M^{\otimes m} \Big[|m| \Big]$ is spanned by
the vector $v_{m_1}\otimes \cdots \otimes v_{m_n}$.
When $k=1$ we see that $\Sing M^{\otimes m} \Big[ |m| -2 \Big]$
is given by

\[   \Bigg\{ \ w=\sum_{j=1}^{n} I_j v_{m_1} \otimes 
     \cdots \otimes fv_{m_j} \otimes \cdots \otimes v_{m_n}
     \ \Bigg| \ \text{eigenspace condition and} \]
  
\[   \sum_{j=1}^n I_j[m_j] q^{\frac{1}{4}\big(
     \sum_{i>j}m_i - \sum_{i<j}m_i \big)} = 0 \ \Bigg\}.  \]

\begin{corol}
The map $v_{m_1}\otimes \cdots \otimes fv_{m_j} \otimes
\cdots \otimes v_{m_n} \overset{\psi}{\longmapsto}
\omega_j = c_j(l_j,s_j)$ and $v_{m_1} \otimes \cdots
\otimes v_{m_n} \overset{\psi}{\longmapsto}
\omega_0 = c_0(p_0,s_0)$ makes the following
diagram commute:

\begin{alignat*}{2}  0 \longrightarrow &M^{\otimes m}\Big[|m| -2 \Big]
   \overset{e}{\longrightarrow} &&M^{\otimes m} \Big[
   |m| \Big] \longrightarrow 0 \\
   &{\scriptstyle \psi} \Big\downarrow &&{\scriptstyle \psi}
   \Big\downarrow \\
   0 \longrightarrow &\bigoplus_{j=1}^n \mathbb{C}\cdot \omega_j
   \quad \overset{d_1}{\longrightarrow} &&\mathbb{C}\cdot \omega_0
   \longrightarrow 0
\end{alignat*}

\end{corol}

\bigskip
\noindent
As an immediate consequence of this corollary we obtain

\[ \Sing M^{\otimes m}\Big[ |m|-2 \Big] \ \cong \ 
   H_1 \Big( \mathbb{C} - \{z_1,\dots,z_n\}, \Phi_{1,n}^{\frac{1}{\kappa}} \Big) \]

\bigskip
\noindent
for $z_1 < \cdots < z_n$ real.

\begin{reftc}
\cite{V2}.
\end{reftc}


\chapter{}
\section{Monodromy of KZ equations}
KZ equations are defined over $U=\CC^n-D$, where
$$D=\{z\in\CC^n\ |\ \exists i,j\mbox{ s.t. }z_i=z_j\}.$$
For any $A\in U$, denote $\sol_A$ the space of all solutions of the
equations in a neighborhood of $A$. For any path $p$ from $A$ to $B$
in $U$ the continuation of solutions along $p$ gives an isomorphism
$$\al_p:\sol_A\to\sol_B.$$
If $A,B$ are fixed, then $\al_p$ does not depend on continuous deformations
of the path. If $p$ is a loop, then the isomorphism becomes an automorphism
$$M(p):\sol_A\to\sol_A$$
and gives the monodromy representation
\begin{equation}\label{rms}
\pi_1(\CC^n-D)\to\GL(\sol_A).
\end{equation}

\section{Topological monodromy}
Consider on $\CC^n-D$ the vector bundle whose fiber on $\nup zn$ is
$H_1(\CC-\{z_1,\ldots,z_n\},\pkun)$. There is a flat connection on it:
we deform the cycles for close points $\nup zn\sim\nup{z'}n$
(this is known as the \emph{Gauss--Manin connection}). In particular,
if we have a loop in $\CC^n-D$ we get a map
$$H_1(\CC-\{z_1,\ldots,z_n\},\pkun)\to H_1(\CC-\{z_1,\ldots,z_n\},\pkun).$$
We get thus another action of the group 
\begin{equation}\label{rmt}
\pi_1(\CC^n-D)\to\Aut(H_1).
\end{equation}

\section{$R$-matrix}
Recall from the previous lecture that we have an action of the pure braid group
\begin{equation}\label{rrm}
PB_n\to\Aut(\smd)
\end{equation}
for $M^{\otimes m}=M_{m_1}\oco M_{m_n}$.
Recall also that $PB_n=\pi_1(\CC^n-D)$.

\section{Relationship}
\begin{claim}
\begin{itemize}
\item[(a)] Via the correspondence of the space of solutions of KZ with $H_1$, we get
	an isomorphism between the two first representations (namely, \eqref{rms} and
	\eqref{rmt}).
\item[(b)] Via the correspondence of $\smd$ with $H_1(\CC^n-\{z_1,\ldots,z_n\},\pkun)$, we get
	an isomorphism between the second and third representations (namely,
	\eqref{rmt} and \eqref{rrm}).
\end{itemize}
\end{claim}
The first part of the claim is clear; the second one is remarkable.
It was proved first by Kohno without hypergeometric functions and then by Drinfeld,
who developed for this purpose the theory of quasi Hopf algebras.

We give a sketch of a third proof.
Consider the KZ equation with values in $\smk$; solutions of it are given, as we know, by
\begin{align*}
\pg\nup zn &=\sum_{|J|=k}\int_\gm\pkkn(t,z,m)A_J(t,z)\;\dkf tk \\
\Phi_{k,n} &=\pzzsk\pttsk\ptzsk
\end{align*}
For $z\in\CC^n$, let $\ckn=\{t\in\CC^k\ |\ t_i\neq z_j,\ t_i\neq t_j\ \forall i,j\}$.
The $k$-cycles $\gm$ are elements of $H_k(\ckn,\pkkn)$ and each $\gm$ of this type defines
a solution of KZ. However, not all $\gm$ give a nonzero solution. Notice that
the differential forms $\eta_J=\pkkn A_J\;\dkf tk$ are skew symmetric with respect to
permutations of $t_1,\ldots,t_k$:
$$\eta(t_{\sg(1)},\ldots,t_{\sg(k)})=(-1)^{|\sg|} \eta\nup tk\quad\forall\sg\in\Ss_k.$$
The symmetric group $\Ss_k$ acts also on the $k$-cells. This action induces
an action of $\Ss_k$ on the homology spaces. Let $H_k(\ckn,\pk)^-$ be the
subspace of skew symmetric elements $\gm$ (i.e., $\sg\gm=(-1)^{|\sg|}\gm$).
Let $H_k^\sim$ be the sum of all other isotypical components, $H_k=H_k^-\oplus H_k^\sim$.
Then, for any $\gm\in H_k^\sim$ we have $\pg\equiv 0$. Thus, solutions of KZ are
labeled by $H_k(\ckn,\pk)^-$.

\begin{theor}
\begin{itemize}
\item[(a)] Let $z\in\CC^n$ be such that $z_1<\cdots<z_n$. Then there is a
	natural isomorphism
$$H_k(\ckn,\pk)^-\simeq\smkq.$$
	Thus, the solutions of KZ are labeled by this space.
\item[(b)] Under this identification, the monodromy of KZ is identified with the
	$R$-matrix representation of $PB_n$ in $\smkq$.
\end{itemize}
\end{theor}
The construction of the isomorphism is analogous to the construction of the
isomorphism for $k=1$. The vector $f_Jv=f^{j_1}v_{m_1}\oco f^{j_n}v_{m_n}$
is identified with the following cell: take $j_i$ loops around each $z_i$
(see Figure \ref{tetas}). This is a product of $k$ $1$-dimensional cells,
and then it is a $k$-dimensional cell. Then, take the antisymmetrization of it.
\begin{figure}[ht]
\input{celdas.pstex_t}
\caption{}\label{tetas}
\end{figure}

This proof was developed later than that of Drinfeld. In it, one realizes
KZ equations geometrically, as equations for hypergeometric functions
associated with the master functions $\Phi_{k,n}$. In this geometric approach
the isomorphism of the KZ monodromy with the $R$-matrix representation is
given explicitly, through the hypergeometric pairing
$$\gm\otimes\pk A_J\;dt\mapsto\int_\gm\pk A_J\;dt.$$

We address now the following
\begin{probl}
Show that the construction we made gives many linearly independent solutions.
\end{probl}
This will be solved in the next lecture; we develop here some tools.

\section{Remark on integration cycles}\ \\
If $\ka,m_1,\ldots,m_n\in\RR$, $m_1,\ldots,m_n<0$, $\ka>0$ and $\ka\gg 1$,
then all exponents of $\Phi_{k,n}$ are big positive integers. Assume $z\in\RR^n$.
We can consider in $\RR^k$ the hyperplanes $\{t_i=z_j\}$ and $\{t_i=t_j\}$.
Call $\ckn(z)=\RR^k-\{t\in\RR^k\ |\ t_i=t_j,\ t_i=z_j\}$ complement to
these hyperplanes. For $k=2$ we can picture it as in Figure \ref{hyperp}.
\begin{figure}[ht]
\input{hyperp.pstex_t}
\caption{Hyperplanes in $\RR^k$ ($k=2$, $n=4$)}\label{hyperp}
\end{figure}
We can take now $\gm$ as any bounded domain in $\ckn(z)$ (we shaded two
of such domains in figure \ref{hyperp}).
\begin{claim}
Under these assumptions, $\pg\nup zn$ is a solution of the KZ equation with
values in $\smk$.
\end{claim}

\section{Example: the Selberg integral}\

\noindent Let $n=2$; consider KZ with values in
$$\Sing M_{m_1}\otimes M_{m_2}[m_1+m_2-2k]
	=\{\varphi=\sum_{j_1+j_2=k}I_{(j_1,j_2)}
	f^{j_1}v_{m_1}\otimes f^{j_2}v_{m_2}\ |\ e\varphi=0\}.$$
Calculating explicitly the condition $e\varphi=0$ we can see that the
space is one dimensional and is generated by
$$\om=\sum_{j_1+j_2=k}\frac{(-1)^{j_1}}{j_1!j_2!}
	\prod_{i=0}^{j_1-1}\frac 1{m_1-i}\prod_{i=0}^{j_2-1}\frac 1{m_2-i}
	f^{j_1}v_{m_1}\otimes f^{j_2}v_{m_2}.$$
The hypergeometric solutions have the form
$$\pg(z_1,z_2)=\sum_{j_1+j_2=k}\int_{\Dt(z_1,z_2)}\Phi^{\frac 1\ka}_{k,2}
	A_{(j_1,j_2)}\;\dkf tk\;f^{j_1}v_{m_1}\otimes f^{j_2}v_{m_2},$$
where $\Dt(z_1,z_2)=\{t\in\RR^k\ |\ z_1\le t_k\le\cdots\le t_1\le z_2\}$.
Then $\varphi(z_1,z_2)$ is proportional to $\om$.
To see that $\varphi\not\equiv 0$, we consider the $(k,0)$ coordinate of
the solution (i.e., $(j_1,j_2)=(k,0)$):
$$(z_1-z_2)^{\frac{m_1m_2}{2\ka}}\int_{\Dt(z_1,z_2)}
	\prod_{i=1}^k(t_i-z_1)^{-\frac{m_1}\ka-1}(t_i-z_2)^{-\frac{m_2}\ka}
	\prod(t_i-t_j)^{\frac 2\ka}\;\dkf tk.$$
Using new variables $t_i=(z_2-z_1)s_i+z_1$, we can write this as
$$(z_1-z_2)^{\frac{\Dt(m_1+m_2-2k)}\ka-\frac{\Dt(m_1)}\ka-\frac{\Dt(m_2)}\ka}
	\frac 1{k!}\int_k\left(-\frac{m_1}\ka,-\frac{m_2}\ka+1,\frac 1\ka\right),$$
where $\Dt(a)=\frac{a(a+2)}4$.
Here $\int_k(a,b,c)$ is the Selberg integral:
$$\int_k(a,b,c)=\int_{0\le s_k<\cdots<s_1\le 1}
	\prod_{i=1}^ks_i^{a-1}(1-s_i)^{b-1}
	\prod_{1\le i<j\le k}(s_i-s_j)^{2c}\;\dkf sk.$$
In 1944 Selberg showed that
$$\int_k(a,b,c)=\prod_{j=0}^{k-1}\frac{\Gm(1+c+jc)}{\Gm(1+c)}\,
	\frac{\Gm(a+jc)\Gm(b+jc)}{\Gm(a+b+(k+j-1)c)}.$$
For $n=1$, we get the beta integral $\int_1(a,b,c)=\frac{\Gm(a)\Gm(b)}{\Gm(a+b)}$.
The Selberg integral is one of the most remarkable hypergeometric functions.
In particular, we see that our hypergeometric solution is not zero for
generic $m_1,m_2,\ka$.

Taking a suitable limit of the Selberg integral we get
$$\int_{\RR^k}e^{-a\sum_{i=1}^ks_i^2}\prod_{i<j}(s_i-s_j)^{2c}\;\dkf sk
	=(2\pi)^{\frac k2}(2a)^{-k\frac{c(k-1)+1}2}
		\prod_{j=1}^k\frac{\Gm(1+jc)}{\Gm(1+c)}.$$

\section{Connection with finite reflection groups}\ 

\noindent Macdonald (1982) observed that $s_i-s_j=0$, $1\le i<j\le k$
are equations of the reflecting hyperplanes of the finite group $A_{k-1}$.
There are other finite groups generated by reflections. In the $k$-dimensional
Euclidean space $\RR^k$ consider a certain number of hyperplanes, all passing
through the origin. If the angles between the hyperplanes are properly chosen,
then the group generated by reflections corresponding to them is finite.
These groups were enumerated and classified by Coxeter as
$$A_k\ (k\ge 1),\ B_k\ (k\ge 2),\ D_k\ (k\ge 4)\ E_6,\ E_7,\ E_8,\ 
	F_4,\ G_2,\ H_3,\ H_4,\ I_2(p)\ (p>5).$$
Let $G$ be one of these reflecting groups. Let $P(s)$ be the product
of distances of the point $\nup sk=s$ from all reflecting hyperplanes
belonging to $G$. Let $N$ be the number of hyperplanes.
Not knowing Selberg integrals, Macdonald conjectured that
$$\int_{\RR^k}e^{-\sum s_i^2/2}|P(s)|^{2c}\;\dkf sk
	=2^{-Nc}(2\pi)^{k/2}\prod_{j=1}^k\frac{\Gm(1+cd_j)}{\Gm(1+c)},$$
where $d_j$ are the degrees of a basis set of the space of
homogeneous polynomials which are invariant with respect to $G$.

The conjecture was proved by Opdam. The $q$-analogs were proved by Cherednik.

\begin{reftc}
\cite{V1,V2,V4}
\end{reftc}


\chapter{}
The goal of this lecture is to see how many independent solutions one can get
by taking cycles on the bounded domains in the case $\ka>0$, $m_i<0\ \forall i$.
\section{Determinant formulas}\ \\
Consider the KZ equation with values in $\smd$:
\begin{align*}
\Phi_{1,n} &=\prod_{1\le i<j\le n}(z_i-z_j)^{\frac{m_im_j}2}\prod_{l=1}^n(t-z_l)^{-m_l}, \\
\pg\nup zn &=\sum_{j=1}^n\int_\gm\Phi_{1,n}^{1/\ka}\frac{dt}{t-z_j}\;
	v_{m_1}\oco fv_{m_j}\oco v_{m_n}
\end{align*}
and call $\I j(z)=\int_\gm\Phi_{1,n}^{1/\ka}\frac{dt}{t-z_j}$.
We have $\I 1(z)=-\frac 1{m_1}\sum_{j=2}^nm_j\I j(z)$, which is the condition for
the vector $\sum_j\I j v_{m_1}\oco fv_{m_j}\oco v_{m_n}$
to be a singular vector.
\begin{corol}
If $m_1\neq 0$, then $\dim\smd=n-1$.
\end{corol}
Let $z_1<\ldots<z_n$ be real. Let $\Dt_i(z)=[z_{i-1},z_i]$, $i=2,\ldots,n$.
Let $\ka,m_1,\ldots,m_n\in\RR$, $\ka>0$, $m_1,\ldots,m_n<0$. Then
$\varphi^{(\Dt_i)}\nup zn$, $i=2,\ldots,n$ are solutions of KZ.
\begin{theor}
These solutions are linearly independent for generic $m_1,\ldots,m_n,\ka$.
\end{theor}
The theorem follows from the following lemma
\begin{lemma}[Determinant Formula]
Let $\varphi^1,\ldots,\varphi^{n-1}$ be a set of solutions. Let $w^1,\ldots,w^{n-1}$
be a basis of $\smd$, and write $\varphi^i=\sum_j\varphi^i_jw^j$. Then
$$\det(\varphi^i_j)=(\const)\prod_{i<j}
	(z_i-z_j)^{(n-1)\frac{m_im_j}{2\ka}-\frac{m_i}\ka-\frac{m_j}\ka}$$
\end{lemma}
\begin{proof}
Notice that if one has a differential equation $\frac{dy}{dx}=A(x)y$, where
$A,y$ are $n\times n$ matrices, then $\det'y=\tr A(x)\det y$
(since $\det(I+\ep A)=1+\ep \tr A+o(\ep^2)$). Our system is
$$\frac{\ptl y}{\ptl z_i}=\frac 1\ka\sum_{j\neq i}\frac{\oij}{z_i-z_j}y,$$
and then
$\frac{\ptl\det y}{\ptl z_i}=\frac 1\ka\sum_{j\neq i}\frac{\tr\oij}{z_i-z_j}\det y$,
whence $\det y=(\const)\prod_{i<j}(z_i-z_j)^{\frac{\tr\oij}\ka}$.

We show now that
$\tr\oij=\tr\oij|_{\smd}=(n-1)\frac{m_im_j}{2\ka}-\frac{m_i}\ka-\frac{m_j}\ka$.
We prove it for $\Omega^{(12)}$, the general case is analogous. We can use the basis
\begin{align*}
&w^1=\frac 1{m_1}fv\otimes v\oco v-\frac 1{m_2}v\otimes fv\oco v, \\
&w^i=\frac 1{m_1+m_2}(fv\oco v+v\otimes fv\oco v)-\frac 1{m_i}v\oco fv\oco v
\intertext{and then}
&\Omega^{(12)}w^1=((\frac 12h\otimes h+e\otimes f+f\otimes e)\otimes 1\oco 1)(w^1) \\
	&\phantom{\Omega^{(12)}w^1}=(\frac 12m_1m_2-m_1-m_2)w^1 \\
&\Omega^{(12)}w^i=\frac 12m_1m_2w^i.
\end{align*}
\end{proof}
We compute now the constant in the determinant.
Let $\bar\Phi=\prod_{l=1}^n(t-z_l)^{\al_l}$, and the differential forms
$\eta_j=\frac{d(t-z_j)}{t-z_j}$, $j=1,\ldots,n$. Let $z_1<\ldots<z_n$ as before.
We fix then the arguments in each interval $\Dt_i(z)$, and then
the integrals $\int_{z_{i-1}}^{z_i}\bar\Phi\eta_j$ are well defined.
\begin{theor}
$$\det_{2\le i,j\le n}(\al_j\int_{z_{i-1}}^{z_i}\bar\Phi\eta_j)
	=\frac{\Gm(\al_1+1)\cdots\Gm(\al_n+1)}{\Gm(\al_1+\cdots+\al_n+1)}
	\prod_{i<j}(z_i-z_j)^{\al_i+\al_j}.$$
\end{theor}
Before giving the proof we notice that for $n=2$, $z_1=0$ and $z_2=1$
the theorem says
$$\al_2\int_0^1t^{\al_1}(1-t)^{\al_2-1}dt
	=\frac{\Gm(\al_1+1)\Gm(\al_2+1)}{\Gm(\al_1+\al_2+1)},$$
and this is Euler's formula.

\begin{proof}
By induction. For $n=2$, we saw that this is Euler's formula. Consider then
$z_1,\ldots,z_{n-1}$ fixed and let $z_n\to\infty$; we will see that the
ratio RHS/LHS tends to $1$.

We remark that
\begin{enumerate}
\item For $i,j\le n-1$,
$$\int_{z_{i-1}}^{z_i}\al_j\bar\Phi\frac{dt}{t-z_j}=(-z_n)^{\al_n}
	\left[\int_{z_{i-1}}^{z_i}\al_j\prod_{l=1}^{n-1}(t-z_l)^{\al_l}
	\frac{dt}{t-z_j}+\Oc(\frac 1{z_n})\right]$$
\item $\int_{z_{i-1}}^{z_i}\al_n\bar\Phi\frac{dt}{t-z_n}
	=\Oc((-z_n)^{\al_n-1})$
\item $\int_{z_{n-1}}^{z_n}\al_n\bar\Phi\frac{dt}{t-z_i}
	=\Oc((-z_n)^{\al_1+\cdots+\al_n})$
\item $\int_{z_{n-1}}^{z_n}\al_n\bar\Phi\frac{dt}{t-z_n}
	=(-z_n)^{\al_1+\cdots+\al_n}\int_0^1\al_nt^{\al_1+\cdots+\al_{n-1}}(1-t)^{\al_n-1}\;dt
	+\Oc(z_n^{\al_n})$.
\end{enumerate}
We have then a matrix
\begin{center}
\begin{tabular}{l|c|c|c|c|}
 & $\bar\Phi\frac{dt}{t-z_1}$ & \hspace*{1cm} $\cdots\cdots$ \hspace*{1cm}
 	& $\bar\Phi\frac{dt}{t-z_{n-1}}$ & $\bar\Phi\frac{dt}{t-z_n}$ \\ \cline{2-5}
$\Dt_1$ & & & & $\Oc((-z_n)^{\al_n-1})$ \\ \cline{2-5}
	& & & & $\vdots$ \\ \cline{2-5}
$\Dt_{n-1}$ & & & & $\Oc((-z_n)^{\al_n-1})$ \\ \cline{2-5}
$\Dt_n$ & $\Oc((-z_n)^{\al_1+\cdots+\al_n})$ & $\cdots\cdots$
	& $\Oc((-z_n)^{\al_1+\cdots+\al_n})$ & \\ \cline{2-5}
\end{tabular}
\end{center}
Thus, dividing the $i$-th row ($i<n$) by $(-z_n)^{\al_n}$ and the
$n$-th row by $(-z_n)^{\al_1+\cdots+\al_n}$, we get
$$\det_{(n)}\nup zn=(-z_n)^{(n-1)\al_n}\det_{(n-1)}\nup z{n-1}
	\frac{\Gm(\al_1+\cdots+\al_{n-1}+1)\Gm(\al_n+1)}{\Gm(\al_1+\cdots+\al_n+1)}.$$

Now, it is easy to proof that RHS/LHS is holomorphic with respect to $z_n$. Since
its limit at $\infty$ is $1$, we finished.
\end{proof}

\begin{remar}
We proved that we constructed all solutions for generic $m_1,\ldots,m_n,\ka$.
The case of degenerate parameters is very explicit: it happens in the poles of
$\Gm(\al_1+\cdots+\al_n+1)$.
\end{remar}
\begin{remar}
Considering solutions with values in $\smk$, one can similarly show that the set
of hypergeometric solutions generate all solutions for generic $m_1,\ldots,m_n,\ka$.
For $n=2$, one has to replace Euler's formula by the Selberg integral.
\end{remar}

There are remarkable general determinant formulas by Terao and Douai.
Let $A$ be a finite collection of hyperplanes in $\RR^k$ and chose
for each $H\in A$ a number $\lm_H$ and a linear function $f_H$. Then
$$\det\left(\int_\Dt\prod_{H\in A}f_H^{\lm_H}\om\right)
	=\left(\text{Alternating product of $\Gm$ functions}\right)
		\prod\left(\text{values of $f_H^{\lm_H}$ at vertices}\right),$$
where $\Dt$ runs through bounded domains and $\om$ through a suitable basis of
logarithmic differential forms.
In our case there is an action of the symmetric group and one must consider
isotypical components.

Another determinant formula: let $A$ be again an arrangement of hyperplanes in
$\RR^k$, now in generic position. It can be proved that the number of bounded
domains is $\lcomb{n-1}k$. On the other hand, the number of monomials
$t^a=t_1^{a_1}\cdots t_k^{a_k}$ such that $a_1+\cdots+a_k\le n-k-1$ is also
$\lcomb{n-1}k$. One can consider $\det(\int_\Dt t^a\;dt_1\cdots dt_k)$, where
$\Dt$ runs on the bounded domains and $t^a$ runs on the preceding monomials.
We have
$$\det=\pm [(n-1)!]^{\lcomb{n-2}{k-1}}
	\prod_{S\subset\RR^k}(k!\vol(S)),$$
where $S$ runs on the $k$-simplices with faces in $A$.

Yet another example: let $z_1<\cdots<z_n\in\RR$. Consider $\det(\int_{z_i}^{z_{i+1}}t^j\;dt)$.
This is a Vandermonde determinant, and we have
$$\det(\int_{z_i}^{z_{i+1}}t^j\;dt)=\const\prod_{j<i}(z_i-z_j).$$

\section{Resonances}
Consider the case $k=1$. Then
$$\det(\text{KZ Solutions})=\frac{\prod_j\Gm(-\frac{m_j}\ka+1)}
	{\Gm(-\frac{\sum_im_i}\ka+1)}\prod(z_i-z_j)^{\fbox{*}}.$$
When $-\frac{\sum_im_i}\ka+1=0$ then we have a pole in the denominator
and then we get linearly dependent solutions. In this case, the bundle
$\smd$ has a nontrivial subbundle which is invariant by the KZ equation.
We want to describe this subbundle. It turns out that in CFT the appearance
of such subspace was predicted, and it was described in terms of
representation theory.

For a general $k$, let $m_1,\ldots,m_n\in\ZZ_{>0}$ and consider
$L^{\otimes m}=L_{m_1}\oco L_{m_n}$ the tensor product of irreducible
f.d. representations. Assume that $\ka\in\ZZ$, $\ka>2$,
$0\le m_1,\ldots,m_n$ and $|m|-2\le\ka-2$. Consider the KZ equation
with values in
$$\slk=\{v\in L^{\otimes m}\ |\ hv=(|m|-2k)v,\ ev=0\}.$$

Set $p=\ka-1-|m|+2k$, this is a positive integer. Take
$W\nup zn\subset\slk$ given by
$$W(z)=\{w\in L^{\otimes m}\ |\ hw=(|m|-2k)w,\quad
	ew=0,\quad (ze)^pw=0\},$$
where $ze:L^{\otimes m}\to L^{\otimes m}$ is the operator
$w_1\oco w_n\mapsto\sum_iz_iw_1\oco ew_i\oco w_n$.

\begin{theor}
The subspaces $W(z)\subset\slk$ are preserved by KZ.
\end{theor}
\begin{defin}
These subspaces are called \emph{spaces of conformal blocks}.
\end{defin}
\begin{theor}
All the hypergeometric solutions $\pg(z)$ belong to the bundle of
conformal blocks.
\end{theor}
This means that hypergeometric functions ``know" about CFT and its conformal blocks.

Consider the case $k=1$. We have
$$\sld=\{w=\sum_{j=1}^nI_jv_{m_1}\oco fv_{m_j}\oco v_{m_n}\ |\ 
	\sum_{j=1}^nI_jm_j=0\}.$$
Here $p=\ka+1-|m|\in\ZZ_{>0}$. If $p>1$ then $W(z)=\sld$
(by a degree argument $(ze)^2=0$). If $p=1$ then
$$W(z)=\{w=\sum I_jv_{m_1}\oco fv_{m_j}\oco v_{m_n}\ |\ 
	\sum m_jI_j=0,\quad\sum z_jm_jI_j=0\}.$$
Notice that $p=1\iff\ka=\sum m_j$. Then
\begin{theor}
For any hypergeometric solution $\pg(z)$ we have $\sum z_jm_jI_j=0$.
\end{theor}
\begin{proof}
It is similar to the proof that $\sum m_jI_j=0$. In that case, we have
$$-\ka d_t(\Phi^{1/\ka})=\sum m_j\Phi^{1/\ka}\frac{dt}{t-z_j},$$
and integrating over a cycle $\gm$ we get the identity.
Now, we have $\ka d_t(t\Phi^{1/\ka})=\sum z_jm_j\Phi^{1/\ka}\frac{dt}{t-z_j}$, i.e.,
\begin{align*}
d_t(t\prod(z_i-z_j)^{\mbox{*}}\prod(t-z_l)^{-m_l/\ka})
	&=\Phi^{1/\ka}dt+\Phi^{1/\ka}\sum-\frac{m_l}\ka\frac{(t-z_l)+z_l}{t-z_l} \\
	&=(1-\sum_{l=1}^n\frac{m_l}\ka)\Phi^{1/\ka}dt
		-\sum z_l\frac{m_l}{\ka}\Phi^{1/\ka}\frac{dt}{t-z_l},
\end{align*}
and, again, integrating over a cycle $\gm$ we get the identity.
\end{proof}

\section{Dynamical equations}
Let $\pkkn=\pzz\ptt\ptz$. Let
\begin{align*}
\pg\nup zn &=\sum_{|J|=k}\I Jf_Jv, \\
\I J\nup zn &=\int_\gm\pk A_j\;\dkf tk.
\end{align*}
be solutions of KZ with values in $\smk$, $\gm\in H_k(\Cc_{k,n}(z),\pkkn)$.
Shift $m_1\mapsto m_1+\ka$. Then the integrand is multiplied by
$\frac{\const(z)}{(t_1-z_1)\cdots (t_k-z_1)}$. This is a rational univalued
function; the same cycle $\gm$ is good for the new integral.
This defines a map
$$\{\I J(z,m_1,\ldots,m_n)\}\to\{\I J(z,m_1+\ka,\ldots,m_n)\}$$
in $H^1(\Cc_{k,n}(z),\pkkn)$.
\begin{probl}
Describe this map.
\end{probl}

Again, the answer will be given in terms of $\slt$-representation theory.
Relations of this type will be called \emph{dynamical equations}. Before going over
the problem, we make a

\begin{remar}
Since $\pg\in\smk$, there are many linear relations between $\{\I J\}$.
\end{remar}
\begin{lemma}
If $m_1\notin\{0,1,\ldots,k-1\}$, then for any $J_0$, $|J_0|=k$, the integral
$\I{J_0}$ is a linear combination of integrals $\{\I J\}$ with $J$ having the
form $(0,j_2,\ldots,j_n)$, $j_1+\cdots+j_n=k$. The coefficients of this linear
combination do not depend on $\gm$.
\end{lemma}
For example, for $k=1$ we have $I_1=\frac 1{m_1}\sum_{j=2}^nm_jI_j$. We will
study the transformation $m_1\mapsto m_1+\ka$ for ``restricted vector valued functions"
$$\psi^{(\gm)}\nup zn=\sum_{\stackrel{J=(0,j_2,\ldots,j_n)}{|J|=k}}\I Jf_Jv.$$
\begin{remar}
Having the restricted function $\psi^{(\gm)}$ one can reconstruct the initial
function $\pg$. Therefore one can rewrite the KZ differential equation
for $\pg(z)$ as a suitable differential equation for the restricted function.
This new differential equation is called \emph{trigonometric KZ equation}
\end{remar}

\begin{reftc}
\cite{V1,V2,V3,V4}, \cite{MTV}, \cite{DT}, \cite{FeSV}
\end{reftc}


\chapter{}
The goal of this lecture is to describe the trigonometric KZ equation and
dynamical equation and see their interaction.

\section{Dynamical equations}
Let
\begin{align*}
\pkkn&=\pzz\ptt\ptz \\
\I J\nup zn &=\int_\gm\pk A_j\;\dkf tk, \\
\pg\nup zn &=\sum_{|J|=k}\I Jf_Jv.
\end{align*}
be solutions of KZ with values in $\smk$, $\gm\in H_k(\Cc_{k,n}(z),\pkkn)$.
Shift $m_1\mapsto m_1+\ka$. Then the integrand is multiplied by
$\frac{\const(z)}{(t_1-z_1)\cdots (t_k-z_1)}$. This is a rational univalued
function; the same cycle $\gm$ is good for the new integral.
This defines a map in cohomology, given by
$$\{\I J(z,m_1,\ldots,m_n)\}\to\{\I J(z,m_1+\ka,\ldots,m_n)\}.$$
\begin{probl}
Describe this map.
\end{probl}
The answer will be given in terms of $\slt$-representation theory.
Relations of this type will be called \emph{dynamical equations}. Before going over
the problem, we make a

\begin{remar}
Since $\pg\in\smk$, there are many linear relations between $\{\I J\}$.
\end{remar}
\begin{lemma}
If $m_1\notin\{0,1,\ldots,k-1\}$, then for any $J_0$, $|J_0|=k$, the integral
$\I{J_0}$ is a linear combination of integrals $\{\I J\}$ with $J$ having the
form $(0,j_2,\ldots,j_n)$, $j_1+\cdots+j_n=k$. The coefficients of this linear
combination do not depend on $\gm$.
\end{lemma}
For example, for $k=1$ we have $I_1=\frac 1{m_1}\sum_{j=2}^nm_jI_j$. We will
study the transformation $m_1\mapsto m_1+\ka$ for ``restricted vector valued functions"
$$\psi^{(\gm)}\nup zn=
	\sum_{\begin{smallmatrix}J=(0,j_2,\ldots,j_n)\\|J|=k\end{smallmatrix}}\I Jf_Jv.$$
\begin{remar}
Having the restricted function $\psi^{(\gm)}$ one can reconstruct the initial
function $\pg$. Therefore one can rewrite the KZ differential equation
for $\pg(z)$ as a suitable differential equation for the restricted function.
This new differential equation is called \emph{trigonometric KZ equation}
\end{remar}

\section{Trigonometric $r$-matrix}
Set
$$\left.\begin{array}{l}
\Om^+=\frac 14h\otimes h+e\otimes f \\
\Om^-=\frac 14h\otimes h+f\otimes e
\end{array}\right\}\in\slt^{\otimes 2}$$
We have $\Om^++\Om^-=\Om$, the Casimir element.
Define the \emph{trigonometric $r$-matrix} $\gm(z)=\frac{\Om^+z+\Om^-}{z-1}$.
A more symmetric form is got when one considers $\gm(z_i/z_j)=\frac{\Om^+z_i+\Om^-z_j}{z_i-z_j}$.
Let $\hh=\CC h\subset\slt$ be the Cartan subalgebra. Let $V=V_2\oco V_n$
be a tensor product of $\slt$ modules, $\ka\in\CC$ and $\vlm\in\hh$.

\begin{defin}
The \emph{trigonometric KZ equations} on a $V$-valued function $u(z_2,\ldots,z_n,\vlm)\in V$
are the equations
$$\ka z_i\ddx {z_i}u=\sum_{j\neq i}\frac{(\Om^+z_i+\Om^-z_j)^{(i,j)}}{z_i-z_j}u+\vlm^{(i)}u,
	\quad i=2,\ldots,n,$$
where $\vlm^{(i)}=1\oco 1\otimes\vlm\otimes 1\oco 1$ in the $i$-th factor.
$\vlm$ and $\ka$ are parameters of the equation.

We introduce the \emph{trigonometric KZ operators}
$$\nb_i(\ka,\vlm)=\ka z_i\ddx{z_i}-\sum_{j\neq i}r^{(i,j)}(z_i/z_j)-\vlm^{(i)},
	\quad i=2,\ldots,n.$$
Then the trigonometric KZ equations have the form
$$\nb_i(\ka,\vlm)u=0,\quad i=2,\ldots,n.$$
\end{defin}

\section{Hypergeometric solutions of trigonometric equations}
We consider solutions with values in $V=M_{m_2}\oco M_{m_n}[\sum_{j=2}^nm_j-2k]$.
This is so because the trigonometric KZ equation does not preserve singular spaces,
though it does preserve weight spaces (there is a correspondence between singular
spaces in $V_{m_1}\oco V_{m_n}$ and weight spaces of $V_{m_2}\oco V_{m_n}$).
Take $\vlm=\lm\frac h2$, $\lm\in\CC$. Set $\nu=\sum_{j=2}^nm_j-2k$ and put
\begin{multline*}
\Phi(t_1,\ldots,t_k,z_2,\ldots,z_n) = \prod_{2\le i<j\le n}(z_i-z_j)^{\frac{m_im_j}2}
	\prod_{i=2}^nz_i^{\frac 12m_i(\lm+(m_i-\nu)/2)} \\
	\times \prod_{i=1}^k\prod_{j=2}^n(t_i-z_j)^{-m_j}
		\prod_{1\le i<j\le k}(t_i-t_j)^2\prod_{i=1}^kt_i^{-(\lm-1-\nu/2)}
\end{multline*}
For any $J=(j_2,\ldots,j_n)$, $|J|=k$, introduce the function $A_J$ by
the same formulas as before. Set $f_Jv=f^{j_2}v_{m_2}\oco f^{j_n}v_{m_n}$.
\begin{remar}
The complement to the singularities of $\Phi$ with respect to $t_1,\ldots,t_k$
is the space $\Cc_{k,n}\subset\CC^k$, which is also the complement to the singularities
of $\Phi_{k,n}$ with $z=(0,z_2,\ldots,z_n)$.
\end{remar}
\begin{theor}
For any $\gm(z)\in H_k(\Cc_{k,n}(z),\pk)$, the function
$$\pg\nup zn=\sum_{|J|=k}\int_{\gm(z)}\pk A_J\;dt f_Jv$$
is a solution of the trigonometric KZ equation with parameter
$\vlm=\lm\frac h2$.
\end{theor}
Let us describe then the Dynamical Equations $\vlm\mapsto\vlm+\ka\frac h2$
(this corresponds to $\lm\mapsto\lm+\ka$).
For $\lm\in\CC$, set
$$P(\lm)=\sum_{i=0}^{\infty}f^ie^i\frac 1{i!}\prod_{j=1}^i\frac 1{\lm-\frac h2-j}.$$
\begin{remar}
If $V$ is a f.d. $\slt$ module or a Verma module, then $P(\lm):V\to V$ is well
defined. Similarly, $P(\lm)$ is well defined on a tensor product of such modules.
\end{remar}
\begin{defin}
The \emph{Dynamical Equation} for $u\nupp z2n\in V_2\oco V_n$ is
$$u(z_2,\ldots,z_n,\lm+\ka)=\prod_{i=2}^n(z_i^{\frac h2})^{(i)}P(\lm)u(z_2,\ldots,z_n,\lm).$$
Denote $K(z_2,\ldots,z_n,\lm)=\prod_{i=2}^n(z_i^{\frac h2})^{(i)}P(\lm)$.
It is clear that $K(z,\lm)$ preserves the weight decomposition of $V$.
\end{defin}

\begin{theor}
The Dynamical Equation is compatible with the trigonometric KZ equation.
Explicitly,
$$\nb_i(\ka,\vlm+\ka\frac h2)K(z,\lm)=K(z,\lm)\nb_i(\ka,\vlm).$$
\end{theor}
\begin{corol}
If $u(z_2,\ldots,z_n,\lm)$ is a solution of the trigonometric KZ equation with
parameter $\vlm=\lm\frac h2$, then $K(z,\lm)u(z,\lm)$ is a solution of the
trigonometric KZ equation with parameter $\vlm+\ka\frac h2$.
\end{corol}
For hypergeometric solutions we can say more: consider a hypergeometric
solution of the trigonometric KZ equation with values
in $M_{m_2}\oco M_{m_n}[\sum_{i=2}^nm_i-2k]$,
$$\pg(z_2,\ldots,z_n,\lm)=\sum_{\begin{smallmatrix}J=\nupp j2n \\ |J|=k\end{smallmatrix}}
	\int_\gm\pk A_J\;\dkf tk.$$
\begin{theor}
The function $\pg$ satisfies also the Dynamic Equation
$$\pg(z_2,\ldots,z_n,\lm+\ka)=K(z,\lm)\pg(z_2,\ldots,z_n,\lm).$$
\end{theor}

\begin{remar}
The function $\pg(z_2,\ldots,z_n,\lm+\ka)$ is well defined since the corresponding integrand
differs from the integrand of $\pg(z_2,\ldots,z_n,\lm)$ by a univalued factor.
\end{remar}
\begin{examp}
Take $V=M_m$ and consider a weight subspace $V[m-2k]=\CC f^kv_m$.
Let $\vec u(z,\lm)\in V[m-2k]$, $\vec u(z,\lm)=u(z,\lm)f^kv_m$ for a scalar
$u(z,\lm)$. We describe trigonometric KZ, Dyn. Eqn. and its hypergeometric solutions
(here $\lm,z\in\CC$).

\noindent\textbf{Trigonometric KZ:}
$\ka z\ddx z\vec u=\vlm\vec u$, $\vlm=\lm\frac h2$.
Then the equation reads as
$$\ka z\ddx zu=\lm\frac{m-2k}2u.$$
Since $K(z,\lm)f^kv_m=z^{h/2}P(\lm)f^kv_m$, we must compute $P(\lm)f^kv_m$.
\end{examp}
\begin{theor}
We have
$$P(\lm)f^kv_m=\frac{(\lm+\frac m2)(\lm+\frac m2-1)\cdots(\lm+\frac m2-k+1)}
	{(\lm-\frac m2)(\lm-\frac m2+1)\cdots(\lm-\frac m2+k-1)}f^kv_m.$$
\end{theor}
\begin{proof}
\begin{align*}
P(\lm)f^kv_m &=\sum_{j\ge 0}\frac{k(k-1)\cdots(k-j+1)(m-k+1)\cdots(m-k+j)}
		{j!(\lm-\frac{m-2k}2-1)\cdots(\lm-\frac{m-2k}2-j)} f^kv_m \\
	&= F(-k,m-k+1,-\lm+\frac{m-2k}2+1,z=1)f^kv_m, 
\end{align*}
and by Gauss' formula we have
$F(a,b,c,1)=\frac{\Gm(c)\Gm(c-a-b)}{\Gm(c-a)\Gm(c-b)}$.
\end{proof}
\begin{corol}
The Dynamical Equation is
$$u(z,\lm+\ka)=z^{\frac{m-2k}2}\prod_{j=0}^{k-1}
	\frac{(\lm+\frac m2-j)}{(\lm-\frac m2+j)}u(z,\lm).$$
The hypergeometric solution is given by
\begin{align*}
u(z,\lm) &= \int_{0\le t_k\le\cdots\le t_1\le z}\pk A\;\dkf tk,\quad\text{where} \\
\pk &= z^{\frac 1{2\ka}(\lm+k)m}\prod_{i=1}^kt_k^{-\frac 1\ka(\lm-1-\frac{m-2k}2)}
	(t_i-z)^{-m/\ka}\prod_{i,j}(t_i-t_j)^{2/\ka} \\
A &= \frac 1{(t_1-z)\cdots(t_k-z)}.
\end{align*}
Changing variables $t_i=zs_i$ we get
$$u(z,\lm)=z^{\lm\frac{m-2k}{2\ka}}\int_k
	(-\frac 1\ka(\lm-1-\frac{m-2k}2)+1,-\frac m\ka,\frac 1\ka),$$
where as before $\int_k(a,b,c)=\prod_{j=0}^{k-1}\frac{\Gm(1+(j+1)c)}{\Gm(1+c)}
	\frac{\Gm(a+jc)\Gm(b+jc)}{\Gm(a+b+(k+j-1)c)}$ is the Selberg integral,
and we have the relation
$$\int_k(a+1,b,c)=\prod\frac{a+jc}{a+b+(k+j-1)c}\int_k(a,b,c),$$
from where we can see that $u(z,\lm)$ satisfies both equations.
\end{corol}

\section{Quasiclassical asymptotics of hypergeometric functions and Bethe ansatz}
Let $V=V_1\oco V_n$ be a tensor product of $\slt$ modules.
Set $H_i(z)=\sum_{j\neq i}\frac{\Om^{(ij)}}{z_i-z_j}$ as operators on $V$.
\begin{theor}
For any $z$, the operators $H_i(z)$ commute.
\end{theor}
\begin{proof}
Follows from the identity $[\Om,x\otimes 1+1\otimes x]=0$ $\forall x\in\slt$.
\end{proof}

The operators $\{H_i(z)\}$ are called \emph{the Hamiltonians} of the Gaudin model
of an ``inhomogeneous magnetic chain". In integrable models of quantum mechanics,
usually one has a vector space $V=V_1\oco V_n$ called the \emph{space of states}.
One may think about $n$ particles at the points $z_1,\ldots,z_n\in\CC$, and
the spaces $V_1,\ldots,V_n$ are their spin spaces. Usually one has a family of
commutative linear operators on $V$ called Hamiltonians. The main problem is to
diagonalize them and to find their eigenvalues. The next problem is to find the
limit when $n\to\infty$. The Bethe ansatz is a method for diagonalizing the
Hamiltonians. One looks for eigenvectors of the form
$$v\nup tk\in V,$$
where $v\nup tk$ is some suitable function of parameters $\vec t$. One then tries
to find the values of the parameters such that $v\nup tk$ is an eigenvector of the
Hamiltonians. One proves that if $t_1,\ldots,t_k$ satisfy some explicit system
of equations
\begin{equation}\label{eq1}
F_j\nup tk=0,\quad j=1,\ldots,k,
\end{equation}
then $v\nup tk\in V$ is an eigenvector.

The equations \eqref{eq1} are called \emph{the Bethe equations}, $v(t)$ is called
\emph{the Bethe vector}. The method is called \emph{the Bethe ansatz method}.

\titulo{The standard conjectures} are
\begin{enumerate}
\item The number of solutions properly counted $=\dim$ of the space of states.
\item Bethe vectors from a basis
\end{enumerate}
\titulo{Important observation:} In all known examples there exists a single (master)
function
$$\Phi\nup tk$$
such that $\frac{\ptl\Phi}{\ptl t_j}=F_j,\ j=1,\ldots,n$. Thus, with a
given integrable model of quantum mechanics there is associated a (master) function
$\Phi\nup tk$ such that
\titulo{Conjectures (reformulated)}
\begin{enumerate}
\item The number of critical points of $\Phi=\dim$ of the space of states,
\item the Bethe vectors $\{v(t_0)\}|_{t_0\in\mbox{crit. pts.}}$ form a basis,
\item the Hessian $\det(\frac{\ptl^2\ln\Phi}{\ptl t_i\ptl t_j})(t_0)$
	``=" square of the length of the Bethe vector $v(t_0)$
\end{enumerate}
\begin{remar}
For f.d. irreducible representations of $\slt$ this is recently known to be true.
\end{remar}

\section{Asymptotic solutions of KZ and eigenvectors of
$H_i(z)=\sum_{j\neq i}\frac{\Om^{(ij)}}{z_i-z_j}$}\

\noindent Let $V=V_1\oco V_n$ be a tensor product of $\slt$ modules.
The KZ equation for $u\nup zn\in V$ is
\begin{equation}\label{eq2}
x\ddx{z_i}u=H_i(z)u,\quad i=1,\ldots,n.
\end{equation}
\begin{defin}
Let
\begin{equation}\label{eq3}
u=e^{\frac {S(z)}{\ka}}(f_0(z)+\ka f_1(z)+\ka^2f_2(z)+\cdots)
\end{equation}
for certain functions $S,f_0,\ldots$.
We say that $u$ is an \emph{asymptotic solution of KZ} if the substitution of
$u$ in \eqref{eq2} gives $0$ in the expansion as a formal power series in $\ka$.
\end{defin}
\begin{theor}
Let $u(z)$ be an asymptotic solution of KZ. Then $\forall i$ the vector $f_0(z)$ is an
eigenvector of $H_i(z)$ with eigenvalue $\frac{\ptl S}{\ptl z_i}(z)$.
\end{theor}
\begin{proof}
$$\ka\big[\frac 1\ka\dfx S{z_i}e^{\frac{S(z)}\ka}(f_0+\ka f_1+\cdots)+
	e^{\frac{S(z)}\ka}(\dfx{f_0}{z_i}+\ka\dfx {f_1}{z_i}+\cdots)\big]
	=H_ie^{\frac{S(z)}\ka}(f_0+\ka f_1+\cdots),$$
and then
$$\dfx S{z_i}f_0=H_if_0,\quad\dfx S{z_i}f_1+\dfx{f_0}{z_i}=H_if_1,\ldots$$
\end{proof}
\begin{remar}
With the following powers we can in general recover $f_i$ from $f_0$ $\forall i>0$.
\end{remar}

\subsection{Shapovalov form}
Let $V=V_m$ generated by $v$. There is a unique bilinear symmetric form $S$ on $V$
such that $S(v,v)=1$, $S(hx,y)=S(x,hy)$, $S(ex,y)=S(x,fy)$ $\forall x,y\in V$.
For example $S(fv,fv)=S(v,efv)=S(v,hv)=m$, and in general we get
$$S(f^kv,f^kv)=k!m(m-1)\cdots(m-k+1).$$
The form $S$ is called \emph{Shapovalov form}.

Consider now $V_1\otimes V_2$ and $S=S_1\otimes S_2$ the tensor product of the Shapovalov
forms on $V_1$ and $V_2$, i.e., $S(x_1\otimes x_2,y_1\otimes y_2)=S(x_1,y_1)S(x_2,y_2)$.
We have
\begin{lemma}
$S(\Om(x_1\otimes x_2),y_1\otimes y_2)=S(x_1\otimes x_2,\Om(y_1\otimes y_2))$.
\end{lemma}
\begin{corol}
For $V=V_1\oco V_n$, $S=S_1\oco S_n$, we have $\forall x,y\in V$
$$S(H_i(z)x,y)=S(x,H_i(z)y).$$
\end{corol}
\begin{lemma}
Let $u_1(z)$ be a solution of KZ with parameter $\ka$. Let $u_2(z)$ be a solution
with parameter $-\ka$. Then $S(u_1(z),u_2(z))$ does not depend on $z$.
\end{lemma}
\begin{proof}
$$\ddx{z_i}S(u_1,u_2)=S(\ddx{z_i}u_1,u_2)+S(u_1,\ddx{z_i}u_2)
	=S(\frac 1\ka H_iu_1,u_2)+S(u_1,-\frac 1\ka H_iu_2)=0.$$
\end{proof}
\begin{lemma}
Let $u(z,\ka)=e^{\frac {S(z)}{\ka}}(f_0(z)+\ka f_1(z)+\cdots)$ be an asymptotic
solution of KZ with parameter $\ka$. Then $S(f_0(z),f_0(z))$ does not depend on $z$.
\end{lemma}
\begin{proof}
Let $u_2(z,\ka)=u(z,-\ka)$. It is an asymptotic solution to KZ with parameter $-\ka$, and
then $S(u,u_2)$ does not depend on $z$. But
$$S(e^{\frac{S(z)}\ka}(f_0(z)+\ka f_1(z)+\cdots),e^{-\frac{S(z)}\ka}(f_0(z)-\ka f_1(z)+\cdots))
	=S(f_0(z),f_0(z))+\Oc(\ka).$$
\end{proof}

\begin{reftc}
\cite{FMTV}, \cite{TV3}, \cite{MaV}, \cite{EV1}, \cite{F}, \cite{G}, \cite{RV}, \cite{TV4}
\end{reftc}


\chapter{}
During the last lecture we proved the following result:\\
given
\begin{equation}
 u=e^{\frac{S(z)}{\kappa}}(f_{0}(z)+\kappa f_{1}(z)+\kappa^{2}f_{2}(z)+...) 
\end{equation}
asymptotic solution of the KZ equation, $f_{0}(z)$ is a common eigenvector
for the set of commutating operators $H_i=H_i(z)$, $i=1,...,n$, where:
$$ H_i(z)=\sum_{i\neq j}\frac{\Omega^{(ij)}}{z_i-z_j} $$
with eigenvalues $\frac{\partial S}{\partial z_i}(z)$.
We also introduced the Shapovalov form $S$: given $V=V_m$, generated by $v$, 
$S$ is the unique bilinear symmetric form such that:
\begin{equation}
S(f^kv,f^kv)=k!m(m-1)...(m-k+1)
\end{equation}
and $0$ otherwise, and we proved that if $V=V_1\otimes...\otimes V_n$, $S=S_1\otimes...\otimes S_n$ 
and $u=e^{\frac{S(z)}{\kappa}}(f_{0}(z)+\kappa f_{1}(z)+\kappa^{2}f_{2}(z)+...)$ asymptotic 
solution of the KZ equation, then $S(f_0(z), f_0(z))$ does not depend on $z$.
During this lecture we will discuss how to find explicit asymptotic solutions of KZ equation
and hence eigenvectors for the set of hamiltonian operators $H_i(z)$.\\
Let's consider the space of $\smk$. We know that given:
\begin{equation}
\Phi_{k,n}(t,z,m)=\prod_{1\leq i<j\leq n} (z_i-z_j)^{\frac{m_i m_j}{2}}\prod_{1\leq i<j\leq k} 
(t_i-t_j)^2 \prod_{l=1}^{n}\prod_{i=1}^{k} (t_i-z_l)^{-m_l}
\end{equation}
\begin{equation}
A_{J}(t,z)=\frac{1}{j_1!...j_n!}Sym_{t}\Big[ \prod_{l=1}^{n}\prod_{i=1}^{j_{l}}\frac{1}{t_{j_{1}}+...+t_{j_{l-1+i}}-z_l}\Big],
\end{equation}
\begin{equation}
u^{\gamma}=\Sigma\int_{\gamma} \Phi^{\frac{1}{\kappa}}_{k,n}(t,z)A_{J}(t,z)dt_1\wedge...\wedge dt_k
\end{equation}
is a solution of the KZ equation.\\
We will construct asymptotic solutions taking the limit of 
$u^{\gamma}$ for $\kappa\longrightarrow 0$.\\ 

Let's consider the following examples where we introduce  
the {\it Steepest Descend Method}. 
\begin{examp}
Let's consider the following function:
\begin{equation}
I(\kappa)=\int_{-\infty}^{+\infty} e^{-\frac{s^2}{\kappa}}a(s) ds.
\end{equation}
If $a(0)\neq 0$ we can conclude that:
\begin{equation}
I(\kappa)=\int_{-\infty}^{+\infty} e^{-\frac{s^2}{\kappa}}a(s) ds \sim \sqrt{\pi},
\end{equation}
since:
\begin{equation}
\int_{-\infty}^{+\infty} e^{-\frac{s^2}{\kappa}}a(s) ds=\sqrt{\pi\kappa}\big(a(0)+\kappa\frac{a''(0)}{4}+...\big)
\end{equation}
In the same way we could analyze the function:
\begin{equation}
I(\kappa)=\int_{-\infty}^{+\infty} e^{i\frac{s^2}{\kappa}}a(s) ds
\end{equation}
\end{examp}

\begin{defin}
the function $a(s)$ is called {\it amplitude} while the function 
$\phi(s)=-\frac{s^2}{\kappa}$ is called {\it phase}.
\end{defin}
\begin{examp}
Given a phase function of several variables, $S=S(t_1,...,t_n)$ whose critical points are isolated and non
degenerate, we can reduce the analysis of the behavior of
\begin{equation}
I(\kappa)=\int_{{\RR}^n} e^{-\frac{S(t_1,...,t_n)}{\kappa}}a(t_1,...,t_n) dt_1...dt_n
\end{equation}
to that studied in the previous example. In fact in a neighborhood of any of the critical 
point of $S$, we can find a new set of coordinates $s_1,...,s_n$, such that:
\begin{equation}
S(s_1,...,s_n)=A\pm s_1{^2}\pm s_2{^2}...\pm s_n{^2}
\end{equation}
\end{examp}

Applying this method to:
 
\begin{equation}
\int \Phi^{\frac{1}{\kappa}}_{k,n}(t,z)A_{J}(t,z)dt_1\wedge...\wedge dt_k
\end{equation}
we can conclude that it localizes at the
 critical points of the function $\Phi^{\frac{1}{\kappa}}_{k,n}(t,z)$.\\

Let's define $S=\frac{\log\Phi_{k,n}(t,z)}{\kappa}$,
let's suppose that $z^{o}\in {\CC}^{n}$ has distinct coordinates and let $t^{o}\in {\CC}^{k}$
be a non-degenerate critical point of $\Phi_{k,n}$ w.r.t $t$, i.e 
\begin{equation}\frac{\partial\Phi_{k,n}}{\partial t_i}(t^o, z^o)=0
\end{equation}
 and
\begin{equation}
\det\frac{{\partial}^2\Phi_{k,n}}{\partial t_i\partial t_j}(t^o, z^o)\neq 0.
\end{equation}

\begin{theor}{\rm (Implicit Function.)}
There exist a unique holomorphic ${\CC}^k$ valued function $t=t(z)$,
such that:
\begin{enumerate}
\item $\frac{\partial\Phi_{k,n}}{\partial t_i}(t(z),z)=0\hspace{5pt} i=1,...,k;$\\
\item $t^{o}=t(z^{o}).$
\end{enumerate}
\end{theor}

Under these hypothesis we have the following:

\begin{theor} The KZ equation with values in $\smk$ has an asymptotic solution, $u$,
in a neighborhood of $z^{o}$ given by the following expression:

\begin{equation}
u(z)=\frac{\Phi_{k,n}(t(z),z)^{\frac{1}{\kappa}}}{\sqrt{\det\frac {{\partial}^2}{\partial t_i\partial t_j}
\ln\Phi_{k,n}(t(z),z)}}\big(\sum_{\vert J\vert=k}{\mathcal A}_{J}(t(z),z)f_{J}v+o(\kappa)\big).
\end{equation}
\end{theor}

\begin{corol} If $t^{o}\in {\CC}^k$ is a non degenerate critical point of $\Phi_{k,n}=\Phi_{k,n}(t,z^{o})$ then 
the vector:
\begin{equation}
\omega(t^{o},z^{o})=\sum_{\vert J\vert=k}{\mathcal A}_{J}(t(z),z)f_{J}v
\end{equation}

is an eigenvector of  the set of commutating hamiltonians $H_i=H_i(z)$, $i=1,...,n$.
\end{corol}

\begin{defin} 
$\omega(t^{o},z^{o})$ is called the {\it Bethe} vector.
\end{defin}

We also have the following:

\begin{corol}

$$\frac{S(\omega ((t(z),z),\omega (t(z),z))}{\det\frac{\partial^2}{\partial t_i\partial t_j}
\ln\Phi_{k,n}(t(z),z)}=cost$$

\end{corol}

\begin{quest}
What is this constant?
\end{quest}

\begin{remar} Even if this question makes sense for arbitrary Kac-Moody algebras the answer is known only
for $\slt$; in this case we have:

\begin{theor}
$$ S(\omega (t(z),z),\omega(t(z),z))=\det{\frac{\partial^2}{\partial t_i\partial t_j}}
\ln\Phi_{k,n}(t(z),z)$$
i.e. the constant is $1$.
\end{theor}

\end{remar}

\begin{remar} The set of critical points $t^{o}=(t_1,...,t_k)$ of $\Phi_{k,n}=\Phi_{k,n}(t,z^{o})$ 
is invariant w.r.t the permutation of the coordinates, moreover the Bethe vector corresponding to the point
belonging to the same orbit are equal.
\end{remar}

we will address the following:

\begin{quest} Is it true that:
\begin{enumerate}
\item the number of the orbits of critical points of $\Phi_{k,n}$ equals the dim. of $\smk$?
\item the Bethe vectors form a basis?
\end{enumerate}
\end{quest}

To answer these question we introduce the following general result;
given an arrangement of hyperplanes $\mathcal A=\{H\}$ in ${\RR}^k$, defined
by the equations $f_{H}=0$, set:
\begin{equation}
\Phi^{\lambda}_{\mathcal A}:=\prod_{H\in\mathcal A} f^{\lambda_{H}}_{H}, \lambda\in{\CC}.
\end{equation}

\begin{quest} What is  the number of critical points of the function $\Phi^{\lambda}_{\mathcal A}$ for generic 
$\{\lambda_{H}\}$?
\end{quest}

The equation that gives us the critical points is $\frac{d\Phi}{\Phi}=0$ or in other words:

$$\sum\lambda_{H}\frac{\frac{\partial f_{H}}{\partial t_i}}{f_H}=0,\hspace{6pt} i=1,...,k.$$

\begin{remar} Considering a real arrangement in ${\RR}^k$, and positive $\lambda_{H}$, then each bounded 
component has a critical point. It is also possible to show that each of these critical points is non degenerate.
\end{remar}

Let ${\mathcal A}_R\in{\RR}^k$ be a real arrangement, ${\mathcal A}\subset {\CC}^k$ its complexification
and let's define $M={\CC}^k-{\bigcup}_{H\in\mathcal A}H$. Then for positive $\{\lambda_{H}\}$ it is possible to 
prove the following:

\begin{theor}
$$\sharp\Big(\text{critical points of }\Phi^{\lambda}_{\mathcal A}\text{ in }M \Big)
	= \sharp\Big(\text{bounded components of }{\RR}^k-{\bigcup}_{H\in\mathcal A}H \Big)
	=\vert\chi (M)\vert$$
where  $\chi (M)$ is the Euler characteristic of the open complex manifolds $M$.
\end{theor}

\begin{remar} It was conjectured by {\it Varchenko} that for generic $\{\lambda_{k}\}$ this result holds 
for any complex arrangement ${\mathcal A}\subset {\CC}^{k}$. This conjecture was proved by Orlik and Terao.
\end{remar}

From the theorem follows:

\begin{corol} If the highest weights $m_1,...,m_n$ are negative, i.e. the modules $M_{m_{j}}$ are
 all irreducible\\
$\sharp\Big($ orbits of critical points  of $\Phi_{k,n}$ $\Big)=\dim\smk$.
\end{corol}

\begin{remar} Varchenko and Reshetikhin checked that for generic points $z_{1},...,z_{n}$ and negative
$m_1,...,m_n$ the Bethe vectors form a basis in $\smk$. In this case we have:
$$\dim\smk={k+n-2\choose n-2}=\frac{{\mathcal B}}{k!}$$
where $\mathcal B$ in the number of the bounded components of the arrangement.
\end{remar}

Let's now consider the following:

\begin{examp} Let's consider the Selberg integrand:

$$\tilde\Phi (t;\alpha ,\beta ,\gamma)=\prod_{j=1}^{k}t_{j}^{\alpha}(1-t_{j})^{\beta}\prod_{1\leq i<j\leq k}(t_{i}-t_{k})^{2\gamma}$$
where $\alpha, \beta ,\gamma >0$. This function has $k!$ critical points, one on each bounded domain. The 
symmetric group acts on the set of critical points and we have exactly one orbit.

Let's consider the elementary symmetric functions $\lambda_{1}=\sum_{i=1}^{k}t_{i}$, 
$\lambda_{2}=\sum_{i<j}t_{i}t_{j}$,....,$\lambda_{k}=\prod_{i=1}^{k} t_{i}$. We have the following:
The critical points of the Selberg integrand satisfy the following equations:
\begin{equation}
\frac{\alpha}{t_i}+\frac{\beta}{t_i-1}+\sum_{j\neq i}\frac{2\gamma}{t_i-t_j}=0,
\end{equation}
for $i=1,...,k$. If $(t_1,...,t_k)$ is a critical point then:
\begin{theor}
\begin{equation}
\lambda_{j}={k\choose j}\prod_{i=1}^{j}\frac{\alpha (k-i)\gamma}{\alpha+\beta +(2k-i-1)\gamma},
\end{equation}
for $j=1,...,k$.
\end{theor}
\begin{proof}
The equation for the critical point can be rearranged in the following way:
$$(j+1)(\alpha +j\gamma)\lambda_{i+1} = (k-j)(\alpha -\beta +(k+j-1)\gamma)\lambda_{j}.$$
\end{proof}
\end{examp}

Now let $t\in {\CC}^k$ be a singular point of the function $\tilde\Phi (t;-m_{1},-m_{2} ,1)$ and 
let's consider the Bethe vector:
\begin{equation}
\omega =\sum_{p=0}^{k}f^{k-p}v_{m_1}\otimes f^{p}v_{m_2}A_{k-p,p}(t).
\end{equation}
where:
$$A_{k-p,p}(t)=\sum_{1\leq i_{1}<...<i_{j}\leq k}\frac{1}{(t_{i_1}-1)...(t_{i_p}-1)}\prod_{i\notin (i_1,...,i_p)}\frac{1}{t_i}$$

This vector belongs to $\Sing M_{m_1}\otimes M_{m_2}[m_{1}+m_{2}-2k]$ and w.r.t. the standard
generator
$$\omega_{m{_1}+m{_2}-2k}=\sum_{p=0}^{k}\frac{f^{k-p}v_{m_1}\otimes f^{p}v_{m_2}}{(k-p)!m_{1}(m_1-1)...(m_1-k+p+1)m_{2}(m_2-1)...
(m_2-k+p+1)}$$
of $\Sing M_{m_1}\otimes M_{m_2}[m_{1}+m_{2}-2k]$ we have the following:
\begin{theor}
$$\omega_{Bethe}=k!\prod_{j=1}^{k}(m_{1}+m_{2}-(2k-j-1))\omega_{m{_1}+m{_2}-2k}.$$
\end{theor}
\begin{theor} Let $S=S_{m_1}\otimes S_{m_2}$ the Shapovalov form on $M_{m_1}\otimes M_{m_2}$ and $t^{o}$ a 
critical point of the function $\tilde\Phi (t;-m_{1},-m_{2} ,1)$. Then:
$$S(\omega_{Bethe},\omega_{Bethe})=\det\big(\frac{\partial^{2}}{\partial t_i\partial t_j}\ln\tilde{\Phi}(t)\big)$$
\end{theor}

The theorem follows from the following lemmas:

\begin{lemma}
$$S(\omega_{Bethe},\omega_{Bethe})=k!\prod_{l=0}^{k-1}\frac{(m_{1}+m_{2}-2k+l+2)^3}{(m_{1}-l)(m_{2}-l)}$$
\end{lemma}
\begin{proof}
The result is based on:
$$F(a,b,c;1)=\frac{\Gamma (c)\Gamma (c-a-b)}{\Gamma (c-a)\Gamma (c-b)}.$$
\end{proof}
\begin{lemma} 
$$\det\big(\frac{\partial^{2}}{\partial t_i\partial t_j}\ln\tilde{\Phi}(t^{o},-m_{1},-m_{2},1)\big)=k!\prod_{l=0}^{k-1}\frac{(m_{1}+m_{2}-2k+l+2)^3}{(m_{1}-l)(m_{2}-l)}.$$
\end{lemma}
\begin{proof}
This follows comparing the asymptotic expansion of the Selberg integral:
$$ I(\kappa)=\int_{0\leq t_{k}\leq ...\leq t{1}\leq 1}
	{\tilde\Phi}^{\frac{1}{\kappa}}dt_{1}\wedge...\wedge dt_{k}$$ 
with its expression in term of $\Gamma$-functions.
\end{proof}

The important example for applications is the case of the tensor product of finite
dimensional representation and not for generic Verma modules. Given 
$m_{1},...,m_{n}\in {\ZZ}_{>0}$ let's consider the function:
$$\Phi_{k,n}(t,z,m)=\prod_{1\leq i<j\leq k}(t_i-t_j)^2 \prod_{l=1}^{n}
\prod_{i=1}^{k} (t_i-z_l)^{-m_l}$$

\begin{genpr}
Study the critical set of this function.
\end{genpr}

\begin{remar} It's natural to conjecture that the number of orbits of the
critical points must be related to the dimension of 
$S\smk$.
\end{remar}

\begin{remar} Let $m_{1},...,m_{n}\in {\ZZ}_{>0}$ and let's consider 
$$L^{\otimes m}=L_{m_1}\otimes...\otimes L_{m_n}=\oplus L_{a}$$
where each $a$ has the form:
\begin{enumerate}
\item $a=\vert m\vert -2k$;
\item $a\ge 0$.
\end{enumerate}
Let's denote $\omega (m,k)$ the multiplicity of $L_{\vert m\vert -2k}$ in 
$L^{\otimes m}$.

\begin{examp}
Let take $m=(1,1,1)$, $\vert m\vert =3$, so:
$$L_{1}\otimes L_{1}\otimes L_{1}=L_{3}\oplus L_{1}\oplus L_{1}.$$
In this case we have $\omega (m,0)=1$ and $\omega (m,1)=2$.
\end{examp}

Note that $\dim M^{\otimes m}[\vert m\vert -2k]=\omega (m,k)$ and 
each direct summand $L_{a}$ contains one singular vector of weight $a$.
\end{remar}

Now we are ready to describe the critical set of the function 
$$\Phi_{k,n}(t,z,m)=\prod_{1\leq i<j\leq k}(t_i-t_j)^2 \prod_{l=1}^{n}
\prod_{i=1}^{k} (t_i-z_l)^{-m_l}$$
for $m_{1},...,m_{n}\in {\ZZ}_{>0}$.
Let  $\lambda_{1}=\sum_{i=1}^{k}t_{i}$, 
$\lambda_{2}=\sum_{i<j1}t_{i}t_{j}$,....,$\lambda_{k}=\prod_{i=1}^{k} t_{i}$ the standard symmetric functions
of $t_1,...,t_k$ and denote ${\CC}^{k}_{\lambda}$ the space with coordinates $\lambda_{1},...,\lambda_{k}$.

\begin{theor} Let $m_{i}$, $i+1,...,n$, as before and $k\in {\ZZ}_{>0}$.
\begin{enumerate}
\item If $\vert m\vert +1-k>k$, then for generic $z$ all critical points of $\Phi_{k,n}(t)$ are non degenerate and
the critical set consists of $\omega (m,k)$ orbits.
\item If  $\vert m\vert +1-k=k$, then for any $z$ the function $\Phi_{k,n}(t)$ has no critical points.
\item If $0\leq\vert m\vert +1-k<k$ then for generic $z$ the function $\Phi_{k,n}(t)$ may have only non isolated
critical points. Written in symmetric coordinates $\lambda_{1},...,\lambda_{k}$, the critical set consists of
$\omega (m, \vert m\vert +1-k)$ straight lines in the space ${\CC}^{k}_{\lambda}$.
\item If $\vert m\vert +1-k<0$, then for any $z$ the function  $\Phi_{k,n}(t)$ has no critical points.
\end{enumerate}
\end{theor}
\begin{examp} Let $m=(1,1,1)$, $\vert m\vert =3$. In this case we have:
$$\prod_{1\leq i<j\leq k}(t_i-t_j)^2 \prod_{l=1}^{3}\prod_{i=1}^{k} (t_i-z_l)^{-1}.$$
If:
\begin{enumerate}
\item $k=1$, $\vert m\vert +1-k=3$ so $\Phi_{1,3}(t)$ has two non degenerate critical points;
\item $k=2$, $\Phi_{2,3}(t)$ has no critical points;
\item $k=3$, $\Phi_{3,3}(t)$ has two straight lines of critical points;
\item $k=4$, $\Phi_{4,3}(t)$ has one straight line of critical points;
\item $k\ge 5$, $\Phi_{k,3}(t)$ has no critical points.
\end{enumerate}
\end{examp}

\begin{reftc}
\cite{F}, \cite{G}, \cite{RV}, \cite{OT}, \cite{V1}, \cite{MV}, \cite{TV4}, \cite{ScV}
\end{reftc}


\chapter{}

During the last lecture we studied the critical points of  the function:

$$\Phi_{k,n}(t,z,m)=\prod_{1\leq i<j\leq k}(t_i-t_j)^2 \prod_{l=1}^{n}
\prod_{i=1}^{k} (t_i-z_l)^{-m_l}$$

We studied the action on the set of the critical points of the symmetric group and we discussed the following:
\begin{theor}
$$\dim\smk={k+n-2\choose n-2}.$$
\end{theor}

and also:

\begin{theor} Let $m_{i}\in {\ZZ}_{>0}$ for $i=1,...,n$, $k\in {\ZZ}_{>0}$ and let
$\omega(m=(m_{1},...,m_{n}), k)$ be the multiplicity of $L_{\vert m\vert -2k}$ in
$L_{m_1}\otimes ...\otimes L_{m_n}$. Then:
\begin{enumerate}
\item If $\vert m\vert +1-k>k$, then for generic $z$ all critical points of $\Phi_{k,n}(t)$ are non degenerate and
the critical set consists of $\omega (m,k)$ orbits.
\item If  $\vert m\vert +1-k=k$, then for any $z$ the function $\Phi_{k,n}(t)$ has no critical points.
\item If $0\leq\vert m\vert +1-k<k$ then for generic $z$ the function $\Phi_{k,n}(t)$ may have only non isolated
critical points. Written in symmetric coordinates $\lambda_{1},...,\lambda_{k}$, the critical set consists of
$\omega (m, \vert m\vert +1-k)$ straight lines in the space ${\CC}^{k}_{\lambda}$.
\item If $\vert m\vert +1-k<0$, then for any $z$ the function  $\Phi_{k,n}(t)$ has no critical points.
\end{enumerate}
\end{theor}

Today we will discuss the relation of the previous theorem (in particular the statements 2., 3., and 4.) with
the theory of the Fuchsian differential equations.

\section{Differential Equations with Regular Singular Points}

We will focus to the case of second order equations:

\begin{equation}
u''+p(z)u'+q(z)u=0 
\end{equation}

Let's start with the following:

\begin{defin} A point $z_{o}\in {\CC}$ is called {\rm ordinary} for the differential equation if the
functions $p(z)$ and $q(z)$ are holomorphic at a neighborhood of this point.
\end{defin}

To study the case $z=\infty$ put $z=\frac{1}{\xi}$:\\

$$\frac{du}{dz}=\frac{d\xi}{dz}\frac{du}{d\xi}=-{\xi}^2\frac{du}{d\xi}$$ and
$$\frac{d^2u}{d\xi^2}=\xi^4\frac{d^2u}{d\xi^2}+2\xi^3\frac{du}{d\xi}.$$

With this change of variables the equation turns into:
\begin{equation}
\frac{d^2u}{d\xi^2}+\Big[\frac{2}{\xi}-\frac{1}{\xi^2} p\big(\frac{1}{\xi}\big)\Big]
\frac{du}{d\xi}+\frac{1}{\xi^4}q\big(\frac{1}{\xi}\big)u.
\end{equation}

So we have the following:

\begin{defin} The point $z=\infty$ is an ordinary point of the equation if $\xi =0$ 
is an ordinary point of the transformed equation, i.e if:
$$\Big[\frac{2}{\xi}-\frac{1}{\xi^2} p\big(\frac{1}{\xi}\big)\Big]$$
and
$$\frac{1}{\xi^4}q\big(\frac{1}{\xi}\big)$$
are holomorphic functions at $\xi =0$.
\end{defin}

\begin{defin} If $z_{o}$ is called {\rm singular} point for the equation if it is not ordinary.
\end{defin}
\begin{defin} The point $z_{o}$ is called {\rm regular singular} if the following conditions are
 verified:\\
 1) $z_{o}$ is singular;\\
 2) at $z_{o}$ the function $p(z)$ at has a pole of order $\leq 1$;\\
 3) at $z_{o}$ the function $q(z)$ has a pole of order $\leq 2$.
 \end{defin}
 
 \begin{remar} We have analogous definitions for the point $z=\infty$ and local coordinate $z=\frac{1}{\xi}$.
 \end{remar}
 
 In the following two section we will discuss the solutions of the equation at an ordinary point and at 
 a singular point.
 
\section{Solutions at a neighborhood of an ordinary point}
 
 Let $z_{o}$ an ordinary point, so we can write:
 
$$p(z)=\sum_{n=0}^{\infty} b_{n}(z-z_{o})^n$$ 
and 
$$q(z)=\sum_{n=0}^{\infty} d_{n}(z-z_{o})^n$$ 
 
for $\vert z-z_{o}\vert <R$. We look for solutions of the form:
$$u(z)=\sum_{n=0}^{\infty} a_{n}(z-z_{o})^n.$$
 
Under these hypothesis we have the following:

\begin{theor}
For any choice of $a_{0}, a_{1}$ there exist a unique series 
$$u(z)=\sum_{n=0}^{\infty} a_{n}(z-z_{o})$$
with a non zero radius of convergence, satisfying the differential equation:
$$u''+p(z)u'+q(z)u=0.$$
\end{theor}
\begin{proof}
Plugging $u(z)=\sum_{n=0}^{\infty} a_{n}(z-z_{o})$ into the equation 
$u''+p(z)u'+q(z)u=0$ we get:\\
for $n=0$, $0a_{0}=0$, for $n=1$, $0a_{1}+0a_{0}b_{0}=0$, while for $n\geq 2$,
$$n(n-1)a_{n}+\sum_{k=0}^{n-1} ka_{k}b_{n-k-1} +\sum_{k=0}^{n-2} a_{k}b_{n-k-2}.$$
\end{proof}

\section{Solutions at a neighborhood of a regular singular point}
 
 Let's assume $z_{o}$ be a regular singular point so that the coefficients of the equation can be written
 in the following way:
$$p(z)=\sum_{n=0}^{\infty} b_{n}(z-z_{o})^{n-1}$$ 
and 
$$q(z)=\sum_{n=0}^{\infty} d_{n}(z-z_{o})^{n-2},$$

for $\vert z-z_{o}\vert <R$. In this case we look for solutions of the form:
$$u(z)=\sum_{n=0}^{\infty} a_{n}(z-z_{o})^{n+c}$$
where $a_{0}\neq 0$ and $c$ is a constant that has to be determined.
Plugging the function $u(z)=\sum_{n=0}^{\infty} a_{n}(z-z_{o})^{n+c}$ into the equation, we get the following
equation:
$$\sum_{n=0}^{\infty}(n+c)(n+c-1)a_{n}(z-z_{0})^{n+c-2}+\Big[\sum_{n=0}^{\infty}p_{n}(z-z_{0})^{n-1}\Big]
\Big[\sum_{n=0}^{\infty}(n+c)a_{n}(z-z_{0})^{n+c-1}\Big]+$$
$$+\Big[\sum_{n=0}^{\infty}q_{n}(z-z_{0})^{n-2}\Big]
\Big[\sum_{n=0}^{\infty}a_{n}(z-z_{0})^{n+c}\Big]=0.$$
From this equation we get the following recurrence relations:\\
for $n=0$: 
\begin{equation}
[c(c-1)+cp_{0}+q_{0}]a_{0}=0,
\end{equation}
while for $n\geq 1$:
\begin{equation}
[(n+c)(n+c-1)+(n+c)p_{0}+q_{0}]a_{n}=-\sum_{k=0}^{n-1}a_{k}[(k+c)p_{n-k}+q_{n_k}].
\end{equation}
Since $a_{0}\neq 0$, we have:
\begin{equation}
c^2+(p_{0}-1)+q_{0}=0.
\end{equation}

\begin{defin}
The equation $c^2+(p_{0}-1)+q_{0}=0$ is called {\rm indicial} equation, while the roots $c_{1}$ and $c_2$
are called the {\rm exponents} of the equation at the point $z=z_{0}$.
\end{defin}

Let's write $F(c)=c^2+(p_{0}-1)+q_{0}$ so that we have:
$$F(n+c)a_n=-\sum_{k=0}^{n-1}a_{k}[(k+c)p_{n-k}+q_{n_k}]$$
for $n\geq 1$, and $F(c)=0$, for the indicial equation. 
Let's order the roots $c_1, c_2$ of the polynomial $F(c)$ in such a way $\Re c_1\geq\Re c_2$ and set 
$s=c_1-c_2$. Now if we start with $c=c_1$, we get:
$$F(n+c_1)=n(n+c_1-c_2)=n(n+s)$$
from which we see that the recurrence equations are satisfied; so we can state the following:
\begin{theor} If $z_{0}$ is a regular singular point with exponentials $c_1$ and $c_2$ such that 
$\Re c_1\geq\Re c_2$, then the series:
\begin{equation}
u_{1}(z)=(z-z_{0})^{c_1}+\sum_{n=1}^{\infty} a_{n}(z-z_{0})^{n+c_{1}}
\end{equation}
satisfies the equation and it is convergent in a region of the form $0<\vert z-z_{0}\vert <r$.
\end{theor} 
\begin{quest} How to find the second solution?
\end{quest}
The answer to this question depends on $s=c_{1}-c_{2}$. Put $c=c_2$, then the equation for $n=0$ is 
satisfied; for $n\geq 1$ we get:
$$n(n-s)a_n =\sum_{k=0}^{n-1}q_{k}\big[(k+c_2)p_{n-k}+q_{n-k}\big]:$$
\begin{theor} If $s\notin {\ZZ}$ then there exists a series:
\begin{equation}
u_{2}(z)=(z-z_{0})^{c_2}+\sum_{n=1}^{\infty} a_{n}(z-z_{0})^{n+c_{1}}
\end{equation}
that satisfies the equation and converges in a region of the form $0<\vert z-z_{0}\vert <r.$
\end{theor}
In the case $s=0$ we don't get any new solution, while if $s\in {\ZZ_{>0}}$, for $n=s$, we get:
$$0=\sum_{k=0}^{s-1}q_{k}\big[(k+c_2)p_{s-k}+q_{s-k}\big].$$
This gives a  non trivial constraint: given arbitrary $a_0$ and $a_s$ we can get:
$$u_{1}(z)=(z-z_{0})^{c_1}+\sum_{n=1}^{\infty} a_{n}(z-z_{0})^{n+c_{1}}$$
and 
$$u_{2}(z)=(z-z_{0})^{c_2}+\sum_{n=1}^{\infty} a_{n}(z-z_{0})^{n+c_{2}}$$
if $a_0$ and $a_s$ automatically satisfy the recurrence relation or:
$$u_{2}(z)=u_{1}\ln(z-z_0)+\sum_{n=0}^{\infty}d_{n}(z-z_{0})^{n+c_2}$$
if not.

\section{Fuchsian Equations}

Let's start with the following:
\begin{defin}{\rm [Lazarus Fuchs (1833-1902)]}
A linear differential equation is called {\rm fuchsian} if all its singular points are regular singular.
\end{defin}

\begin{remar} During this section we will focus only on second order differential equations.

\end{remar}

\begin{remar} If the equation:
$$u''+p(z)u'+q(z)u=0$$
is fuchsian, then the coefficients $q(z)$ and $p(z)$ are rational functions.
\end{remar}

We also have the following:
\begin{propo}
A fuchsian equation has at least a singular point.
\end{propo}

In what follows we will give a description of the second order fuchsian differential equations with one
two and three singular points. Let's start with the first case. Let's assume that $u''+p(z)u'+q(z)u=0$ 
has only one singular point at $z=z_1$, then:
\begin{propo} The equation has the form:
\begin{equation}
u''+\frac{2}{z-z_1} u'=0.
\end{equation}
and its solutions are of the form $u(z)=a+\frac{d}{z-z_1}$ for given $a$ and $b$.
\end{propo}
\begin{proof}
From the hypothesis we have that:
$$p(z)=\frac{p_{1}(z)}{z-z_1}$$
and 
$$q(z)=\frac{q_{1}(z)}{(z-z_1)^2}$$
where $p_{1}(z)$ and $q_{1}(z)$ are polynomials. Since $z=\infty$ is an ordinary point we have:
$$2z-z^2 p(z)=2z-\frac{z^2 p_{1}(z)}{z-z_1}=\frac{2z^2-z^2p_{1}-2z_{1}z}{z-z_1}$$
and 
$$z^4q(z)=\frac{z^4 q_{1}(z)}{(z-z_1)^2}$$
are holomorphic at $z=\infty$. From the last equation we deduce  that $q_1\equiv 0$
(otherwise deg$z^4 q_{1}(z)\geq 4$) while from the first one follows that $p_{1}=2$. Clearly 
$u_{1}(z)=1$ and $u_{2}(z)=\frac{1}{z-z_{1}}$ are solutions.
\end{proof}

We also have the following result:

\begin{theor} If $z=\infty$ is the only singular point of a fuchsian differential equation, then the equation is:
\begin{equation}
u''=0.
\end{equation}
\end{theor}
\begin{proof}
In this case the coefficients $p(z)$ and $q(z)$ are two polynomials. So we have that the function
$2z-z^2p(z)$ has poles of order $\leq 1$ while the function $z^4q(z)$ has poles of order $\leq 2$. From this 
observation follows that $p(z)=q(z)=0$.
\end{proof}

Now let's discuss the case of fuchsian equations with two singular points $z_1$ and $z_2$. In this case 
we have two different cases:
\begin{enumerate}
\item  $z_2=\infty$, or:
\item  $z_1$ and $z_2$ $\in {\CC}$.
\end{enumerate}
In the first case we have:
\begin{propo} The fuchsian equation can be written as:
\begin{equation}
u''+\frac{k_1}{z-z_1}u'+\frac{k_2}{(z-z_1)^2}u=0,
\end{equation}
with $k_1, k_2\in {\CC}$. Moreover the solutions are:
\begin{equation}
u(z)=A(z-z_1)^{c_1}+B(z-z_1)^{c_2}
\end{equation}
if $c_1\neq c_2$, or:
\begin{equation}
u(z)=(z-z_1)^{c}\big(A+B\ln(z-z_1)\big)
\end{equation}
if $c_1=c_2=c$.
\end{propo}
For the second case, when both $z_1$ and $z_2$ belong to $\CC$, it is possible to prove the following:
\begin{propo}
The fuchsian equation has the form:
\begin{equation}
u''+\frac{2z+k_3}{(z-z_1)(z-z_2)}u'+\frac{k_4}{(z-z_1)^2(z-z_2)^2}u=0.
\end{equation}
\end{propo}
Let's now start to discuss the case of three singular point.
\begin{remar} We have the following general:
\begin{theor}
A second order fuchsian equation with $n+1$ singular points $z_1,...,z_n,\infty$ has the form:
\begin{equation}
u''+\frac{T_{n-1}}{(z-z_1)...(z-z_n)}u'+\frac{T_{2n-2}}{(z-z_1)^2...(z-z_n)^2}u=0,
\end{equation}
where $T_{n-1}$ and $T_{2n-2}$ are polynomials of degree less or equal to $n-1$ and $2n-2$.
\end{theor}
\end{remar}

If we call $\alpha_{1,k}$ and $\alpha_{2,k}$ the exponents of the equation at the point $z=z_k$ and 
$\alpha_{1,\infty}$ and $\alpha_{2,\infty}$ the exponents at the point $z=\infty$, then:
\begin{theor}{\rm [Fuchsian Invariants.]} 
\begin{equation}
\alpha_{1,\infty}+\alpha_{2,\infty}+\sum_{k=1}^{n}(\alpha_{1,k}+\alpha_{2,k})=n-1.
\end{equation}
\end{theor}
Let's focus on the case of three singular points $z_1$, $z_2$ and $\infty$. In this case the exponents are:
$(\alpha_{1,\infty},\alpha_{2,\infty})$ $(\alpha_{1,1}, \alpha_{2,1})$ $(\alpha_{1,2}, \alpha_{2,2})$ and 
we have:
\begin{equation}
\alpha_{1,\infty}+\alpha_{2,\infty}+\sum_{k=1}^{2}(\alpha_{1,k}+\alpha_{2,k})=1.
\end{equation}

\begin{claim} the fuchsian equation can be written in the following way:
$$u''+\Big[\frac{1-\alpha_{1,1}-\alpha_{2,1}}{z-z_1}+\frac{1-\alpha_{1,2}-\alpha_{2,2}}{z-z_2}\Big]u'+$$
$$\Big[\frac{\alpha_{1,1}\alpha_{2,1}}{(z-z_1)^2}+
\frac{\alpha_{2,1}\alpha_{2,2}}{(z-z_2)^2}+
\frac{\alpha_{1,\infty}\alpha_{2,\infty}-\alpha_{1,1}\alpha_{2,1}-\alpha_{1,2}\alpha_{2,2}}{(z-z_1)(z-z_2)}
\Big]u=0$$
\end{claim}

\begin{remar} The previous equation is the form of the general Fuchsian equation of the second order and with
three regular singular points written in terms of the exponents.
\end{remar}

Using the notation $z=a, b, c$ for the singular points (all of them different by $\infty$)
and $(a',a'')$, $(b',b'')$, $(c',c'')$ for the correspondents exponents, we have: 

\begin{theor} The previous equation can be written as:
$$u''+\Big[\frac{1-a'-a''}{z-a}+\frac{1-b'-b''}{z-b}+\frac{1-c'-c''}{z-c}\Big]u'+$$
$$\Big[\frac{a'a''(a-b)(a-c)}{z-c}+\frac{b'b''(b-a)(b-c)}{z-b}+\frac{c'c''(c-a)(c-b)}{z-c}\Big]
\frac{u}{(z-a)(z-b)(z-c)}=0.$$
\end{theor}

\begin{remar} This equation is called {\rm Riemann-Papperitz equation.}
\end{remar}

\begin{defin}{\rm [Riemann P-symbol.]}

Let's write:

\begin{equation}
u=P\left (
\begin {array}{cccc}
a   & b   & c       \\
a'  & b'  & c' & z  \\
a'' & b'' & c''     
\end{array}
\right )
\end{equation}

This is called the {\rm Riemann P-symbols}. This denotes a solution of a Fuchsian differential equation
with singular points at $a,b,c$ (distinct but can be $\infty$) and exponentials $(a',a'')$, $(b',b'')$, 
$(c',c'')$.
\end{defin}

\begin{remar} Under these hypothesis we have:
$$a'+a''+b'+b''+c'+c''=1.$$
\end{remar}
\begin{remar} The informations given in the P-symbol determines the equation uniquely.
\end{remar}

\begin{theor} Let's suppose $a, c\in {\CC}$ and:
\begin{equation}
u=P\left (
\begin {array}{cccc}
a   & b   & c       \\
a'  & b'  & c' & z  \\
a'' & b'' & c''     
\end{array}
\right ).
\end{equation}
If 
$$u=\Big(\frac{z-a}{z-c}\Big)^kw$$
then:

\begin{equation}
u=P\left (
\begin {array}{cccc}
a   & b   & c       \\
a'-k  & b'  & c'+k & z  \\
a''-k & b'' & c''+k     
\end{array}
\right ).
\end{equation}
\end{theor}

\begin{theor} If:
\begin{equation}
u=P\left (
\begin {array}{cccc}
a   & b   & \infty       \\
a'  & b'  & c' & z  \\
a'' & b'' & c''     
\end{array}
\right )
\end{equation}
and $u=(z-a)^k$ then:

\begin{equation}
u=P\left (
\begin {array}{cccc}
a   & b   & \infty       \\
a'-k  & b'  & c'+k & z  \\
a''-k & b'' & c''+k     
\end{array}
\right ).
\end{equation}
\end{theor}

\begin{theor} If 
\begin{equation}
u=P\left (
\begin {array}{cccc}
a   & b   & c       \\
a'  & b'  & c' & z  \\
a'' & b'' & c''     
\end{array}
\right )
\end{equation}
the fractional transformation:
$$v=\frac{Az+B}{Cz+D}$$
with $AD-CB\neq 0$
maps $z=a, b, c$ into $v=a_{1}, b_{1}, c_{1}$ and the solution $u$ can be written
in terms of $v$ in the following way:
\begin{equation}
u=P\left (
\begin {array}{cccc}
a_{1}   & b_{1}   & c_{1}      \\
a'  & b'  & c' & v  \\
a'' & b'' & c''     
\end{array}
\right )
\end{equation}
\end{theor}

We also have the following general:
\begin{lemma} Any linear fractional transformation:
$$v=\frac{Az+B}{Cz+D}$$
with $AD-BC\neq 0$ is the composition of elementary transformations:\\
$v=\gamma z$, $v=\gamma +z$ and $v=\frac{1}{z}$.
\end{lemma}

We have the following important:
\begin{theor} Any P-symbols can be reduced to a P-symbols of the form:
\begin{equation}
y=P\left (
\begin {array}{cccc}
0   & 1   & \infty       \\
0  &  0 & \alpha & x \\
1-\gamma & \gamma -\alpha -\beta & \beta    
\end{array}
\right )
\end{equation}
\end{theor}

\begin{remar} The corresponding differential equation is the following:

\begin{equation}
\frac{d^2 y}{dx^2}+\Big[\frac{1-(1-\gamma)}{x}+\frac{1-(\gamma -\alpha -\beta)}
{x-1}\Big]\frac{dy}{dx}+\frac{\alpha \beta}{x(x-1)}y=0.
\end{equation}
Multiplying both sides by $x(x-1)$ we get:
\begin{equation}
x(x-1)y''(x)+[\gamma -(\alpha +\beta +1)x]y'-\alpha\beta y=0,
\end{equation}
i.e. the {\rm Gauss} differential equation.
\end{remar}

Let's start the  discussion the connection of the theory of the Fuchsian differential equations with the
critical points of the master function.
Let's consider the following:
\begin{equation}
F(x)u''(x)+G(x)u'(x)+H(x)u(x)=0
\end{equation}
where $F$, $G$ and $H$ are polynomials of degree $n$, $n-1$ and $n-2$ respectively. 
\begin{claim} If the polynomial $F$ has no multiple roots then the previous equation is Fuchsian.
\end{claim}
In this case let's write:
$$F(x)=\prod_{j=1}^{n}(x-z_j)$$
$$\frac{G(x)}{F(x)}=\sum_{j=1}^{n}-\frac{m_j}{x-z_j}$$
for suitable complex numbers $m_j$ and $z_j$.
\begin{claim} It is possible to prove that $0$ and $m_j+1$ are exponents of the above equation and if $-k$ is 
one of the exponents at $z=\infty$ then the other is $k-l(m)-1$.
\end{claim}
\begin{probl}
Given polynomials $F(x)$ and  $G(x)$ as above:
\begin{enumerate}
\item find a polynomial $H(x)$ of degree at most $n-2$ such that the equation has a polynomial solution of a 
preassigned degree $k$;
\item find the number of solution of 1.
\end{enumerate}
\end{probl}

we have the following classical results (reference: Szego G., ``Orthogonal Polynomials"):
\begin{theor}
\begin{enumerate}
\item Let $u(x)$ a polynomial solution of degree $k$ of the equation, with roots $t_{1}^{0},...,t_{k}^{0}$ of
multiplicity one. Then $t^{0}=(t_{1}^{0},...,t_{k}^{0})$ is a critical point of the function 
$\Phi_{k,n}(t;z,m)$, where $z=(z_{1},...,z_{n})$ and $m=(m_{1},...,m_{n})$.\\
\item Let $t^{0}$ be a critical point of the function $\Phi_{k,n}(t;z,m)$, then the polynomial 
$u(x)=(x-t_{1}^{0})...(x-t_{k}^{0})$ of degree $k$ is a solution of the equation with 
$$H(x)=\frac{\big(-F(x)u''(x)-G(x)u'(x)\big)}{u(x)}$$
being a polynomial of degree at most $n-2$.
\end{enumerate}
\end{theor}
\begin{remar} Recall that the critical points of the function $\Phi_{k,n}(t;z,m)$ are given by:
\begin{equation}
\sum_{j=1}^{n}\frac{-m_j}{t_i-z_j}+\sum_{j\neq i}\frac{2}{t_i-t_j}, i=1,...,k.
\end{equation}
\end{remar}
\begin{remar} This can be seen as an equivalent version of the Bethe ansatz problem.
\end{remar}
\begin{proof}(of part (1))
Let $u(x)=(x-t_{1}^{0})...(x-t_{k}^{0})$ be a solution of the equation:
$$Fu''+Gu'+Hu=0,$$
and set $x=t_{i}^{0}$; then:
$$F(t_{i}^{0})u''(t_{i}^{0})+G(t_{i}^{0})u'(t_{i}^{0})+H(t_{i}^{0})u(t_{i}^{0})=0.$$
From this we can deduce:
$$\frac{u''(t_{i}^{0})}{u'(t_{i}^{0})}+\frac{G(t_{i}^{0})}{F(t_{i}^{0})}=0,$$
but by definition we have:
$$\frac{G(t_{i}^{0})}{F(t_{i}^{0})}=\sum\frac{-m_i}{t_{i}^{0}-z_i}.$$
To conclude we need the following:

\begin{lemma} If  $u(x)=(x-t_{1}^{0})...(x-t_{k}^{0})$, then:
\begin{equation}
\frac{u''(t_{i}^{0})}{u'(t_{i}^{0})}=\sum_{j\neq i}\frac{2}{t_{i}^{0}-t_{j}^{0}}.
\end{equation}
\end{lemma}
\end{proof}

\begin{examp} If:

$$ u(x)=(x-t_{1}^{0})(x-t_{2}^{0})$$
then:
$$ u(x)=(x-t_{1}^{0})+(x-t_{2}^{0}),$$

$$u''=2$$
so:
$$\frac{u''}{u'}(t_1^{0})=\frac{2}{t_{1}^{0}-t_{2}^{0}}$$
\end{examp}

This give us part i. of the previous theorem.\\

\begin{remar} The critical points label solutions of the problem 2., and a critical point defines a 
differential equation and a polynomial solution of that equation.
\end{remar}

\begin{remar} In 19-th century Heine and Stieltjes showed that for fixed real $z_1,...,z_n$ and negative
$m_1,...,m_n$ the problem 2. has:
$${k+n-2\choose n-2}$$
solutions (reference: {\rm Szego G., ``Orthogonal Polynomials"}).
\end{remar}

Let $m_1,...,m_n\in {\ZZ_{>0}}$ and $z\in {\ZZ_{>0}}$, such that $\vert m\vert -2k\geq 0$, then
$k<\vert m\vert -k+1$. Assume that $t^{0}$ is a critical point of $\Phi_{k,n}(t;z,m)$, then it is possible to
construct a line of critical points of the function $\Phi_{\vert m\vert -k+1,n}(t;z,m)$. 

\begin{examp} Let $m=(1,1,1)$ and
$$\Phi_{1,3}(t;z,m)=(t-z_1)^{-1}(t-z_2)^{-1}(t-z_3)^{-1}.$$
it is possible to construct a line of critical points of the function:
Given the two critical points of the function $\Phi_{1,3}(t;z,m)$ we get two lines of critical points of the function.
$$\Phi_{3,3}(t;z,m)=\prod_{i=1}^{3}\prod_{l=1}^{3}(t_{i}-z_l)^{-1}\prod_{1\leq i<j\leq 3}(t_{i}-t_{j})^{2}.$$
Moreover it is possible  to show that for generic $z_1, z_2, z_3$ the points
on those two lines are the only critical points.
\end{examp}

\section{Enumerative Algebraic Geometry}

It is possible to reformulate the main theorem in terms of Enumerative Algebraic Geometry. Let $V$ a two 
dimensional space of polynomials of one variable with coefficients in $\CC$. Let $k_1$ the degree of a typical
polynomial and $k_2$ the degree of a special polynomial $(k_1<k_2)$

\begin{defin}
Given two functions $f(x)$ and $g(x)$ let's denote with:
\begin{equation}
W(f,g)=\Big\vert
\begin {array}{cccc}
f   & g  \\
f'  &  g' \\
\end{array}
\Big\vert
\end{equation}
their Wronskian.
\end{defin}

\begin{remar} If $(f,g)$ and $(\tilde {f},\tilde {g})$ are two basis of $V$ then there exits a constant $c$ 
such that:
$$W(f,g)=cW(\tilde {f},\tilde {g}).$$
\end{remar}

\begin{lemma}
If $deg f=k_1$ and $deg g=k_2$ then $degW(f,g)=k_1+k_2-1$.
\end{lemma}
\begin{proof}
This follows from direct calculation:
\begin{equation}
W(f,g)=\Big\vert
\begin {array}{cccc}
x^{k_1}+..   & x^{k_2}+..  \\
k_{1}x^{k_1-1}+..  & k_{2}x^{k_2-1}+..  \\
\end{array}
\Big\vert=(k_{2}-k_{1})x^{k_1+k_2-1}+...
\end{equation}
\end{proof}

Let $W_{V}=x^{k_1 +k_2 -1}+...$ the Wronskian of $V$ and let's suppose that $W_{V}(x)=\prod_{l=1}^{n}(x-z_l)^{m_l}$,
where $m_l\geq 1$.
\begin{defin}
We will say that $V$ is not singular if $\forall x_{0}\in {\CC}$ there is $f\in V$ such that $f(x_{0})\neq 0$.
\end{defin}

Now we can state the following:
\begin{probl} Assume that $k_{1}<k_{2}$ and $W_{V}(x)=\prod_{l=1}^{n}(x-z_l)^{m_l}$ are fixed.
What is the number of non singular vector spaces $V$ of polynomials with such data?
\end{probl}

It is possible to get an answer to question in terms of representation of $\slt$:
\begin{theor}
For generic $z_1,...,z_n$ the number of non singular vector spaces with such data is equal to:\\
mult$L_{k_{1}-k_{1}-1}\subset L_{m_1}\otimes... \otimes L_{m_n}.$
\end{theor}

\begin{reftc}
\cite{R}, \cite{ScV}
\end{reftc}


\chapter{}
During this lecture we will study the elliptic generalization of the hypergeometric functions. Motivations
come from CFT (Conformal Field Theory) where to any Riemann surfaces $\Sigma$ with marked points decorated with
finite dimensional representations of a simple Lie algebra is assigned the (finite dimensional) vector space
of {\it conformal blocks}. The set of all conformal blocks forms a vector bundle over the moduli space of the
marked Riemann surfaces, this vector bundle carries a projectively flat connection.
We will discuss in detail the case of genus 1:
$$\Sigma =T={\CC}/{\ZZ}+\tau {\ZZ},\hspace{10pt} Im\tau >0.$$

{\it Notations:} we will consider the Lie algebra $\mathfrak{g}$=$\mathfrak{sl}_{2}$,
the Cartan subalgebra $\mathfrak{h}$=${\CC}h$, the module:
$V=V_{1}\otimes...\otimes V_{n},$
and we define: $V[0]=\{v\in V\vert hv=0\}.$ The conformal block will be given by the function:
$$u=u(z_{1},...,z_{n},\lambda,\tau)\in V[0],\hspace{10pt}\lambda ,z_{i}\in{\CC},\hspace{5pt}i=1,..,n
\hspace{5pt}\tau\in{\HH},$$
i.e will be a function of $\lambda$ and will depend on the parameter $z_{1},...,z_{n},\tau$.\\
We will discuss the following equations:
\begin{equation}
1.\hspace{5pt}\kappa\frac{\partial u}{\partial z_{i}}
	=\sum_{j\neq i}r^{(i,j)}(z_{i}-z_{j}, \lambda ,\tau)
		u-h^{(i)}\frac{\partial u}{\partial\lambda}, 
	\hspace{5pt}i=1,...,n;
\end{equation}

\begin{equation}
2.\hspace{5pt} 4\pi i\kappa\frac{\partial u}{\partial\tau}
	=\frac{\partial^{2}u}{\partial\lambda^{2}}+
		\sum_{i,j}H^{(i,j)}(z_{i}-z_{j},\lambda,\tau)u,
	\hspace{5pt}i=1,...n.\label{ciccio3}
\end{equation}

\begin{remar} The equation \ref{ciccio3} is called {\it KZB heat equation}.
\end{remar}

\begin{examp} For $n=1$ we have $V=V_{1}=L_{2p}$ $(\dim L_{2p}=2p+1)$ and 
$V[0]=L_{2p}[0]={\CC}f^{p}v_{2p}$. In this case the function $u=u(\lambda ,z, \tau)$ is a scalar 
valued function
and we have that:\\

\begin{align}
1.\hspace{25pt}\frac{\partial u}{\partial z_{1}} &= 0,\\
2.\hspace{5pt}2\pi i\kappa\frac{\partial u}{\partial\tau} &= \frac{\partial^{2}u}{\partial\lambda^{2}}
+p(p+1)\rho'(\lambda ,\tau)u,\label{ciccio4}
\end{align}
\end{examp}
where:
$$\rho(\lambda ,\tau)=\frac{\theta'(\lambda ,\tau)}{\theta(\lambda ,\tau)}$$
and 
$$\theta(\lambda ,\tau)=-\sum_{j\in {\ZZ}}e^{\pi i(j+\frac{1}{2})^{2}\tau+2\pi i(j+\frac{1}{2})(\lambda +\frac{1}{2})}$$

\begin{remar} The derivative in the equation \ref{ciccio4} is taken w.r.t $\lambda$, i.e  
$('=\frac{\partial}{\partial\lambda})$.
\end{remar}
\begin{remar}
The function $\theta$ is the first Jacobi theta function.
\end{remar}

\begin{remar} We have that:
$$\rho'(\lambda ,\tau)=-\wp(\lambda ,\tau)+c(\tau)$$
where:
$$\wp(z,\tau)=\frac{1}{z^2}+\sum_{n,m\in {\ZZ}^2\backslash (0,0)}\Big(\frac{1}{(z-m-n\tau)^2}-\frac{1}{(m+n\tau)^2}\Big)$$
is the {\it Weierstrass} $\wp$-function.
\end{remar}

\begin{examp} For $g=\mathfrak{sl}_{N}$,
$n=1$, $V=S^{pN}{\CC}^{N}$ and $h={\CC}/{\CC}(1,...,1)$, the {\it KZB} equation becomes:
\begin{equation}
4\pi i\kappa\frac{\partial u}{\partial\tau}=\sum^{N}_{i=1}\frac{\partial^2 u}{\partial\lambda^{2}_{i}}+
2p(p+1)\sum_{i<j}\rho'(\lambda_{i}-\lambda_{j},\tau)u.
\end{equation}
The RHS of the previous equation is the hamiltonian operator of the {\rm quantum elliptic Calogero-Moser N-body 
system.}
\end{examp} 

\begin{remar}
A remarkable fact about all form of KZ equation is that them can be realized geometrically. Their solutions
can be expressed in terms of hypergeometric integrals depending on parameters. The KZ equation(s) can also 
be realized as Gauss-Manin connection.
\end{remar}

\begin{examp}
Let $n=1$, $p=1$ so $V=L_2$ $(\dim L_2=3)$. In this case the solution of the KZB equation can 
be written in the following form:
$$u(\lambda ,\tau)=u_{g}(\lambda ,\tau)=\int_{0}^{1}\Big(\frac{\theta (t,\tau)}{\theta'(0,\tau)}\Big)^{-\frac{2}{\kappa}}
\frac{\theta(\lambda -t,\tau)\theta'(0,\tau)}{\theta(\lambda ,\tau)\theta(t,\tau)}g(\lambda-\frac{2}{\kappa}t,\tau)dt,$$
where $g=g(\lambda,\tau)$ is any holomorphic solution of:
$$2\pi i\kappa\frac{\partial g}{\partial\tau}=\frac{\partial^2 u}{\partial\lambda^2}.$$
For instance we can take the function:
$$ g=e^{\lambda\mu+\frac{\mu^2}{2\pi i\kappa}\tau}\hspace{10pt}\mu\in{\CC}.$$
\end{examp}

\begin{remar} The simplest KZ equations have the classical Gauss hypergeometric function as their solution.
So it is natural to consider solutions of all KZ type equations as generalized hypergeometric functions.
The Gauss hypergeometric function is:
$$\frac{\Gamma(b)\Gamma(c-b)}{\Gamma(c)}F(a,b,c;z)=\int_1^{\infty}t^{a-c}(t-1)^{c-b-1}(t-z)^{-a}dt$$
We can see that the function:
$$\Big(\frac{\theta (t,\tau)}{\theta'(0,\tau)}\Big)^{\alpha}$$
is the analogue of the function $t^{\alpha}$, while:
$$\frac{\theta(\lambda -t,\tau)\theta'(0,\tau)}{\theta(\lambda ,\tau)\theta(t,\tau)}$$
is the analogue of the function $\frac{1}{t}$.
\end{remar}

\section{Asymptotic Solutions}
Let's consider the asymptotic behavior of the function $u=u(\lambda ,\tau)$, i.e let's study the behavior
of this function for $\kappa\longrightarrow 0$.

\begin{remar}
Under this hypothesis we are dealing with the following eigenvalue problem:
$$\frac{d^2}{d\lambda^2}v(\lambda)-2{\wp}(\lambda ,\tau)v(\lambda)=Ev(\lambda)$$
where $E$ is a number. The equation:
\begin{equation}
\frac{d^2}{d\lambda^2}v(\lambda)-p(p+1){\wp}(\lambda ,\tau)v(\lambda)=Ev(\lambda)
\end{equation}
is called {\rm Lame'} equation and its solutions are called Lame' functions.
\end{remar}

We are looking for solutions of the form:
$$u(\lambda ,\tau)=e^{\frac{S(\tau)}{\kappa}}(f_{0}(\lambda ,\tau)+\kappa f_{1}(\lambda ,\tau)+...).$$
Plugging this expression into the equation:

\begin{equation}
2\pi i\kappa\frac{\partial u}{\partial\tau}=\frac{\partial^2 u}{\partial\lambda^2}-p(p+1)\rho'(\lambda ,\tau)u
\end{equation}
we get recurrence relations for the functions $f_{i}(\lambda ,\tau)\hspace{10pt}n=0,1,....$.

\begin{defin} We say that the function $u=u(\lambda ,\tau)$ is a formal asymptotic solution if all the
recurrence relations are satisfied.
\end{defin}

\begin{corol} If $u=u(\lambda ,\tau)$ is an asymptotic solution then for every fixed $\tau$ the function
$f_{0}=f_{0}(\lambda ,\tau)$ is an eigenvector of the operator:
$$\frac{\partial^2}{\partial\lambda^2}+p(p+1)\rho'$$
with eigenvalue $\frac{\partial S}{\partial\tau}$.
\end{corol}
Applying the stationary phase method to the integral:

$$u(\lambda ,\tau)=\int_{0}^{1}\Big(\frac{\theta (t,\tau)}{\theta'(0,\tau)}\Big)^{-\frac{2}{\kappa}}
\frac{\theta(\lambda -t,\tau)\theta'(0,\tau)}{\theta(\lambda ,\tau)\theta(t,\tau)}
e^{\lambda\mu+\frac{\mu^2}{2\pi i\kappa}\tau}dt,$$
we get the following:

\begin{theor}
For any $\mu\in{\CC}$, the function:
$$v(\lambda)=e^{\lambda\mu}\frac{\theta(\lambda -t_{0}, \tau)}{\theta (\lambda ,\tau)}$$
is an eigenfunction of:
$$\frac{\partial^2}{\partial\lambda^2}-2\wp(\lambda ,\tau)$$
if $t_{0}$ is a critical point, w.r.t $t$, of the function:
$$\theta(t,\tau)e^{-\mu t}.$$
\end{theor}

\begin{remar} This is an example of Bethe ansatz. Hermite developed this method for the Lame' equation (1872).
\end{remar}
\begin{remar} Similar formulas are available for the eigenfunctions of  the Many-Body Elliptic Calogero-Moser 
hamiltonian.
\end{remar}

\section{Special Solutions and Conformal Blocks}
In what follows we will consider the differential equation:
\begin{equation}
2\pi i\kappa\frac{\partial u}{\partial\tau}=\frac{\partial^{2} u}{\partial\lambda^2}+p(p+1)\rho'(\lambda ,\tau)u.
\label{oleg}\end{equation}
Let's fix $\kappa\in{\ZZ}_{>0}$ and $p\in{\ZZ}_{\geq 0}$ and assume that $\kappa\geq 2p+2$.
We will consider holomorphic solutions $u=u(\lambda ,\tau)$ of the KZB equation with the following properties:
\begin{enumerate}
\item $u(\lambda +2,\tau)=u(\lambda ,\tau)$;
\item $u(\lambda +2\tau,\tau)=e^{-2\pi i\kappa(\lambda+\tau)}u(\lambda,\tau)$;
\item $u(-\lambda ,\tau)=(-1)^{p+1}u(\lambda ,\tau)$;
\item $u(\lambda ,\tau)={\mathcal O}((\lambda -m-n\tau)^{p+1})$, as
	$\lambda\longrightarrow m+n\tau, \hspace{5pt} m,n\in{\ZZ}$.
\end{enumerate}

\begin{defin}
Holomorphic solutions of \ref{oleg} satisfying the conditions 1., 2., 3., and 4. are called {\rm conformal blocks} associated
with the elliptic curve $T=\frac{\CC}{\ZZ +\tau{\ZZ}}$ and the point $z=0$.
\end{defin}

\begin{remar} It known that under this hypothesis the space of conformal blocks has dimension equal to 
$\kappa -2p-1$.
\end{remar}

\section[Theta functions of level $\kappa$.]{Theta functions of level $\kappa$.}

Before describing the space of conformal blocks it is useful discuss some elementary properties of theta 
functions.
\begin{defin}
Let $\kappa\in{\ZZ}_{>0}$. A holomorphic function $u=u(\lambda ,\tau)$ is called a theta function of level 
$\kappa$ w.r.t the lattice $2{\ZZ}+2\tau{\ZZ}$ if:
\begin{enumerate}
\item $u(\lambda +2,\tau)=u(\lambda ,\tau)$;
\item $u(\lambda +2\tau,\tau)=e^{-2\pi i\kappa(\lambda+\tau)}u(\lambda,\tau)$.
\end{enumerate}
\end{defin}

\begin{lemma}
The space of theta functions of level $\kappa$ is $2\kappa$ dimensional. A basis is given by:
$$\theta_{n,\kappa}(\lambda ,\tau)=\sum_{j\in {\ZZ}}e^{2\pi i\kappa(j+\frac{n}{2\kappa})^2+2\pi i\kappa(j+\frac{n}{2\kappa})\lambda}$$
$n=0,...,2\kappa -1$
\end{lemma}
\begin{proof}
Follows from the definition and from the Fourier expansion.
\end{proof}

These functions have the following remarkable:
\begin{prope} Let $q=e^{\frac{\pi i}{\kappa}}$:
\begin{enumerate}
\item $\theta_{n,\kappa}(\lambda ,\tau)$ is a solution of:
$$2\pi i\kappa\frac{\partial g}{\partial\tau}=\frac{\partial^2 g}{\partial\lambda^2};$$
\item $\theta_{n,\kappa}(\lambda ,\tau)=\theta_{n+2\kappa,\kappa}(\lambda ,\tau)$;
\item $\theta_{n,\kappa}(-\lambda ,\tau)=\theta_{-n,\kappa}(\lambda ,\tau)$;
\item $\theta_{n,\kappa}(\lambda+\frac{2}{\kappa},\tau)=q^{2n}\theta_{n,\kappa}$;
\item $\theta_{n,\kappa}(\lambda+\frac{2t}{\kappa},\tau)e^{-2\pi it-\frac{2\pi i\tau}{\kappa}}\theta_{n+2,\kappa}(\lambda ,\tau)$.
\end{enumerate}
\end{prope}

\subsection{Modular Transformations}

Let's start with:

\begin{defin}
The {\rm Modular Group} is $M=\SL(2,{\ZZ})/(\pm Id)$. The generators are:
$$T=\mdpd 1101\hspace{10pt} S=\mdpd 01{-1}0$$
and they satisfy the relations: 
$$S^2=1,\hspace{10pt}(ST)^3=1.$$
\end{defin}
$M$  acts on $\HH$ by $T:\tau\longrightarrow\tau +1\hspace{10pt}S:\tau\longrightarrow -\frac{1}{\tau}$.
\begin{remar}
Elliptic curves corresponding to parameters $\tau,\tau +1,-\frac{1}{\tau}$ are isomorphic so it natural to
expect relations among conformal blocks associated to these curves.
\end{remar}
Let's introduce the following maps:
\begin{defin}
\begin{equation}
(Tu)(\lambda ,\tau)=u(\lambda ,\tau+1)
\end{equation}
and
\begin{equation}
(Su)(\lambda ,\tau)=e^{-\pi i\kappa\frac{\lambda^2}{2\tau}}\tau^{-1/2}u({\lambda}/{\tau},-1/{\tau}).
\end{equation}
\end{defin}

\begin{defin} $\Theta_{\kappa}$ will denote the space of theta functions of level $\kappa$.
\end{defin}

\begin{lemma} The maps $T$ and $S$ preserve the space of $\Theta_{\kappa}$-functions:
\begin{equation}
1.\hspace{60pt}\theta_{n,\kappa}(\lambda ,\tau +1)=q^{\frac{n^2}{2}}\theta_{n,\kappa}(\lambda ,\tau);
\end{equation}
\begin{equation}
2.\hspace{5pt}\theta_{n,\kappa}\big(\frac{\lambda}{\tau}, -\frac{1}{\tau}\big)=\sqrt{-\frac{i\tau}{2\kappa}}
e^{\pi i\kappa\frac{\lambda^2}{2\tau}}\sum_{m=0}^{2\kappa -1}q^{-mn}\theta_{m, \kappa}(\lambda ,\tau).
\end{equation}
\end{lemma}

\begin{remar}
The maps $T$ and $S$ define a projective representation of the Modular group in the space of
$\Theta_{\kappa}$-functions.
\end{remar}

\subsection{Symmetric Theta Functions}

\begin{defin}
Let's define:
\begin{equation}
\theta_{n,\kappa}^{S}(\lambda ,\tau)=\theta_{n,\kappa}(\lambda ,\tau)+\theta_{n,\kappa}(-\lambda ,\tau)=
\theta_{n,\kappa}(\lambda ,\tau)+\theta_{-n,\kappa}(\lambda ,\tau).
\end{equation}
$\theta_{n,\kappa}^{S}$ is called symmetric theta function of level $\kappa$.
\end{defin}

It is possible to prove the following:
\begin{lemma}
The space of symmetric theta functions of level $\kappa$ has dimension equal to $\kappa +1$ and it is generated
by $\theta_{0,\kappa}^{S},....,\theta_{\kappa,\kappa}^{S}$.
\end{lemma}

Let's go back to the general theory of theta function and discuss the following:

\begin{examp}[First Jacobi theta function]
The first Jacobi theta  function is defined by the following formula:
\begin{equation}
\theta(\lambda ,\tau)=-\sum_{j\in{\ZZ}}e^{\pi i(j+\frac{1}{2})^{2}\tau+2\pi i(j+\frac{1}{2})(\lambda+\frac{1}{2})}.
\end{equation}
This function satisfies the following {\rm product formula}:
\begin{equation}
\theta(\lambda ,\tau)=ie^{\pi i(\frac{\tau}{4}-\lambda)}(\kappa ,\tilde{q})(\frac{\tilde{q}}{\kappa}, \tilde{q})(\tilde{q},\tilde{q})
\end{equation}
where $\tilde{q}=e^{2\pi i\tau}$, $\kappa =e^{2\pi i\lambda}$ and $(\kappa ,\tilde{q})=\prod_{n=0}^{\infty}(1-\kappa {\tilde q}^{n})$.
Moreover this function has the following:
\begin{prope}
\begin{enumerate}
\item $\theta(-\lambda ,\tau)=-\theta(\lambda ,\tau);$
\item $\theta(\lambda +1,\tau)=-\theta(\lambda ,\tau);$
\item $\theta(\lambda +\tau ,\tau)= e^{-\pi i\tau -2\pi i\lambda}\theta(\lambda ,\tau);$
\item $\theta(\lambda ,\tau)$ has simple zeros at ${\ZZ}+\tau{\ZZ}$;
\item $\theta(\lambda ,\tau)$ is a theta function of level 2;
\item $\theta(\lambda ,\tau)$ is a solution of the heat equation:
$$2\pi i\kappa\frac{\partial g}{\partial\tau}=\frac{\partial^2 g}{\partial\lambda^2};$$
\item{\rm (Modular property)} The function $\theta(\lambda ,\tau)$ satisfies the following two equations:
\end{enumerate}

\begin{equation}
\theta(\lambda ,\tau +1)=e^{\frac{\pi i}{4}}\theta(\lambda ,\tau);
\end{equation}

\begin{equation}
\theta\big(\frac{\lambda}{\tau},-\frac{1}{\tau}\big)\frac{\sqrt{-i\tau}}{\tau}e^{\pi i\lambda^2}{\tau}
\theta(\lambda ,\tau)
\end{equation}
From the property we see that $\theta$ generates a one dimensional subspace, invariant w.r.t the action 
of the modular group.
\end{prope}
\end{examp}

Let's go back to the general theory and let's start to investigate the relation between conformal blocks 
and theta functions. 
\begin{lemma}
For fixed $\tau$ the space of conformal blocks has dimension $\kappa -2p+1$ and consists of linear 
combinations of functions of the form:
$$\theta^{p+1}(\lambda ,\tau)\theta_{n,\kappa -2p}^{S}(\lambda ,\tau),\hspace{10pt}n=0,...,\kappa -2p.$$
\end{lemma}

\begin{examp} Let's suppose that $\kappa =2p+2$; then the space of conformal blocks is one dimensional. This 
implies that the holomorphic solutions of the KZB equation must be of the form:
\begin{equation}
u(\lambda ,\tau)=A(\tau)\theta^{p+1}(\lambda ,\tau).
\end{equation}
\end{examp}
\begin{lemma}
The function $u(\lambda ,\tau)=\theta^{p+1}(\lambda ,\tau)$ is a solution of the KZB heat equation for 
$\kappa =2p+2$.
\end{lemma}

\section{Integral Formulas}

In this section we will introduce the integral formulas for the conformal blocks. 
\begin{notat}
Let's define:
\begin{equation}
\sigma_{\lambda}(t,\tau)=\frac{\theta(\lambda -t,\tau)\theta'(0,\tau)}{\theta(\lambda ,\tau)\theta(t,\tau)};
\end{equation}
\begin{equation}
E(t,\tau)=\frac{\theta (t,\tau)}{\theta'(0,\tau)}.
\end{equation}
Let's also define the Master function:
\begin{equation}
\Phi(t_1,...,t_p)=\prod_{j=1}^{p}E(t_j ,\tau)^{-\frac{2p}{\kappa}}\prod_{1\leq i<j\leq p}
E(t_i-t_j,\tau)^{\frac{2}{\kappa}} \label{ciccio}.
\end{equation}
\end{notat}

\begin{remar}
\ref{ciccio} is the elliptic analogous of the Selberg integral.
\end{remar}

We also  introduce the integration cycles:

$${\Delta}_{k}\subset{\RR}^{k}\subset{\CC}^{k}$$
where
$${\Delta}_{k}=\{(t_{1},...,t_{k}\vert 0\leq t_{k}\leq ...\leq t_{1}\leq 1\}.$$
We define $\tilde{\Delta_{k}}$ the image of the standard real simplex under the map:
\begin{equation}
(t_{1},...,t_{k})\longrightarrow (\tau t_{1},...,\tau t_{k}). \label{ciccio2}
\end{equation}

\begin{remar} 
\ref{ciccio2} define a rotation of the standard simplex.
\end{remar}

For $0\leq k\leq p$ let's define:

$$
J^{[k]}_{n,\kappa}(\lambda ,\tau)=\int_{\Delta_{k}}\prod_{j=1}^{k}\sigma_{\lambda}(t_{j},\tau)dt_{j}
\int_{\tilde\Delta_{p-k}}\prod_{j=k+1}^{p}\sigma_{\lambda}(t_{j},\tau)dt_{j}\Phi(t_{1},...,t_{p},\tau)
\theta_{\kappa ,n}(\lambda+\frac{2}{\kappa}\sum_{j=1}^{p}t_{j},\tau).
$$

\begin{defin}
\begin{equation}
u^{[k]}_{n,\kappa}(\lambda ,\tau)=J^{[k]}_{n,\kappa}(\lambda ,\tau)+(-1)^{p+1}J^{[k]}_{n,\kappa}
(-\lambda ,\tau).
\end{equation}
\end{defin}
\begin{theor}
For any $0\leq k\leq p$ and any $n$
$$u_{n,\kappa}^{[k]}(\lambda ,\tau)$$
is a conformal block.
\end{theor}
We have the following:
\begin{defin}
The
$$u_{n,\kappa}^{[k]}(\lambda ,\tau)$$
are called {\rm elliptic hypergeometric functions.}
\end{defin}
Now we state some elementary properties of the elliptic hypergeometric functions:
\begin{prope}
\begin{enumerate}
\item $u_{n,\kappa}^{[k]}(\lambda ,\tau)=u_{n,2+\kappa}^{[k]}(\lambda ,\tau)$;
\item $u_{n,\kappa}^{[k]}(\lambda ,\tau)=-q^{2k(n+p-k)}u^{[k]}_{\kappa ,-n-2(p-k)}$.
\end{enumerate}
\end{prope}

\begin{theor}[on Basis]
\begin{enumerate}
\item The set $u_{n,\kappa}^{[p]}(\lambda ,\tau),\hspace{5pt}n=p+1,...,\kappa -p-1$, is a basis in the space
of conformal blocks.
\item $u_{n,\kappa}^{[p]}\equiv 0$ if $0\leq n\leq p$ or $\kappa -p\leq n\leq \kappa$.
\end{enumerate}
\end{theor}
\begin{examp}
Let $\kappa=2p+2$ so the dimension of the space of conformal blocks is equal to 1. In this case we have that:
$$u_{n,\kappa}^{[p]}(\lambda ,\tau)\equiv 0\hspace{5pt}if\hspace{5pt}n\neq p+1\hspace{5pt}mod\hspace{5pt}\kappa.$$
\end{examp}
In this case the $p$-dimensional integral can be calculate explicitly. This integral can be called the elliptic
Selberg integral.
Recall that:

$$B_{p}(\alpha ,\beta ,\gamma)=\int_{\Delta_{p}}\prod_{j=1}^{p}t_{j}^{\alpha -1}(1-t_{j})^{\beta -1}
\prod_{1\leq i<j \leq p}(t_{j}-t_{k})^{2\gamma}dt_{1}\wedge ...\wedge d_{k}=$$
$$\frac{1}{p!}\prod_{j=0}^{p-1}\frac{\Gamma(1+\gamma (1+j))}{\Gamma (1+\gamma)}\frac{\Gamma (\alpha +j\gamma)
\Gamma (\beta +j\gamma)}{\Gamma(\alpha +\beta +(p+j-1)\gamma)}.$$

Now we can state the following:
\begin{theor}

$$u_{2p+2}^{[p]}(\lambda ,\tau)=C_{p}B_{p}\big(\frac{1}{2}+\frac{1}{2(p+1)}, -1+\frac{1}{p+1}, 
\frac{1}{2(p+1)}\big)\theta(\lambda ,\tau)^{p+1},$$
where:
$$C_{p}=(2\pi )^{\frac{p}{2}}e^{-\pi i\big[\frac{p(3p-1)}{4(p+1)}+\frac{p+1}{2}\big]}\prod_{j=1}^{p}
\big(1-e^{-2\pi i\frac{j}{2p+2}}\big)$$
\end{theor}

\begin{reftc}
\cite{FV}, \cite{FSV1,FSV2}
\end{reftc}


\chapter{}
During the last lecture we introduced the elliptic analogous of the KZ differential 
equations. In particular, given $\tau\in {\HH}_{+}$ and $\lambda\in {\CC}$, we studied the 
following equation:
\begin{equation}
2\pi i\kappa\frac{\partial u}{\partial\tau} = \frac{\partial^{2}u}{\partial\lambda^{2}}
+p(p+1)\rho'(\lambda ,\tau)u, \label{ciccio104}
\end{equation}
where $\kappa ,p\in{\ZZ}_{+}$ and $\kappa\geq 2p+2$ and $\rho$ is the logarithmic derivative of the first 
Jacobi 
theta function:
$$\rho =\frac{\theta'(\lambda ,\tau)}{\theta'(\lambda ,\tau)}.$$
We also introduced the space of conformal blocks, i.e the space of holomorphic solution of \ref{ciccio104}
that satisfy the following properties:
\begin{enumerate}
\item $u(\lambda +2,\tau)=u(\lambda ,\tau)$;
\item $u(\lambda +2,\tau)=e^{-2\pi i\kappa (\lambda +\tau)}u(\lambda ,\tau)$;
\item $u(-\lambda,\tau)=(-1)^{p+1}u(\lambda ,\tau)$;
\item $u(\lambda ,\tau)={\mathcal O}((\lambda -m-n\tau)^{p+1}))$ as
	$\lambda\longrightarrow n+m\tau\in{\ZZ}+\tau {\ZZ}$.
\end{enumerate}
During this lecture we will study the action of the modular group on the space of conformal blocks. To this
end we remind the integral formulas introduced at the end of the last lecture.
We defined:
$$E(t,\tau)=\frac{\theta(t,\tau)}{\theta'(0,\tau)},$$
$$\sigma_{\lambda}(t,\tau)=\frac{\theta(\lambda -t,\tau)\theta'(0,\tau)}{\theta(\lambda ,\tau)\theta 
(t,\tau)},$$
and:
$$\Phi(t_1,...,t_p)=\prod_{j=1}^{p}E(t_j ,\tau)^{-\frac{2p}{\kappa}}\prod_{1\leq i<j\leq p}
E(t_i-t_j,\tau)^{\frac{2}{\kappa}} \label{ciccio10}.$$
If ${\Delta}_{k}=\{(t_{1},...,t_{k}\}\vert 0\leq t_{k}\leq ...\leq t_{1}\leq 1\}$ is the standard k-simplex 
and  $\tilde{\Delta}_{k}$ is the image of $\Delta_{k}$ under the map:
$$(t_{1},...,t_{k})\longrightarrow (\tau t_{1},...,\tau t_{k}),$$
we have that:

\begin{equation}
J^{[k]}_{n,\kappa}(\lambda ,\tau)=\int_{\Delta_{k}}\prod_{j=1}^{k}\sigma_{\lambda}(t_{j},\tau)dt_{j}
\int_{\tilde\Delta_{p-k}}\prod_{j=k+1}^{p}\sigma_{\lambda}(t_{j},\tau)dt_{j}\Phi(t_{1},...,t_{p},\tau)  
\theta_{\kappa ,n}(\lambda+\frac{2}{\kappa}\sum_{j=1}^{p}t_{j},\tau), \label{ciccio105}
\end{equation}

where $\theta_{\kappa ,n},\hspace{5pt}$ is a theta function of level $\kappa$.\\
We also defined the functions:
$$u^{[k]}_{n,\kappa}(\lambda ,\tau)=J^{[k]}_{n,\kappa}(\lambda ,\tau)+(-1)^{p+1}J^{[k]}_{n,\kappa}
(-\lambda ,\tau)$$
we observed that:
\begin{prope}
\begin{enumerate}
\item $u_{n,\kappa}^{[k]}(\lambda ,\tau)=u_{n,2+\kappa}^{[k]}(\lambda ,\tau)$;
\item $u_{n,\kappa}^{[k]}(\lambda ,\tau)=-q^{2k(n+p-k)}u^{[k]}_{\kappa ,-n-2(p-k)}$
\end{enumerate}
\end{prope}
and we stated the following important:
\begin{theor}[on Basis]
\begin{enumerate}
\item The set $u_{n,\kappa}^{[p]}(\lambda ,\tau),\hspace{5pt}n=p+1,...,\kappa -p-1$, is a basis in the space
of conformal blocks.
\item $u_{n,\kappa}^{[p]}\equiv 0$ if $0\leq n\leq p$ or $\kappa -p\leq n\leq \kappa$.
\end{enumerate}
\end{theor}

In what follows it will be convenient use the following:

\begin{notat}
\begin{align}
q   & = e^{\frac{\pi i}{\kappa}},\\
[n] & =\frac{q^{n}-q^{-n}}{q-q^{-1}}.
\end{align}
\end{notat}
The following theorem describes some important relations satisfied by the functions $u_{\kappa ,n}^{[k]}$ 
with different $k$.
\begin{remar}
We recall that $k$ represents the number of horizontal integrations performed in \ref{ciccio105}.
\end{remar}
\begin{theor}[Stokes Theorem] For $0\leq k\leq p-1$ and any $n$, we have:
\begin{equation}
[p-k]\big(q^{m+p-k}-q^{-m-p+k} \big)u_{\kappa ,m}^{[k]}=q^{-m-k-1}[k+1]\big(q^{-2(k+1)}
u_{\kappa ,m+2}^{[k+1]}-u_{\kappa ,m}^{[k+1]}\big).
\end{equation}
\end{theor}
\begin{proof}
We will give  the idea of the proof of the theorem. Let $p=1$ and $k=0$. In these hypothesis the
integrals are 1 dimensional:
\begin{equation}
\int_{interval}\Phi (t,\tau) \sigma_{\lambda}(t,\tau)\theta_{\kappa ,m}(\lambda +\frac{2}{\kappa}t,\tau)dt.
\end{equation}
Now will perform the integration on the fundamental domain $(1,\tau)$. Since the integrand is holomorphic
we get that:
$$I_{1}+I_{2}+I_{3}+I_{4}=0$$ 
and from the periodicity properties of the $\theta$-functions, we get $I_{1}=c_{1}u_{m}^{[1]}$, 
$I_{4}=c_{4}u_{m}^{[0]}$, $I_{2}=c_{2}u_{m}^{[0]}$, $I_{3}=c_{1}u_{m+2}^{[1]}$.
\end{proof}
\begin{remar} The Stokes theorem allows to express an arbitrary function $u_{\kappa ,m}^{[k]}$ in terms of 
basic ones.
\end{remar}

\section{Modular Transformations of Conformal Blocks}

Let's introduce the following two transformations:

\begin{align}
(Tu)(\lambda ,\tau) &= u(\lambda ,\tau +1),\\
(Su)(\lambda ,\tau) &= e^{-\pi i\kappa\frac{\lambda^2}{2\tau}}\tau^{-\frac{1}{2}}u(\frac{\lambda}{\tau},-\frac{1}{\tau}).
\end{align}
We have the following:
\begin{theor} If $u=u(\lambda ,\tau)$ is a solution of the KZB heat equation the $Tu$ and $Su$ a also solutions
of the same equation. Moreover the maps $T$ and $S$ preserve the space of conformal blocks.
\end{theor}
We also have the following:
\begin{propo}
The restrictions of the maps $T$ and $S$ to space of conformal blocks satisfy the following identities:
\begin{align}
S^2      &= (-1)^{p}iq^{-p(p+1)}Id, \\
(ST)^{3} &=(-1)^{p}iq^{-p(p+1)}Id,
\end{align}
i.e they define a projective representation of the modular group on the space of conformal blocks.
\end{propo}
Now we want to discuss the following problem: find an explicit description of the operators $S$ and $T$.

The solution of this problem is contained in the following two lemmas:
\begin{lemma}
$$Tu_{\kappa ,n}^{[p]}(\lambda ,\tau)=q^{\frac{n^2}{2}}u_{\kappa ,n}^{[p]},$$
\end{lemma}
and

\begin{lemma}
\begin{equation}
Su_{\kappa ,n}^{[p]}(\lambda ,\tau)=\frac{e^{-\frac{\pi i}{4}}}{\sqrt {2\kappa}}\sum_{m=0}^{2\kappa -1}
q^{-mn}u_{\kappa ,m}^{[0]}(\lambda ,\tau). \label{ciccio106}
\end{equation}
\end{lemma}

\begin{remar}
The proofs of the previous lemmas are straightforward. We observe that in \ref{ciccio106} the RHS is given
in terms of functions obtained taking p-vertical integrations while the LHS is a function obtained taking
p-horizontal integrations.
To get the formula \ref{ciccio106} we need to apply p-times the Stokes theorem to get the form:
\begin{equation}
Su_{\kappa ,n}^{[p]}=\sum_{m}S_{m,n}u_{\kappa ,m}^{[p]} \label{ciccio107}.
\end{equation}
\end{remar}
To describe the coefficients $S_{n,m}$ in \ref{ciccio107} of the previous remark we need the following:
\begin{theor}
\begin{equation}
S_{n,m}=\sqrt{\frac{-i}{2\kappa}}q^{p(n-p)-\frac{p(p+1)}{2}}(q^{-m}-q^{m})\prod_{j=1}^{p}(q^{-n+j}-q^{n-j})
P_{n-p-1}^{(p+1)}(m).
\end{equation}
\end{theor}

\begin{defin} The function $P_{n-p-1}^{(p+1)}(x)$ is the {\it Macdonald polynomial} of level $p+1$ and degree
$n-p-1$ associated with the Lie algebra $\slt$.
\end{defin}

\section{Macdonald Polynomials}
The Macdonald polynomials of level $k$ associated with $\slt$ are $x$-even polynomials of $q^{mx}$, $m\in{\ZZ}$
defined by the following rules:
\begin{enumerate}
\item $P_{0}^{(k)}(x)\equiv 1$ and $P_{n}^{(k)}(x)=q^{nx}+q^{-nx}$+lower order terms;
\item $\langle P_{m}^{(k)}(x), P_{n}^{(k)}(x) \rangle=0$ if $m\neq n$, where:
$$ \langle f(x),g(x) \rangle
	=constant\hspace{5pt}term\hspace{5pt}\big(f(x)g(x)\prod_{j=0}^{k-1}(1-q^{2(j+x)})(1-q^{2(j-x)}).$$
\end{enumerate}

\begin{examp} For $n>0$ we have that:
$$P_{n}^{(0)}(x)=q^{nx}+q^{-nx}.$$
\end{examp}

In the following we will define an operator that play an important role in the theory of Macdonald 
polynomials:

\begin{defin}[Shift Operator] The {\it shift operator} $D$, in the $\slt$-case, is defined by:
\begin{equation}
Df(x)=\frac{f(x-1)+f(x+1)}{q^{x}+q^{-x}}.
\end{equation}
where $f=f(x)$ is a given function.

\end{defin}

Given the shift operator it is possible to obtain Macdonald polynomials of higher level knowing polynomials
of lower level. In fact we have the following:
\begin{theor}[Askey-Ismail]
\begin{equation}
(DP_{n}^{(k)})(x)=(q^{-n}-q^{n})P_{n-1}^{(k+1)}(x).\label{ciccio1010}
\end{equation}
\end{theor}

\begin{examp}
For $n>0$ we have that:
\begin{equation}
P_{n-1}^{1}(x)=\frac{q^{nx}+q^{-nx}}{q^{x}+q^{-x}}.\label{ciccio108}
\end{equation}
This follows applying the definition of the shift operator and previous theorem to
$P_{n}^{(0)}(x)=q^{nx}+q^{-nx}$. 
From \ref{ciccio108} we get that $P_{0}^{(2)}(x)=1$, $P_{1}^{(2)}(x)=q^{x}+q^{-x}$ and that
$$P_{2}^{(2)}(x)=q^{2x}+\frac{1+2q^2+q^4}{1+q^2+q^4}+q^{-2x}.$$
\end{examp}
The following theorem describes the relation between the shift operator and the modular transformations of
the basic elliptic hypergeometric functions:
\begin{theor}
For $0\leq k\leq p$ we have:
$$
Su_{\kappa ,n}^{[p]}(\lambda ,\tau)=\sqrt{\frac{-i}{2\kappa}}\sum_{m=-p+2k+1}^{\kappa -p-1}
q^{pn-km-\frac{k(k+1)}{2}}{p\brack k}^{-1}_q(q^{m-p+k}-q^{m+p-k})$$
$$\prod_{j=1}^{k}(q^{-n+j}-q^{n-j})P_{n-k+1}^{(k+1)}(m+p-k)u_{m}^{[k]}(\lambda ,\tau).
$$
\end{theor}

\begin{remar} Comparing this formula and the Stokes theorem we see that the Stokes formula gives the shift
operator relation \ref{ciccio1010}.
\end{remar}

\section{Trace of Intertwining Operators}

In what follows $q=e^{\frac{\pi i}{\kappa}}$. The quantum group $U_{q}(\slt)$ has generators $E$, $F$ and 
$q^{ch}$, where $c\in {\CC}$, with relations:
\begin{enumerate}
\item $q^{ch}q^{c'h}=q^{(c+c')h}$;
\item $q^{ch}Eq^{-ch}=q^{2ch}E$;
\item $q^{ch}Fq^{-ch}=q^{-2ch}F$;
\item $EF-FE=\displaystyle\frac{q^h-q^{-h}}{q-q^{-1}}$. The comultiplication is defined by the 
following:
\end{enumerate}
$\Delta (E)=E\otimes q^h+1\otimes E$,\hspace{5pt} $\Delta (F)=F\otimes 1+q^{-h}\otimes F$,\hspace{5pt} 
$\Delta (q^{ch})=q^{ch}\otimes q^{ch}$.
Let's identify weights for $U_{q}(\slt)$ with complex number as follows and let's say that a vector $v$ in
a $U_{q}(\slt)$-module has weight $\nu\in{\CC}$ if $q^{h}v=q^{\nu}v$. Let $M_{\mu}$ the Verma module over
$U_{q}(\slt)$ with highest weight $\mu$ and let $v_\mu$ be its highest vector. Let $k$ be a non negative
integer such that $\kappa\geq 2k+2$, $U$ an irreducible finite dimensional representation of $U_{q}(\slt)$
of weight $2k$ and let's denote with $U[0]$  the zero weight subspace of $U$. Let $u\in U[0]$.
For generic $\mu$ let:
\begin{equation}
\Phi_{\mu}^{\nu}:M_{\mu}\longrightarrow M_{\mu}\otimes U
\end{equation}
be the intertwining operator defined by:
\begin{equation}
\Phi_{\mu}^{\nu}v_{\mu}=v_{\mu}\otimes u+\frac{Fv_{\mu}}{[-\mu]}\otimes Eu+...+
\frac{F^{j}v_{\mu}}{[j]!(-\mu ,q)_{j}}\otimes E^{j}u+...
\end{equation}
where $(-\mu ,q)_{j}=[n][n+1]...[n+j-1]$.\\
Introduce now and End(U[0])-valued function $\psi^{(k)}(q,\nu ,\mu)$ defined by:
\begin{equation}
\psi^{(k)}(q,\nu ,\mu)u=Tr\vert_{M_{\mu}}(\Phi_{\mu}^{\nu}q^{\nu h}).
\end{equation}
Since $U[0]$ is 1 dimensional, this function is scalar. It is is possible give an explicit description of this 
function:
\begin{theor}[Etingof-Varchenko] The function $\psi^{k}(q,\nu ,\mu)$ is given by the formula:
$$
\psi^{(k)}(q,\nu ,\mu)=q^{\nu\mu}\sum_{j=0}^{k}(-1)^jq^{\frac{j(j-3)}{2}}(q-q^{-1})^{-j-1}
\frac{[k+j]!}{[j]![k-j]!}\frac{q^{-j\mu-(j-1)\nu}}{\prod_{l=0}^{j-1}[\mu -l]\prod_{l=0}^{j}[\nu -l]}.$$
\end{theor}
Let's introduce now the renormalized trace functions:
\begin{defin}
$$\Psi^{(k)}(q,\nu ,\mu)=\prod_{j=1}^{k}\Big(\frac{q^{\mu +1-j}-q^{-\mu -1+j}}{q^{\nu +j}-q^{-\nu -j}}\Big).$$
\end{defin}
\begin{remar}
The functions $\Psi^{(k)}$ are holomorphic functions of $\mu$.
\end{remar}
The relation between Macdonald polynomials and the trace functions is given by the following: 
\begin{theor}[Etingof-Styrkas]
$$
\Psi^{(k)}(q^{-1},-m-p+k,n-1)-\Psi^{(k)}(q^{-1},-m-p+k,-n-1)=P_{n-k-1}^{(k+1)}(x)\prod_{j=1}^{k}(q^{-n+2j}-q^n).
$$
\end{theor}
Now we will explain how the trace functions enter in the picture.
Our starting point was the formula \ref{ciccio106}:
$$Su_{\kappa ,n}^{[p]}(\lambda ,\tau)=\frac{e^{-\frac{\pi i}{4}}}{\sqrt {2\kappa}}\sum_{m=0}^{2\kappa -1}
q^{-mn}u_{\kappa ,m}^{[0]}(\lambda ,\tau).$$
Then using the Stokes theorem:
$$[p-k]\big(q^{m+p-k}-q^{-m-p+k} \big)u_{\kappa ,m}^{[k]}=q^{-m-k-1}[k+1]\big(q^{-2(k+1)}
u_{\kappa ,m+2}^{[k+1]}-u_{\kappa ,m}^{[k+1]}\big),$$
we obtained for any $0\leq k\leq p$ an expression of the form:
$$
Su_{\kappa ,n}^{[p]}(\lambda ,\tau)=\sqrt{\frac{-i}{2\kappa}}
\Big[\sum_{m=-p+2k+1}^{k-p-1}f^{(k)}_{m,n}u_{\kappa ,m}^{[k]}(\lambda ,\tau)+
\sum_{m=\kappa -p+2k+1}^{2\kappa -p-1}f^{(k)}_{m,n}u_{\kappa ,m}^{[k]}(\lambda ,\tau)\Big]$$
for suitable numbers $f^{(k)}_{m,n}$ given by the following:
\begin{theor}
\begin{equation}
f^{(k)}_{m,n}=q^{pn-km-k(k+1)}(q^{m+p-k}-q^{-m-p+k}){p\brack k}_{q}^{-1}\Psi^{(k)}(q^{-1},-m-p+k,n-1).
\end{equation}
\end{theor}

\section{$q$-KZ Equations}

In what follows we will discuss the theory of the {\it q-hypergeometric} functions and the related KZ-equations.
These are {\it difference equations} that can be thought as a quantized version of the previous ones.
In the ``classical" case we had that given $V=V_{1}\otimes...\otimes V_{n}$, tensor product of $\slt$ modules, we could
define the following system of differential KZ-equations:

$$\kappa\frac{\partial u}{\partial z_i}=\sum_{j\neq i}\frac{{\Omega}^{(i,j)}}{z_{i}-z_{j}}u$$

where $u=u(z_1,...,z_n)$ is a $V$-valued function and $\Omega=e\otimes f+f\otimes e+\frac{1}{2}h\otimes h$ 
is the Casimir operator.
These equations define a flat connection on the flat vector bundle $V\times {\CC}^n\longrightarrow {\CC}^n$ , that 
has regular singularities along the diagonals.
The monodromy of this connection can be described in terms of the universal $R$-matrix associated
to the quantum group $U_{q}(\slt)$, where $q=e^{\frac {2\pi i}{\kappa}}$, acting on the tensor product
$V^{q}=V_{1}^{q}\otimes...\otimes V_{n}^{q}$, where the $V_{i}^{q}$ are $U_{q}(\slt)$-modules associated to the $\slt$
modules $V_{i}$ for $i=1,...,n$. Moreover the KZ equations can be also realized geometrically:\\
if we consider the map $p:{\CC}^{k+n}\longrightarrow {\CC}^{n}$ and the Master function 
$\Phi_{k,n}:{\CC}^{k+n}\longrightarrow {\CC}$, then the KZ-equations can be identified with the Gauss-Manin 
connection of the vector bundle $\mathcal H$ whose whose fibers are the homology groups $H_{k}(p^{-1}(z), \Phi)$. 
Now we want to study the discrete version of all these objects. There are three versions of the classical
KZ-equation:\\
{\it 1. rational;}\\
{\it 2. trigonometric} and;\\
{\it 3. elliptic}.\\
In what follows we will describe the discretization of the rational KZ-system. Let's consider the tensor
product of $\slt$ modules $V=V_{1}\otimes...\otimes V_{n}$.

\begin{defin}
The rational q-KZ equation on a $V$-valued function $\omega =\omega(z_1 ,...,z_n)$ has the form:
\begin{equation}
\omega(z_1 ,..,z_{m}+p,..,z_n)=K_{m}(z_1 ,...,z_n,p)\omega(z_1 ,...,z_n),\hspace{5pt} m=1,...,n.
\end{equation}
p is called the {\it step} of the equations and $K_{m}: V\longrightarrow V$ are linear operators.
\end{defin}

Now we have the following data: the step map, a trivial fibration 
$\pi :V\times {\CC}^{n}\longrightarrow {\CC}^{n}$ and the set of linear maps
$\big\{K_m(z,p)\big\}_{m=1,...,n}$. These operators define and identification of the neighboring fibers w.r.t
the lattice generated by $p{\ZZ}^{n}\subset {\CC}^{n}$:

\begin{defin} We'll call this structure a {\rm discrete connection} and we will say that such connection
is flat if the translation along the sides of any elementary square gives the identity operator, i.e:
\begin{equation}
K_{i}(z_1 ,..,z_j+p,..,z_n)K_{j}(z_1 ,...,z_n)K_{j}(z_1 ,..,z_i+p,..,z_n)K_{i}(z_1 ,...,z_n)
\end{equation}
\end{defin}

\begin{remar} The q-KZ discrete connection is flat and the q-KZ equations are the equations for flat 
sections of the discrete connection.
\end{remar}

In the next section we will discuss the explicit formulas for the q-KZ operators $K_m$.

\section{Yangian $Y(gl_2)$ and Rational R-matrices}

\begin{defin}
The {\it Yangian} $Y(gl_2)$ is an associative algebra with unit and generators $T_{i,j}^{(k)}$, $i,j=1,2$
$k=1,2,3,....$.
\end{defin}
To write the relations let's introduce the generating series:
\begin{equation}
T_{i,j}(u)=\delta_{i,j}+\sum_{k=1}^{\infty} T_{i,j}^{(k)}u^{-k}.
\end{equation}
The relations are given by the following:
\begin{equation}
(u-v)\big[T_{i,j}(u),T_{k,l}(v)\big]=T_{j,k}(v)T_{i,l}(u)-T_{k,j}(u)T_{i,l}(v).
\end{equation}

The Yangian is a Hopf algebra whose coproduct $\Delta :Y(gl_2)\longrightarrow Y(gl_2)\otimes Y(gl_2)$
is given by:
\begin{equation}
(\Delta T_{i,j})(u)=\sum_{k=1}^{2}T_{i,k}(u)\otimes T_{k,j}(u).
\end{equation}
On $Y(gl_2)$ are also defined:
\begin{enumerate}
\item a family of homomorphisms $\rho_{z}:Y(gl_2)\longrightarrow Y(gl_2)$, where:
\begin{equation}
(\rho_{z}T_{i,j})(u) = T_{i,j}(u-z),\qquad
\rho_{z}\big(\frac{1}{u^k}\big) = \frac{1}{u^k}\frac{1}{\big(1-\frac{z}{u}\big)^{k}},
\end{equation}
and:
\item the evaluation homomorphism $\epsilon :Y(gl_2)\longrightarrow U(\slt)$, such that:
\begin{equation}
T_{1,1}^{k}=\delta_{1,k}h \hspace{20pt} T_{1,2}^{k}=\delta_{1,k} f,
\end{equation}

\begin{equation}
T_{2,1}^{k}=\delta_{1,k}e \hspace{20pt} T_{2,2}^{k}=-\delta_{1,k} h.
\end{equation}
\end{enumerate}

If $V$ is any $\slt$-module:
\begin{defin}
let's denote $V(z)$ the Yangian module defined via the homomorphism:
\begin{equation}
\epsilon\circ\rho_{z}:Y(gl_2)\longrightarrow U(\slt).
\end{equation}
\end{defin}
Now let $V_1$ and $V_2$ two highest weight $\slt$ modules with highest weight vectors $v_1\in V_1$ and
$v_2\in V_2$. For generic complex numbers $x, y\in{\CC}$ the Yangian modules $V_1(x)\otimes V_2(y)$ and
$V_2(y)\otimes V_1(x)$ are isomorphic.

\begin{defin}
The rational R-matrix $R_{V_1, V_2}(x-y)$ can be defined as the unique element of $End(V_1\otimes V_2)$, such
that:
\begin{enumerate}
\item $R_{V_1, V_2}:v_1\otimes v_2\longrightarrow v_1\otimes v_2$
(identity operator on the tensor product of the highest weight vectors),
\item $\begin{CD}V_1(x)\otimes V_2(y)@>{R_{V_1, V_2}}>> V_1(x)\otimes V_2(y)
		@>{P}>>V_2(y)\otimes V_1(x)\end{CD}$
(the composition $P\circ R_{V_1, V_2}$ is an isomorphism of Yangian modules).
\end{enumerate}

\end{defin}

\begin{remar}
$P$ is the permutation operator.
\end{remar}

\begin{examp} Let $V_{1}=V_{2}={\CC}^2$, $v_{+}, v_{-}$ a basis for ${\CC}^2$ where $v_{+}$ is highest
weight vector. In  $V_{1}\otimes V_{2}$ we have a basis given by $v_{+}\otimes v_{+}$, $v_{+}\otimes v_{-}$,
$v_{-}\otimes v_{+}$ and $v_{-}\otimes v_{-}$.
Under these hypotheses we have that the R-matrix with respect to the ordered basis written above is:
\begin{equation}
R(x)=
\begin{pmatrix} 1 & 0 & 0 & 0 \\ 0 & \frac{x}{x+1} & \frac{1}{x+1} & 0 \\
 0 & \frac{1}{x+1} & \frac{x}{x+1} & 0 \\ 0 & 0 & 0 & 1 \end{pmatrix}\quad
\end{equation}
\end{examp}

\begin{remar} The R-matrices satisfy the QYBE:
$$
R_{V_1, V_2}^{(12)}(x_1-x_2)R_{V_1, V_3}^{(13)}(x_1-x_3)R_{V_2, V_3}^{(23)}(x_2-x_3)=
R_{V_2, V_3}^{(23)}(x_2-x_3)R_{V_1, V_3}^{(13)}(x_1-x_3)R_{V_1, V_2}^{(12)}(x_1-x_2)
$$
\end{remar}
Now let $V=V_{1}\otimes ...\otimes V_n$ the tensor product of highest weight $\slt$-modules. Then:
\begin{defin}
The q-KZ operators $K_m$ are defined by the following:
$$
K_m=R_{V_m, V_{m- 1}}(z_m-z_{m-1}+p)...R_{V_m, V_{1}}(z_m-z_{1}+p)e^{\mu h^{(m)}}$$
$$R_{V_m, V_{n}}(z_m-z_n)...R_{V_m, V_{m+1}}(z_m-z_{m+1})$$
\end{defin}
\begin{exerc}
Check that the q-KZ discrete connection defined by the operators $K_m$ given in the previous definition
is flat.\\
{\rm Hint.} Use the fact that the R-matrices satisfy the WYE.
\end{exerc}
\begin{remar}
The rational q-KZ difference equations turn into the KZ differential equations under the following limiting
procedure.\\
Let $z_m=SZ_m$, $m=1,...,n$, with $S>>1$ a number and $Z_m$ a new variable and $\mu=\frac{\nu}{S}$, where
$\nu$ is a new parameter. Introduce $u(Z)=\omega (Sz)$. Then the q-KZ equation is:
\begin{equation}
u(Z_1,...,Z_m+\frac{p}{S},...Z_n)=K_m(SZ,p)u(Z), \hspace{5pt}m=1,...,n\label{ciccio109}
\end{equation}
where:
\begin{equation}
K_m(SZ,p)=1+\frac{1}{S}\sum_{j\neq m}\frac{\Omega^{(m,j)}+c_{m,j}}{Z_m-Z_j}+\frac{\nu h^{(m)}}{S},
\end{equation}
for suitable numbers $c_{m,j}$. As $S\rightarrow\infty$ the difference equation \ref{ciccio109} turns into a 
differential equation:
\begin{equation}
p\frac{\partial u}{\partial z_m}=\sum_{j\neq m}\frac{\Omega^{(m,j)}+c_{m,j}}{z_m-z_j}+\nu h^{(m)},\hspace{5pt}m=1,...,n.
\end{equation}
\end{remar}

\begin{reftc}
\cite{FV}, \cite{FSV1,FSV2}, \cite{EV2}, \cite{EFK}
\end{reftc}


\chapter{}
During the last lecture we introduced the quantized version of the KZ-equations and we studied in detail the 
rational case. Let's recall the main features of this construction.
Let $V=V_{1}\otimes...\otimes V_{n}$ a tensor product of $\slt$-modules, then the q-KZ equation is given by:
\begin{equation}
\omega(z_1 ,..,z_{m}+p,..,z_n)=K_{m}(z_1 ,...,z_n,p)\omega(z_1 ,...,z_n),\hspace{5pt} m=1,...,n, \label{pippo1}
\end{equation}
where $u=u(z_1,...,z_n)$ is a $V$-valued function, $p$ is called the step of the equation,
$K_m:V\longrightarrow V,\hspace{5pt}m=1,...,n$ are 
linear operators that are given in terms of R-matrices, i.e:
$$
K_m=R_{V_m, V_{m- 1}}(z_m-z_{m-1}+p)...R_{V_m, V_{1}}(z_m-z_{1}+p)e^{\mu h^{(m)}}$$
$$R_{V_m, V_{n}}(z_m-z_n)...R_{V_m, V_{m+1}}(z_m-z_{m+1}).$$ 
We also described the geometrical interpretation of this discretized version of the KZ-equation: the q-KZ
equation can be interpreted as discrete connection on the trivial vector bundle 
$V\times {\CC}^n\longrightarrow {\CC}^n$. Moreover this connection is flat, i.e: 
$$K_{i}(z_1 ,..,z_j+p,..,z_n)K_{j}(z_1 ,...,z_n)=K_{j}(z_1 ,..,z_i+p,..,z_n)K_{i}(z_1 ,...,z_n).$$

Finally we proved that under a suitable limiting procedure, the q-KZ equation turns into the KZ equation.

\begin{lemma} Let $\vec{z}=(z_1,...,z_n)$ and let: 
\begin{equation}
\omega(\vec{z},p)=e^{\frac{S(\vec{z})}{p}}(f_{0}(\vec{z})+pf_1(\vec{z}).....), \label{pippo2}
\end{equation}
be an asymptotic solution of \ref{pippo1} , as $p\longrightarrow 0$. Then:
$$K_i(\vec{z},p=0)f_{0}(\vec{z})=\lambda_i(\vec{z})f_{0},\hspace{5pt}i=1,...,n,$$
where $\lambda_i(\vec{z})=e^{\frac{\partial S(\vec{z})}{\partial z_i}}.$
\end{lemma}
\begin{proof}
Since $\omega=\omega(\vec{z},p)$ is a solution of \ref{pippo1} we have that:
$$\omega(z_1,..,z_i+p,..,z_n,p)=K_{i}(\vec{z},p)\omega(\vec{z},p).$$
Plugging \ref{pippo2} in the previous equation we get:
\begin{equation}
\begin{split}
e^{\frac{S(z_1,..,z_i+p,..,z_n)}{p}}(f_0(z_1,..,z_i+p,..,z_n) &+ pf_1(z_1,..,z_i+p,..,z_n)+...)         \\
                                                              &= K_i(\vec{z},p)e^{\frac{S(\vec{z})}{p}}
                                                              (f_0(\vec{z})+....).\label{pippo3}        \\
\end{split}
\end{equation}
expanding with respect to the $z_i$ coordinate:
$$e^{\frac{S(z_1,..,z_i+p,..,z_n)}{p}}=e^{\frac{S(\vec{z})}{p}+\frac{\partial S}{\partial z_i}+....},$$
Plugging this result in the LHS of  \ref{pippo3} and taking the limit for $p\longrightarrow 0$, we get:
$$K_i(z,p=0)f_{0}(z)=\lambda_i(z)f_{0}(\vec{z}).$$
From the previous lemma we deduce that the leading term $f_0(z)$ is an eigenvector of commuting  operators
$\big\{K_m(z,p=0)\big\}_{m=1,...,n}$. These operators are the Hamiltonians of the $XXX$-quantum spin chain 
model. So having solutions of the q-KZ equation and computing its quasi-classical asymptotics it is possible
to construct eigenvectors for the operators $\big\{K_m(z,p=0)\big\}_{m=1,...,n}$. These eigenvectors are the
so called {\rm Bethe} vectors and they are constructed by the Bethe ansatz
of the $XXX$-model.
\end{proof}

The main topic of this lecture will be the quantization of the KZ system of differential equation and its 
geometric interpretations.
Let's start recalling the main features of the classical KZ differential equation.

\section{Classical KZ-system}
Given $m_1,...m_n,\kappa\in {\CC}$ consider the following function:
$$\Phi=\prod_{l=1}^{n}(t-z_l)^{\frac{m_l}{\kappa}}.$$
and let $\gamma_m=[z_m,z_{m+1}]$ be the oriented interval.
Let's define the following function:
\begin{equation}
I^{\gamma}(z_1,...,z_n)=\big(\int_{\gamma}\Phi\frac{dt}{t-z_1},...,\int_{\gamma}\Phi\frac{dt}{t-z_n}\big),
\end{equation}
where $\gamma=\gamma_1,...,\gamma_{n-1}$, then:
\begin{equation}
\frac{\partial I^{\gamma}}{\partial z_i}=\frac{1}{\kappa}\sum_{i\neq j}\frac{\Omega_{i,j}}{z_i-z_j}
I^{\gamma}.\label{pippo4}
\end{equation}
\begin{remar} In \ref{pippo4} $\Omega$ depends only on $m_1,...,m_n$.
\end{remar}
\begin{remar} The system \ref{pippo4} is an example of KZ equation.
\end{remar}
\subsection{Differential Forms}
For $z_1,...,z_n$ fixed consider the complex:
\begin{equation}
\begin{CD}
0\longrightarrow\Omega^{0} @>d>>\Omega^{1}\longrightarrow 0,\label{pippo20}
\end{CD}
\end{equation}
where by definition:
\begin{equation}
\Omega^{0}=\{f(t)\Phi(t,z)\vert f\hspace{5pt}\text{is a rational function on}\hspace{5pt}{\CC}\backslash \{{z_1,...,z_n}\}\},
\end{equation}
\begin{equation}
\Omega^{1}=\{\omega=f(t)\Phi(t,z)dt\vert f\hspace{5pt}\text{is a rational function on}\hspace{5pt}{\CC}\backslash \{{z_1,...,z_n}\}\}
\end{equation}
and $d$ is the usual Cartan differential. Under these hypothesis we have the following:
\begin{theor}
For generic $m_1,...,m_n$ and $\kappa$, $H^{0}=0$ and $\dim H^{1}=n-1$. Moreover the differential 
forms $\omega_l=\Phi\displaystyle\frac{dt}{t-z_l},\hspace{5pt}l=1,...,n-1$, form a basis in $H^{1}$.
\end{theor}
\begin{remar}
From the previous theorem follows that $I^{\gamma}\hspace{5pt}\text{for}\hspace{5pt}\gamma=\gamma_1,...,\gamma_{n-1}$, consists of 
integrals of closed forms.
\end{remar}
Let $H_1=(H^1)^{\ast}$ be the dual space of the first cohomology group. Each interval $\gamma_m,
\hspace{3pt}m=1,...,n-1$,
defines a linear functional on $\Omega^1$:
$$[\gamma_m]:f\Phi dt\longmapsto\int_{\gamma_m}f\Phi dt.$$
By the Stokes theorem we have that $\int_{\gamma_m}d(f\Phi)=0$, hence $[\gamma_m]$ defines an element of $H_1$.
We have that:
\begin{theor}
The elements $[\gamma_m]$, $m=1,...,n-1$, form a basis in $H_1$.
\end{theor}
\begin{proof}
The proof of the theorem is based on the following formula:
\begin{equation}
\det_{1\leq l,m\leq n-1}\big(\frac{m_l}{\kappa}\int_{\gamma_m}\Phi\frac{dt}{t-z_l}\big)
\frac{\Gamma(\frac{m_1}{\kappa}+1)\cdot..\cdot\Gamma(\frac{m_n}{\kappa}+1)}{\Gamma(\frac{m_1+...+m_n}{\kappa}+1)}
\prod_{i\neq j}(z_i-z_j)^{\frac{m_j}{\kappa}}.\label{pippo5}
\end{equation}
\end{proof}
\begin{examp}
Let's suppose $n=2$ and $(z_1,z_2)=(0,1)$. Then the formula \ref{pippo5} becomes:
\begin{equation}
b\int_{0}^{1}t^{a}(1-t)^{b-1}dt=\frac{\Gamma(a+1)\Gamma(b+1)}{\Gamma(a+b+1)}.
\end{equation}
\end{examp}
\begin{remar}
In the following table we summarize the main relations between conformal field theory and the geometry of
the classical KZ equation:

\begin{tabular}{|l|l|}
                                                                                                 
\hline 

Conformal Field Theory      &   Geometry of hypergeometric functions                             \\ \hline    

Space of Conformal Blocks   &   Cohomology space $H^1$ associated to the function $\Phi (t,z)$   \\ \hline   

KZ equation                 &   Gauss-Manin connection                                           \\ \hline

Solutions to KZ equation    &   Cycles in $H_1$                                                  \\ \hline

\end{tabular}

\end{remar}

In what follows we proceed to the quantization of the geometry associated to the KZ equation.

\section{Geometric Construction of Quantized KZ Equation}

We know that KZ can be realized as differential equation for integrals of the basic closed differential one
forms over cycles. It turns out that there is a quantization of all those geometric object leading to
a geometric construction of the q-KZ equation. Now we want introduce the p-analogous 
of the main ingredients of the classical picture. Let's start to define the p-analogous of the function 
$\Phi (t,z)=\prod_{l=1}^{n}(t-z_l)^{\frac {m_l}{\kappa}}$.

\underline{The function $\Phi_{p}(t,z)$:} first we observe that the monomial $T^{\frac {m}{\kappa}}$ 
satisfies the following differential equation:
\begin{equation}
\frac{dY}{dT}=\frac{\frac{m}{\kappa}}{T}Y(T).\label{pippo6}
\end{equation}
The analogous of \ref{pippo6} is:
\begin{equation}
y(t+p)=\frac{t+m}{t-m}y(t).\label{pippo8}
\end{equation}
In fact if we pose $t=ST$ and and we define $Y(T)=y(ST)$, then:
\begin{equation}
Y\big(T+\frac{p}{S}\big)=\big(1+\frac{2m}{S}\frac{1}{T-\frac{m}{S}}\big)\label{pippo7}
\end{equation}
and as $S\longrightarrow\infty$ we get that \ref{pippo7} turns into the equation:
\begin{equation}
\frac{dY}{dT}=\frac{\frac{2m}{p}}{T}Y(T).
\end{equation}
\underline{p-Solutions:} the solution of the equation \ref{pippo8} is given by the following:
$$y(t,m)=\Gamma(\frac{t+m}{p})\Gamma(1-\frac{t-m}{p})e^{-\frac{\pi it}{p}},$$
so the p-analogous of the Master function will be given by:
$$\Phi_{(p)}=\prod_{l=1}^{n}y(t-z_l,m_l)=\prod_{l=1}^{n}\Gamma(\frac{t-z_l+m_l}{p})\Gamma(1-\frac{t-z_l+m_l}{p})
e^{-\frac{\pi it}{p}}.$$
\underline{Singularities of the function $\Phi_{p}(t,z)$:} these correspond to the points:
$$t=z_l-m_l-Np$$
$$t=z_l+m_l+(N+1)p$$
where $l=1,...,n$, $N=0,1,2,...$. We will also introduce the following:
\begin{notat}
$\widetilde{\Sing}$=$\{t\in{\CC}\vert\hspace{3pt}t=z_l-m_l+(N+1)p,
	\hspace{3pt}t=z_l+m_l-Np,\hspace{3pt}l=1,...,n,\hspace{3pt}N=0,1,2,...\}$
\end{notat}

\underline{The complex}:
now we define the p-analogous of the complex given in \ref{pippo20}. Let's introduce the space of the 0 and 
1 p-form:\\
$\Omega_{(p)}^{0}=$
$$\Omega_{(p)}^{1}=\{\Phi_{(p)}f(t)\vert\hspace{3pt}\text{$f$ is a rational function
 regular on}\hspace{2pt}{\CC}\backslash\text{$\widetilde{\Sing}$ an with first order poles on it}\},$$

and let the $p$ differential ${\mathcal D}_p:\Omega_{(p)}^{0}\longrightarrow \Omega_{(p)}^{1}$ defined by the following:
\begin{equation}
({\mathcal D}_ph)(t)=h(t+p)-h(t).
\end{equation}
Starting from these data we have the following:
\begin{defin} the $p$ analogous of the complex \ref{pippo20} is: 
\begin{equation}
\begin{CD}
0\longrightarrow\Omega_{p}^{0} @>{\mathcal D}_{p}>>\Omega_{p}^{1}\longrightarrow 0.
\end{CD}
\end{equation}
The first cohomology group is given by:
$$H^{1}_{(p)}=\Omega^{1}_{(p)}\backslash {\mathcal D}_{p}\Omega^{(0)}_{(p)},$$
and the first homology group will be:
$$H_{1}^{(p)}=\big(H^{1}_{(p)}\big)^{\ast}.$$
\end{defin}
\begin{theor} $1$. For generic $z_1,...,z_n, p, m_l$, $H^{1}_{(p)}$ is (n-1)-dimensional and the following 1 forms
form a basis:
$$
\omega_{l}=\frac{1}{t-z_l-m_l}\cdot \frac{t-z_{l-1}+m_{l-1}}{t-z_{l-1}-m_{l-1}}\cdot..\cdot 
\frac{t-z_{1}+m_{1}}{t-z_{1}-m_{1}}dt.$$
$2$. Let $t=ST$, $z_j=SZ_j$, where $S>>1$ and $T$, $Z_j$ are new variables. As $S\longrightarrow\infty$
the p-complex turns into the de Rham complex of $\prod_{l=1}^{n}(T-Z_l)^{\frac{2m_l}{p}}$ and 
$\{\omega_l\}_{l=1,..,n}$ turn into $\big\{\frac{dT}{T-Z_l}\big\}_{l=1,..,n}$ which give a basis in $H^{1}$
of the de Rham complex.
\end{theor}

\subsection{$p$-Homology Theory}
In this subsection we will construct linear functionals on the space $\Omega^{1}_{(p)}\backslash {\mathcal D}_{p}$,
i.e linear functionals:
$$l:\Omega_{(p)}^{1}\longrightarrow {\CC}$$
such that:
$$l(g(t+p)-g(t))=0,\hspace{5pt}\forall g\in\Omega^{1}_{(p)}.$$

\underline{Naive approach}: let $\xi\in{\CC}$ and $h(t)=\Phi_{(p)}(t)f(t)$.
\begin{defin}[Jackson's Integral] The Jackson integral of the function $h=h(t)$ over $[\xi]_{p}$
is given by the following formula:
\begin{equation}
\int_{[\xi]_{p}}h(t)dt=p\sum_{l=-\infty}^{\infty}h(\xi +lp).\label{jack}
\end{equation}
\end{defin}
\begin{prope}
$1$.\hspace{3pt} If $h(t)=g(t+p)-g(t)$ then:\\
$\displaystyle\int_{[\xi]_{p}}h(t)dt=p\sum_{l=\infty}^{\infty}\big(g(\xi +(l+1)p)-g(\xi +lp)\big)=0;$\\
$2$.\hspace{3pt} $\displaystyle{\lim_{p\rightarrow 0}\int_{[\xi]_{p}}h(t)dt=\int_{\xi +\text{line}\hspace{3pt}{\RR}p}h(t)dt.}$
\end{prope}
In our case this definition does not work since the integral \ref{jack} is divergent.\\
\underline{Right approach}: let's choose our data $z_1,...,z_n,m_1,...,m_n,p$ such that the points $z_i$, for
$i=1,...,n$, lie on the imaginary line and $m_1,...,m_n,p\in{\RR}$, with $m_1,...,m_n,p<0$. 
For a given $m\in\{1,...,n-1\}$ let's introduce the following function:

\[
\begin{array}{ccc}
{\GG}_m:{\CC}          &
\longrightarrow        & 
{\CC}                                 \\
                       &             &
                                      \\
    t                  & \longmapsto & 

exp(\frac{2\pi imt}{p})
\end{array}
\]

\begin{remar}
The ${\GG}_m$ are p-periodic functions.
\end{remar}
Now for every $m$ let's consider the following functionals:

\[
\begin{array}{ccc}
 [{\GG}_m]:\Omega_{(p)}^{1}     &
            \longrightarrow     & 
            {\CC}                               \\
                                &             &
                                                \\
         \Phi_{(p)}f            & \longmapsto & 

            \int_{i{\RR}}{\GG}_{m}(t)\Phi_{(p)}(t,z)dt
\end{array}
\]

\begin{remar} 
The functionals $[{\GG}_m]$ are the $p$ analogous of the intervals $[{\gamma}_m]$.
\end{remar}

In what follows we will describe some of the properties of the functionals $[{\GG}_m]$ just introduced.
\begin{prope}
\begin{enumerate}
\item The functionals $[{\GG}_m]$ are well defined if $m\in\{1,...,n-1\}$ and are not defined if
$m\notin\{1,...,n-1\}$.
\item $[{\GG}_m]\vert_{{\mathcal D}_{p}\Omega^{0}_{(p)}}\equiv 0$.
{\rm From the previous property we get the following:}
\begin{corol} The functionals $[{\GG}_m]$, $m\in\{1,..,n-1\}$, define elements in $H_{1}^{(p)}$.
\end{corol}
\item If $z_l=SZ_l$, then: 
$$\lim_{S\rightarrow +\infty}[{\GG}_m]=[z_m,z_{m+1}]$$
in the following sense:
$$\int_{i{\RR}}{\GG}_{m}(t)\Phi_{(p)}(t,z)\omega_ldt\longrightarrow C_m(S,p)\big(\int_{Z_m}^{Z_{m+1}}
\prod_{j=1}^{n}(T-Z_j)^{\frac{2m_j}{p}}\frac{dT}{T-Z_l}+{\mathcal O}(S^{-1})\big), \forall m.$$
Here $C_m=C_m(S,p)$ is some explicitly known function.
\item Let $u_m=u_m(z_1,...,z_n)$, for $m=1,..,n-1$, be the following function:
\begin{equation}
u_m(z_1,...,z_n)=\big(\int_{i{\RR}}{\GG}_{m}(t)\Phi_{(p)}(t,z)\omega_1dt,...,\int_{i{\RR}}{\GG}_{m}(t)\Phi_{(p)}(t,z)\omega_{n-1}dt\big).
\label{pippo30}
\end{equation}
Then for every $m\in\{1,..,n-1\}$, \ref{pippo30} is a solution of the q-KZ equation
with values in $M_{m_1}\otimes,...,\otimes M_{m_n}[\vert m\vert-2]$ with $\mu =0$.
The functions $u_m=u_m(z_1,...,z_n)$ define by \ref{pippo30}
are called {\rm q-hypergeometric} functions associated with $\Phi_{(p)}$.
\item The functionals $[{\GG}_m]$, $m\in\{1,..,n-1\}$, form a basis in the homology space $H_1^{(p)}$.
\end{enumerate}
\end{prope}
\begin{proof}
\begin{enumerate}
\item Follows applying the Stirling formula to the definition of the p-master function $\Phi_{(p)}$:
$$\Phi_{(p)}=\prod_{l=1}^{n}y(t-z_l,m_l)=\prod_{l=1}^{n}\Gamma\big(\frac{t-z_l+m_l}{p}\big)
\Gamma\big(1-\frac{t-z_l+m_l}{p}\big)e^{-\frac{\pi it}{p}}.$$
\item Follows from the periodicity of the functions ${\GG}_m$, in fact since these functions
are p periodic we have:
$$\int_{i{\RR}}{\GG}_{m}(t)\big[\Phi_{(p)}(t+p,z)f(t+p)-\Phi_{(p)}(t,z)f(t)\big]dt=
\int_{i{\RR}+p}{\GG}_{m}(t)\Phi_{(p)}(t,z)f(t)dt$$
$\displaystyle -\int_{i{\RR}}{\GG}_{m}(t)\Phi_{(p)}(t,z)f(t)=0$.
\item Follows from the Stirling formula.
\addtocounter{enumi}{1}
\item This is a consequence of the following formula:
$$
\det_{1\leq l,m\leq n-1}\big(\frac{2m_l}{p}\int_{\RR}{\GG}_{m}\Phi_{(p)}\omega_ldt\big)=$$
$$(2\pi i)^{\frac{n(n-1)}{2}}\frac{\prod_{l=1}^{n-1}\Gamma\big(1+\frac{2m_l}{p}\big)}
{\Gamma\big(1+\frac{2}{p}\sum_{j=1}^{n-1}m_j\big)}
\prod_{j<l}\Gamma\big(\frac{z_j+m_j-z_l+m_l}{p}\big)\Gamma\big(1+\frac{z_l+m_l-z_j+m_j}{p}\big).$$
\end{enumerate}
\end{proof}

\begin{examp} For $n=2$ from the previous determinant we get the following classical formula:
\begin{equation}
\int_{-i\infty}^{i\infty}\Gamma(a+t)\Gamma(a+t)\Gamma(c-t)\Gamma(d-t)dt=2\pi i
\frac{\Gamma(a+c)\Gamma(a+d)\Gamma(b+c)\Gamma(b+d)}{\Gamma(a+b+c+d)}.\label{pippo40}
\end{equation}
\end{examp}

\begin{remar} \ref{pippo40} is the q-analogous of the $\beta$-function and it is called is called 
{\rm Barnes'} formula.
\end{remar}
In the next section we will describe explicit solutions for the q-KZ equation in terms of q-hypergeometric 
functions.

\section{q-KZ Equation: $\slt$ case}

Let's start with some preliminary remarks about the action of the twisted symmetric group $S^k$ on space 
${\mathcal F}_k=\{\text{function of k-variables}\hspace{3pt} t_1,..,t_k\}$.\\
Given $h=h(x)$ a function such that 
\begin{equation}
h(x)h(-x)=1.\label{pippo43}
\end{equation}
let's define:
\begin{equation}
s_j:f(t_1,..,t_k)\longmapsto f(t_1,..,t_{j+1},t_j,..,t_k)h(t_i-t_j).\label{pippo41}
\end{equation}
\begin{defin}
$$S^k=\{s_j:{\mathcal F}_k\longrightarrow {\mathcal F}_k \hspace{2pt}\text{such that \ref{pippo41} is satisfied} \}.$$
\end{defin}

\begin{lemma} The maps $s_j$, $j=1,...,k-1$, given by \ref{pippo41}, define an action of $S^k$ on the space
${\mathcal F}_k$.
\end{lemma}
\begin{proof}
It suffices check that $s_j^2=1$ and $s_js_{j+1}s_{j}=s_{j+1}s_js_{j+1}$ for every 
$j\in\{1,..,k-1\}$.
\end{proof}
\begin{examp}
The following are examples of functions that satisfy the condition \ref{pippo43}.
\begin{equation}
h(x)=\frac{x-1}{x+1},\label{pippo44}
\end{equation}

\begin{equation}
h(x)=\frac{\sin(\frac{\pi(x-1)}{p})}{\sin(\frac{\pi(x+1)}{p})},\label{pippo45}
\end{equation}

\begin{equation}
h(x)=\frac{\theta(x-a,\tau)}{\theta(x+a,\tau)}.\label{pippo46}
\end{equation}
\end{examp}
\begin{remar}
In the formula \ref{pippo46} $\theta=\theta(x,\tau)$ is the first Jacobi $theta$-function.
\end{remar}
\begin{remar}
The functions just described in the previous example enter in the description of the solutions of the
q-KZ equation. In particular the functions \ref{pippo44} and \ref{pippo45} are used in the {\rm rational}
case. For the {\rm trigonometric} case we need the functions \ref{pippo45} and \ref{pippo46}, while for the
elliptic one we need only \ref{pippo46}.
\end{remar}
We need the following:
\begin{defin} Let $F_{m_1,...,m_n}^{k}(z_1,...,z_n)\subset{\mathcal F}_k$ be the ${\CC}$ vector space of the
functions of the form:
$$P(t_1,..,t_k)\prod_{l=1}^{n}\prod_{i=1}^{k}\frac{1}{t_i-z_l-m_l}\prod_{1\leq i<j\leq k}\frac{t_i-t_j}
{t_i-t_j+1},$$
where $P=P(t_1,..,t_k)$ is a symmetric polynomial (in the standard sense), of degree less than $n$ in each 
of the $k$ variables.
\end{defin}
\begin{examp}
For $k=1$ we have that the function in $F_{m_1,...,m_n}^{k}(z_1,...,z_n)$ are of the form:
$$\frac{P(t)}{(t-z_1-m_1)\cdot..\cdot (t-z_n-m_n)},$$
with $deg\hspace{2pt}P<n$.
\end{examp}
\begin{defin}
Let's define:
$$F_{m_1,...,m_n}(z_1,...,z_n)=\bigoplus_{k=0}^{\infty} F_{m_1,...,m_n}^{k}(z_1,...,z_n),$$
where by definition, $F_{m_1,...,m_n}^{0}(z_1,...,z_n)={\CC}$.
\end{defin}
The spaces $F_{m_1,...,m_n}^{k}(z_1,...,z_n)$ have the following:
\begin{prope}
\begin{enumerate}
\item $\dim F_{m_1,...,m_n}^{k}=\dim M_{m_1}\otimes...\otimes M_{m_n}[\vert m\vert-2k]$;
\item $F_{m_1,...,m_n}^{k}$ consists of functions that are symmetric with respect to the action
of $S^{k}$ with $h$ given by \ref{pippo44}.
\end{enumerate}

\begin{examp} If 
$$g=\frac{1}{t_1-z_1-m_1}\cdot\frac{1}{t_2-z_1-m_1}\cdot\frac{t_1-t_2}{t_1-t_2+1},$$
then:
$$s_1(g)=\frac{1}{t_1-z_1-m_1}\cdot\frac{1}{t_2-z_1-m_1}\cdot\frac{t_2-t_1}{t_2-t_1+1}\cdot\frac{t_1-t_2-1}
{t_1-t_2+1}=$$
$$\frac{1}{t_1-z_1-m_1}\cdot\frac{1}{t_2-z_1-m_1}\cdot\frac{t_2-t_1}{t_2-t_1+1}=g.$$
\end{examp}
$3$.\hspace{3pt} 
$F_{m_1,...,m_n}^{k}(z_1,...,z_n)=F_{m_1,..,m_{j+1},m_{j},..,m_n}^{k}(z_1,..,z_{j+1},{j},..,z_n).$

\end{prope}
\begin{examp} For $n=1$ we have that the space $F_{m}^{k}(z)$ is one dimensional and it is spanned by:
\begin{equation}
\omega_{k}(t_1,..,t_k,z)=\prod_{a<b}\frac{t_a-t_b}{t_a-t_b+1}\prod_{a=1}^{k}\frac{1}{t_a-z-m}.
\end{equation}
\end{examp}

Let's introduce the following map:

\[
\begin{array}{ccc}
 F_{m_1,...,m_n'}^{k'}(z_1,...,z_n')\otimes F_{m_{n'+1},...,m_{n'+n''}}^{k''}(z_{n'+1},...,z_{n'+n''}) &
\longrightarrow                                                                                        & 
F_{m_1,...,m_{n'+n''}}^{k'+k''}(z_1,...,z_{n'+n''})                                                                     \\
                                                                                                       &             &
                                                                                                                         \\
    f\otimes g                                                                                         & \longmapsto & 
k(t_1,...,t_{k'+k''})

\end{array}
\]
where

$$k(t_1,...,t_{k'+k''})=\frac{1}{k'!k''!}Sym_{h}\big(f(t_1,..,t_k')g(t_{k'+1},..,t_{k'+k''})\cdot
\prod_{k<j\leq k'+k''}\cdot\prod_{1\leq l\leq n'}\frac{t_j-z_l+m_l}{t_j-z_l-m_l}\big)$$
\begin{examp}
For $k'=0$ and $k''=1$ we have:

\[
\begin{array}{ccc}
 F_{m_1}^{0}(z_1)\otimes F_{m_2}^{1}(z_2)                                           &
\longrightarrow                                                                     & 
F_{m_1,m_2}^{1}(z_1,z_2)                                                                \\
                                                                      &             &
                                                                                        \\
1\otimes \frac{1}{t-z_2-m_2}                                          & \longmapsto & 
\frac{t-z_1+m_1}{t-z_1-m_1}\frac{1}{t-z_2-m_2}

\end{array}
\]
and
\[
\begin{array}{ccc}
 F_{m_1}^{1}(z_1)\otimes F_{m_2}^{0}(z_2)                                         &
\longrightarrow                                                                   & 
F_{m_1,m_2}^{1}(z_1,z_2)                                                                \\
                                                                    &             &
                                                                                        \\
\frac{1}{t-z_1-m_1}\otimes 1                                        & \longmapsto & 
\frac{1}{t-z_1-m_1}

\end{array}
\]
\end{examp}

\begin{claim} For generic values of $z_1,..,z_m$:
$$
\bigoplus_{k'+k''=k} 
F_{m_1,...,m_n'}^{k'}(z_1,...,z_n')\otimes F_{m_{n'+1},...,m_{n'+n''}}^{k''}(z_{n'+1},...,z_{n'+n''})
\longrightarrow F_{m_1,...,m_{n'+n''}}^{k}(z_1,...,z_{n'+n''})$$
is an isomorphism.
\end{claim}

Now we can define:\\

\begin{defin}
\[
\begin{array}{ccc}
 \omega(z): (V_m)^{\ast}        &
            \longrightarrow     & 
            F_{m}(z)                               \\
                                &             &
                                                   \\
             (f^kv)^{\ast}      & \longmapsto & 
\omega_{k}
\end{array}
\]
\end{defin}
The map $\omega_z$ induces a map:
\begin{equation}
\omega(z_1,..,z_n):V_{m_1}^{\ast}\otimes...\otimes V_{m_n}^{\ast}\longrightarrow F_{m_1,...,m_n}(z_1,..,z_n)
.\label{pippo60}
\end{equation}
\begin{claim} For generic $z_1,..,z_n$, the map $\omega=\omega(z_1,...,z_n)$, \ref{pippo60} is an 
isomorphism.
\end{claim}

\begin{reftc}
\cite{V5}, \cite{TV1,TV2}
\end{reftc}


\chapter{}
The last part of the last lecture was devoted to discuss the following problem:
how to solve the q-KZ equation in terms of q-hypergeometric functions.
Let's recall the main points of this construction.
We introduced the action of the twisted symmetric group $S^k$ in the space of 
the function ${\mathcal F}_k=\{\text{functions of k variables}\hspace{2pt}t_1,...,t_k\}$, via the maps:
\begin{equation}
s_j:f(t_1,..,t_k)\longmapsto f(t_1,..,t_{j+1},t_j,..,t_k)h(t_i-t_j),
\end{equation}
where $h=h(x)$ is a function that satisfies the condition: 
$$h(x)\cdot h(-x)=1,$$
and we also defined the space:
$$F_{m_1,...,m_n}^{k}(z_1,...,z_n)\subset{\mathcal F}_k$$
of the functions of the form:
$$P(t_1,..,t_k)\prod_{l=1}^{n}\prod_{i=1}^{k}\frac{1}{t_i-z_l-m_l}\prod_{1\leq i<j\leq k}\frac{t_i-t_j}
{t_i-t_j+1},$$
where $P=P(t_1,..,t_k)$ is a symmetric polynomial (in the standard sense).\\
We also noticed that the space 
$F_{m_1,...,m_n}(z_1,...,z_n)=\bigoplus_{k=0}^{\infty} F_{m_1,...,m_n}^{k}(z_1,...,z_n)$ satisfies the 
following:
\begin{prope}
\begin{enumerate}
\item $\dim F_{m_1,...,m_n}^{k}=\dim M_{m_1}\otimes...\otimes M_{m_n}[\vert m\vert-2k]$;
\item $F_{m_1,...,m_n}^{k}$ consists of functions that are symmetric with respect to the action
of $S^{k}$ with $\displaystyle h=\frac{x-1}{x+1}$;
\item $F_{m_1,...,m_n}^{k}(z_1,...,z_n)=F_{m_1,..,m_{j+1},m_{j},..,m_n}^{k}(z_1,..,z_{j+1},z_{j},..,z_n).$
\end{enumerate}
\end{prope}
Then we defined a tensor product:

\[
\begin{array}{ccc}
 F_{m_1,...,m_{n'}}^{k'}(z_1,...,z_{n'})\otimes F_{m_{n'+1},...,m_{n'+n''}}^{k''}(z_{n'+1},...,z_{n'+n''}) &
\longrightarrow                                                                                        & 
F_{m_1,...,m_{n'+n''}}^{k'+k''}(z_1,...,z_{n'+n''})                                                                     \\
                                                                                                       &             &
                                                                                                                         \\
    f\otimes g                                                                                         & \longmapsto & 
k(t_1,...,t_{k'+k''})

\end{array}
\]
where

$$k(t_1,...,t_{k'+k''})=\frac{1}{k'!k''!}Sym_{h}\big(f(t_1,..,t_{k'})g(t_{k'+1},..,t_{k'+k''})\cdot
\prod_{k<j\leq k'+k''}\cdot\prod_{1\leq l\leq n'}\frac{t_j-z_l+m_l}{t_j-z_l-m_l}\big)$$

and we state the following:
\begin{claim} For generic values of $z_1,..,z_n$:
$$
\bigoplus_{k'+k''=k} 
F_{m_1,...,m_{n'}}^{k'}(z_1,...,z_{n'})\otimes F_{m_{n'+1},...,m_{n'+n''}}^{k''}(z_{n'+1},...,z_{n'+n''})
\longrightarrow F_{m_1,...,m_{n'+n''}}^{k}(z_1,...,z_{n'+n''})$$
is an isomorphism.
\end{claim}
Finally we defined:
\[
\begin{array}{ccc}
 \omega(z): (V_m)^{\ast}        &
            \longrightarrow     & 
            F_{m}(z)                               \\
                                &             &
                                                   \\
             (f^kv)^{\ast}      & \longmapsto & 
\omega_{k}
\end{array}
\]
where 
\begin{equation}
\omega_{k}(t_1,..,t_k,z)=\prod_{a<b}\frac{t_a-t_b}{t_a-t_b+1}\prod_{a=1}^{k}\frac{1}{t_a-z-m},\label{titti1}
\end{equation}
is a generator of the one dimensional vector space $F_{m}^{k}(z)$ and from \ref{titti1} and the tensor 
product we also defined the map:
\begin{equation}
\omega(z_1,..,z_n):V_{m_1}^{\ast}\otimes...\otimes V_{m_n}^{\ast}\longrightarrow F_{m_1,...,m_n}(z_1,..,z_n);
\label{titti2}
\end{equation}
that has the following property:
\begin{claim} For generic $z_1,..,z_n$, the map $\omega=\omega(z_1,...,z_n)$, defined in \ref{titti2}, is an 
isomorphism.
\end{claim}
Now we will proceed introducing a suitable R-matrix to relate basis in different spaces.
Let $z,u\in{\CC}$ such that the map:
\begin{equation}
\omega(z,u):M_{a}^{\ast}\otimes M_{b}^{\ast}\longrightarrow F_{a,b}(z,u)\label{titti0}
\end{equation}
be invertible. Then:
\begin{defin}
The {\rm rational $R$-matrix} $R_{a,b}(z,u)$ is the dual map of the composition:
\begin{equation}
\begin{CD}
M_{a}^{\ast}\otimes M_{b}^{\ast} @>p>> M_{b}^{\ast}\otimes M_{a}^{\ast}
@>\omega(u,z)>> F_{a,b}(z,u) @>\omega^{-1}(z,u)>> M_{a}^{\ast}\otimes M_{b}^{\ast},\label{titti3}
\end{CD}
\end{equation}
\end{defin}
\begin{remar} Alternatively we can think the map $R_{ab}(z,u)$ defined by \ref{titti3} as the matrix that 
express the basis $\omega(u,z)\big((f^mv_b)^{\ast}\otimes((f^lv_a)^{\ast}\big)$ in terms of the basis
$\omega(z,u)\big((f^lv_a)^{\ast}\otimes(f^mv_b)^{\ast}\big)$.
\end{remar}
The rational R-matrix just introduced satisfies the following:
\begin{prope}
$1.$\hspace{3pt}$R_{ab}(z,u)$ depends only on the difference $u-z$;\\
$2.$\hspace{3pt}$R_{ab}(z,u)$ is the standard {\rm R-matrix} associated with the Yangian module structure
of the tensor product $M_a(z)\otimes M_b(u)$.
\end{prope}
In the next section we will describe the trigonometric analogous of the previous construction.
\section{Trigonometric case}
In what follows with $p$ we will denote the step of the q-KZ equation, $q=e^{\frac{\pi i}{p}}$
and $M_{m_j}^{q}$, $j=1,..,n$, will be a ${\mathcal U}_{q}(\slt)$-Verma module.
Let's first introduce the trigonometric analogous of the space
$$F_{m_1,...,m_n}^{k}(z_1,...,z_n)\subset{\mathcal F}_k$$
\begin{defin}Let ${\mathcal F}_{m_1,...,m_n}^{k}(z_1,...,z_n)\subset{\mathcal F}_k$ the space of  function
of the form:
\begin{equation}
{\mathcal P}(\xi_1,..,\xi_l)\prod_{l=1}^{n}\prod_{i=1}^{k}\frac{\exp(\frac{\pi i(z_l-t_a)}{p})}
{\sin(\frac{\pi(t_a-z_l-m_l)}{p})}\prod_{a<b}\frac{\sin(\frac{\pi(t_a-t_b)}{p})}{\sin(\frac{\pi(t_a-t_b+1)}{p})},
\end{equation}
where $\xi_a=\exp(\frac{2\pi it_a}{p})$ and ${\mathcal P}={\mathcal P}(\xi_1,..,\xi_l)$ is a symmetric polynomial
in $\xi_1,...,\xi_l$, in the standard sense of degree $n$ in each of its variable.
\end{defin}
In analogy to the rational case we have:
\begin{defin}
$${\mathcal F}_{m_1,...,m_n}(z_1,...,z_n)=\bigoplus_{k=0}^{\infty}{\mathcal F}_{m_1,...,m_n}^{k},$$
where, by definition, ${\mathcal F}_{m_1,...,m_n}^{0}={\CC}$.
\end{defin}
The spaces ${\mathcal F}_{m_1,...,m_n}^{k}$ satisfy the following:
\begin{prope}
\begin{enumerate}
\item $\dim{\mathcal F}^{k}={n+k-1\choose k}=\dim(M_{m_1}^{q}
\otimes..\otimes M_{m_n}^{q})[\vert m\vert-2k]$;\\
\item ${\mathcal F}^k$ consists of symmetric functions with respect to the action of the group $S^k$
defined by the function:
$${\tilde h}(x)=\frac{\sin(\frac{\pi(x-1)}{p})}{\sin(\frac{\pi(x+1)}{p})};$$
\item ${\mathcal F}_{m_1,...,m_n}^{k}(z_1,...,z_n)
	={\mathcal F}_{m_1,..,m_{j+1},m_{j},..,m_n}^{k}(z_1,..,z_{j+1},z_{j},..,z_n).$
\end{enumerate}
\end{prope}
\begin{examp} As in the rational case we have an explicit description of the generators for the case $n=1$.  
In fact in this hypothesis we have that the space ${\mathcal F}_{m}^{k}(z)$ is generated by:
\begin{equation}
W_k(t_1,...,t_k,z)=\prod_{a=1}^{k}\frac{\exp(\frac{\pi i(z-t_a)}{p})}{\sin(\frac{\pi(t_a-z-m)}{p})}
\prod_{a<b}\frac{\sin(\frac{\pi(t_a-t_b)}{p})}{\sin(\frac{\pi(t_a-t_b+1)}{p})}.
\end{equation}
\end{examp}
Let's introduce the tensor product:

\[
\begin{array}{ccc}
{\mathcal F}_{m_1,...,m_{n'}}^{k'}(z_1,...,z_{n'})\otimes {\mathcal F}_{m_{n'+1},...,m_{n'+n''}}^{k''}(z_{n'+1},...,z_{n'+n''}) &
\longrightarrow                                                                                                     & 
{\mathcal F}_{m_1,...,m_{n'+n''}}^{k'+k''}(z_1,...,z_{n'+n''})                                                                          \\
                                                                                                                    &             &
                                                                                                                                    \\
    f\otimes g                                                                                                      & \longmapsto & 
u(t_1,...,t_{k'+k''})

\end{array}
\]

where:

$$u(t_1,...,t_{k'+k''})=\frac{1}{k'!k''!}Sym_{\tilde {h}}\big(f(t_1,..,t_{k'})g(t_{k'+1},..,t_{k'+k''})
\prod_{k'<j\leq k'+k''}\prod_{1\leq l\leq n'}\frac{\sin(\frac{t_j-z_l+m_l}{p})}
{\sin(\frac{t_j-z_l-m_l}{p})}\big)$$
Also in the this case we have that:
\begin{claim} For generic values of $z_1,..,z_n$:
$$
\bigoplus_{k'+k''=k} 
{\mathcal F}_{m_1,...,m_{n'}}^{k'}(z_1,...,z_{n'})\otimes {\mathcal F}_{m_{n'+1},...,m_{n'+n''}}^{k''}(z_{n'+1},...,z_{n'+n''})
\longrightarrow {\mathcal F}_{m_1,...,m_{n'+n''}}^{k}(z_1,...,z_{n'+n''})$$
is an isomorphism.
\end{claim}
Now let's define:
\vspace{5pt}
\begin{defin}
\begin{equation}
\begin{array}{ccc}
W(z): (M_m)^{q}                 &
            \longrightarrow     & 
            F_{m}(z)                               \\
                                &             &
                                                   \\
             (f^kv)             & \longmapsto & 
C_{k}W_k(t_1,..,t_n,z),
\end{array}\label{titti4}
\end{equation}
where:
$$C_k=\prod_{s=0}^{k-1}\frac{\sin(\frac{\pi(2m-s)}{p})}{\sin(\frac{\pi}{p})}=\prod_{s=0}^{k-1}
\frac{q^{2(m-s)}-q^{-2(m-s)}}{q-q^{-1}}.$$
\end{defin}

The map \ref{titti4} induces a map:
\begin{equation}
W(z_1,..,z_n):M_{m_1}\otimes..\otimes M_{m_n}\longrightarrow {\mathcal F}_{m_1,..,m_n}(z_1,..,z_n).\label{titti5}
\end{equation}
We have:
\begin{claim} For generic $z_1,..,z_n$ the map \ref{titti5} is non degenerate.
\end{claim}
Let's now introduce the trigonometric version of the R-matrix \ref{titti0}. Let $z,u\in{\CC}$ such that
$$W(u,z):V_b\otimes V_a\longrightarrow {\mathcal F}_{a,b}(z,u)$$
is non degenerate. Then:
\begin{defin}
The {\rm trigonometric R-matrix} $R_{a,b}^{\text{trig}}$ is given by following composition:
\begin{equation}
\begin{CD}
M_{a}\otimes M_{b} @> W(z,u)>> {\mathcal F}_{a,b}(z,u)
@> W^{-1}(u,z)>> M_{b}\otimes M_{a} @> p >> M_{a}\otimes M_{b},\label{titti6}
\end{CD}
\end{equation}
\end{defin}
The trigonometric R-matrix defined by \ref{titti6} has the following:
\begin{prope}
$1.$\hspace{3pt} $R_{a,b}^{\text{trig}}$ is a map between $M_a^{q}\otimes M_b^{q}$ and itself;\\
$2.$\hspace{3pt} $R_{a,b}^{\text{trig}}$ depends only on the difference $z-u$, is p-periodic and coincides
with the trigonometric R-matrix $R_{M_{a}^{q},M_{b}^{q}}\big(\exp(\frac{2\pi i(z-u)}{p})\big)$.
\end{prope}
\begin{remar} For any highest weight ${\mathcal U}_{q}(\slt)$-modules $V_1^q$ and $V_2^q$, there is 
$R_{V_1^q,V_2^q}^{q}(z)\in\text{End}(V_1^q\otimes V_2^q)$ satisfying the QYBE:
$$R_{V_1^q,V_2^q}^{q}(\frac{z_1}{z_2})R_{V_1^q,V_3^q}^{q}(\frac{z_1}{z_3})R_{V_2^q,V_3^q}^{q}(\frac{z_2}{z_3})
=R_{V_2^q,V_3^q}^{q}(\frac{z_2}{z_3})R_{V_1^q,V_3^q}^{q}(\frac{z_1}{z_3})R_{V_1^q,V_2^q}^{q}(\frac{z_1}{z_2}).$$
$R_{V_1^q,V_2^q}^{q}(z)$ is called trigonometric R-matrix, it is constructed
using ${\mathcal U}_{q}(\widehat{\slt})$ and normalized requiring that:
$$R_{V_1^q,V_2^q}^{q}(z)(v_1\otimes v_2)=v_1\otimes v_2$$
where $v_i\in V^q_i$ are highest weight vectors.
\end{remar}
\begin{examp}
Let $V_1^q=V_2^q={\CC}^2$, then:
\begin{equation}
R^{\text{trig}}(x)=
\begin{pmatrix} 1 & 0 & 0 & 0 \\ 0 & \frac{x-1}{xq-q^{-1}} & \frac{q-q^{-1}}{xq-q^{-1}} & 0 \\
 0 & \frac{x(q-q^{-1})}{xq-q^{-1}} & \frac{x-1}{xq-q^{-1}} & 0 \\ 0 & 0 & 0 & 1 \end{pmatrix}\quad
\end{equation}
\end{examp}
\begin{remar} There is an analogous construction for the
{\rm elliptic q-KZ equation.}
\end{remar}
\section{Hypergeometric Pairing}
In what follows we will assume that $p\in{\RR}$ $(p<0)$. Let's introduce the following map:
\begin{defin}
$$J(z_1,...,z_n):{\mathcal F}_{m_1,..,m_n}(z_1,..,.z_n)\otimes F_{m_1,..,m_n}(z_1,...,z_n)\longrightarrow {\CC}$$
such that:
\begin{equation}
J(z_1,...,z_n)(f\otimes g)=\int_{\gamma}\Phi_{m_1,..,m_n}(t_1,..,t_k,z_1,..,z_n)f(t_1,..,t_k)g(t_1,..,t_k)dt_1\cdot..\cdot dt_k,
\label{titti7}
\end{equation}
where $\Phi_{m_1,..,m_n}(t_1,..,t_k,z_1,..,z_n)$ is a q-deformation of the Master function, defined by:

$$\Phi_{m_1,..,m_n}(t_1,..,t_k,z_1,..,z_n)=\exp(\mu\sum_{a=1}^{k}t_a/p)\prod_{l=1}^n\prod_{a=1}^k
\frac{\Gamma((t_a-z_l+m_l)/p)}{\Gamma((t_a-z_l-m_l)/p)}\displaystyle\prod_{a<b}
\frac{\Gamma((t_a-t_b-1)/p)}{\Gamma((t_a-t_b+1)/p)},$$
and the cycle $\gamma$ is defined by:
$$\gamma=\{(t_1,..,t_k)\in{\CC}^{k}\vert\hspace{2pt}\Re\hspace{1pt}t_i=0,\hspace{2pt}\forall i=1,..,k\}.$$
\end{defin}
\begin{remar} If $m_1,..,m_n\in{\RR}\hspace{2pt}\text{and}<<0$ and $z_1,..,z_n\in i{\RR}$ the integral 
\ref{titti7} is convergent. If the previous conditions are not satisfied we define \ref{titti7} by analytic
continuation.
\end{remar}
The hypergeometric pairing has the following property:
\begin{propo}
$1.$\hspace{3pt} For generic $z_1,...,z_n$ and generic $m_1,...,m_n$ \ref{titti7} is non degenerate.
\end{propo}
\begin{proof}
The proof follows from the calculation of suitable determinants of hypergeometric integrals.
\end{proof}
\begin{examp}[q-Selberg integral]
Let $n=1$. Then:
$$\int_{-i\infty}^{i\infty}\cdots\int_{-i\infty}^{i\infty}
\prod_{i=1}^{k}\big((u^{2t_i}\Gamma(a+t_i)\Gamma(a-t_i)\big)
\prod_{i<j}\frac{\Gamma(t_i-t_j+x)\Gamma(t_j-t_i+x)}{\Gamma(t_i-t_j)\Gamma(t_j-t_i)}dt_1\cdots dt_k=$$
$$(2\pi \sqrt{-1})^k(u+u^{-1})^{-k(2a+(k-1)x)}\prod_{i=1}^k\frac{\Gamma(1+ix)\Gamma(2a+(i-1)x)}{\Gamma(1+x)}$$
\end{examp}
\begin{remar} The q-Selberg integral defined in the previous example gives the multidimensional 
generalization of the Barnes' formula.
\end{remar}
Now let's consider:
$$J(z_1,...,z_n):{\mathcal F}_{m_1,..,m_n}(z_1,..,.z_n)\otimes F_{m_1,..,m_n}(z_1,...,z_n)\longrightarrow {\CC}$$
and let $\nu$, $\tau\in S^n$, where $S^n$ is the permutation group. Then $J$ defines the following paring:
\begin{equation}
(M^{q}_{m_{\tau_{1}}}\otimes...\otimes M^{q}_{m_{\tau_{n}}})\otimes(M_{m_{\nu_{1}}}\otimes...\otimes M_{m_{\nu_{n}}})^{\ast}\longrightarrow {\CC}
\end{equation}
or, equivalently, the following map:
\begin{equation}
I_{\tau,\nu}(z_1,..,z_n):M^{q}_{m_{\tau_{1}}}\otimes...\otimes M^{q}_{m_{\tau_{n}}}\longrightarrow M_{m_{\nu_{1}}}\otimes...\otimes M_{m_{\nu_{n}}}.
\label{titti8}
\end{equation}
Let's now consider \ref{titti8} for $\nu=id$:
\begin{equation}
I_{\tau,id}(z_1,..,z_n)=I_{\tau}(z_1,..,z_n):M^{q}_{m_{\tau_{1}}}\otimes...\otimes M^{q}_{m_{\tau_{n}}}\longrightarrow M_{m_1}\otimes...\otimes M_{m_n}.
\label{titti9}
\end{equation}
We can now state the following:
\begin{theor} For any $\tau\in S^n$, $I_{\tau}$ defined in \ref{titti9} satisfied the {\rm rational q-KZ equation.}
with values in $M_{m_1}\otimes...\otimes M_{m_n}$.
\end{theor}

\begin{remar}
We can rephrase the content of the previous theorem saying that the solution of the rational q-KZ equation,
with values in $M_{m_1}\otimes...\otimes M_{m_n}$, product of $\slt$-modules, are labeled elements of the
product of ${\mathcal U}_q(\slt)$ corresponding modules.
\end{remar}
We observe that for $\tau\in S^n$, the product $M^{q}_{m_{\tau_{1}}}\otimes...\otimes M^{q}_{m_{\tau_{n}}}$
has a natural basis given by:
\begin{equation}
f_{J}v_{\tau}=f^{j_1}v^{q}_{m_{\tau_1}}\otimes\cdot..\cdot\otimes f^{j_n}v^{q}_{m_{\tau_n}},\label{titti10}
\end{equation}
where $J=(j_1,..,j_n)$ and $v^q_{m_i}\in M_{m_i}^q$ is the generating vector. The basis 
$\big\{f_{J}v_{J}\big\}_{J}$ given in \ref{titti10} corresponds to a solution of the q-KZ equation and for
any element $\tau\in S^n$ we have a basis of solutions.
It turns out that these basis have remarkable asymptotic properties.
\section{Quantization of the Drinfeld-Kohno theorem}
In this section we will discuss the discrete analogue of the Drinfeld-Kohno theorem. Let's start describing
the monodromy of a difference equation. 
\begin{examp} Let's consider the following scalar difference equation:
\begin{equation}
\psi(z+p)=K(z)\psi(z),\hspace{2pt}z\in{\CC}\hspace{3pt}\text{and}\hspace{3pt}p<0.\label{titti11}
\end{equation}
We have that if $\psi=\psi(z)$ is a non zero solution then any other solution $\tilde\psi$ of this equation has the
form:
$${\tilde\psi}(z)=P(z)\psi(z),$$
where $P=P(z)$ is p-periodic, i.e:
$$P(z+p)=P(z)\hspace{3pt}\forall z.$$
Let's now let's suppose that the function $K=K(z)$ satisfies the following properties:
\begin{equation}
K(z)=e^{az}(1+\alpha/z+{\mathcal O}(1/z)),\hspace{3pt}z\rightarrow+\infty,\label{titti12}
\end{equation}
\begin{equation}
K(z)=e^{bz}(1+\beta/z+{\mathcal O}(1/z)),\hspace{3pt}z\rightarrow-\infty,\label{titti13}
\end{equation}
where $a,b,\alpha$ and $\beta$ are some numbers. Then we have that there are solutions $f_{+}$ and
$f_{-}$ of \ref{titti11} such that:
$$f_{+}(z)\sim e^{az/p}z^{\alpha/p}(1+{\mathcal O}(1/z))\hspace{3pt}\text{as}\hspace{3pt}z\rightarrow+\infty$$
and similarly
$$f_{-}(z)\sim e^{bz/p}z^{\beta/p}(1+{\mathcal O}(1/z))\hspace{3pt}\text{as}\hspace{3pt}z\rightarrow-\infty.$$
The function:
\begin{equation}
P_{+,-}=P_{+,-}(z)=\frac{f_+(z)}{f_-(z)}
\end{equation}
is called the {\rm scattering} or {\rm transition} matrix associated to the solutions \ref{titti12} and
\ref{titti13}.
\end{examp}
\begin{remar} The scattering matrix replaces the concept of monodromy in the case of difference equations.
\end{remar}
We want to apply a procedure similar to the one described in the previous example to the q-KZ equation.
\begin{remar} We observe that the q-KZ equation has $n!$ asymptotic zones:
$$A_{\tau},\hspace{3pt}\tau\in S^n:$$
$$A_{\tau}:\hspace{4pt}\Re z_{\tau_1}<<\Re z_{\tau_2}<<....\Re z_{\tau_n},$$
i.e. we can approach $\infty$ in $n!$ different ways.
\end{remar}
\begin{defin}
We will say that a basis of solutions $\psi_1,\hspace{2pt}\psi_2,...$ of the q-KZ equation, forms an 
asymptotic solution in a given zone $A$, if:
\begin{equation}
\psi_j=\exp(\sum_{m=1}^{n}a_{m_j}z_{m}/p)\prod_{m<l}(z_l-z_m)^{b_{jlm}}(v_j+{\mathcal O}(1))\hspace{2pt}\forall j,
\end{equation}
where $a_{m_j}$, $b_{jlm}$ are numbers, $v_1$, $v_2$,... are vectors forming a basis in $V$ and 
${\mathcal O}(1)$ tends to zero as $z\rightarrow\infty$ in $A$.
\end{defin}
\begin{remar} We recall that $K_m=K_m(z_1,..,z_n,p)$ can be written in terms of R-matrices as:
$$
K_m=R_{V_m, V_{m- 1}}(z_m-z_{m-1}+p)...R_{V_m, V_{1}}(z_m-z_{1}+p)e^{\mu h^{(m)}}$$
$$R_{V_m, V_{n}}(z_m-z_n)...R_{V_m, V_{m+1}}(z_m-z_{m+1}).$$ 
In every asymptotic zone $A_\tau$ we have:
$$K_m(z_1,..,z_n)\sim e^{\mu h^{(m)}},$$
i.e at $\infty$ the $K_m(z_1,..,z_n)$ operators are diagonal and the vectors $v_1,v_2,...$ introduced in the
previous definition are common eigenvectors. These vectors can be easily described, namely consider the vectors
$\big\{f^{j_1}v_{m_1}\otimes..\otimes f^{j_n}v_{m_n}\big\}$: they form an eigenbasis for $h^{(m)}$. This basis
will be called the monomial basis in $M_{m_1}\otimes..\otimes M_{m_n}$.
\end{remar}
Now we can state the main result:
\begin{theor}
For every asymptotic zone $A_{\tau}$ the basis in the space of solutions corresponding to the basis:
$$\big\{f_{J}v_{\tau}\big\}\subset M_{m_{\tau_1}}\otimes..\otimes M_{m_{\tau_n}}$$
is asymptotic in the zone $A_{\tau}$.
\end{theor}
We also have the following description of the scattering matrices:
\begin{corol}
The transition functions between the asymptotic solutions corresponding to neighboring zones are given in 
terms of the corresponding trigonometric R-matrices. Namely $\forall\tau\in S^n$, $m\in\{1,..,n\}$, let
\begin{equation}
T_{\tau,m}(z_1,..,z_n):M_{\tau_1}^{q}\otimes..\otimes M_{\tau_{m+1}}^{q}\otimes M_{\tau_{m}}^{q}\otimes..\otimes M_{\tau_{n}}^{q}
\longrightarrow M_{\tau_{1}}^{q}\otimes..\otimes M_{\tau_{n}}^{q}
\end{equation}
be the composition of 
$\displaystyle R_{M^q_{\tau_{m+1}},M^q_{\tau_{m}}}(e^{2\pi i(z_{\tau_{m+1}}-z_{\tau_{m}})/p})$ and:
$$P:M^q_{\tau_{m+1}}\otimes M^q_{\tau_{m}}\longrightarrow M^q_{\tau_{m}}\otimes M^q_{\tau_{m+1}}.$$
Then the transition function between the solutions associated to the zones $A_{\tau}$ and $A_{\sigma_{m}\tau}$
is identified with $T_{\tau,m}(z_1,..,z_n)$.
\end{corol}
\begin{examp}
Let $a,b,c,d,\theta\in{\CC}$ with $\Re \theta>0$.\\
Set $\displaystyle\lambda=\sqrt{a^2-bc}$, $\displaystyle h=\begin{pmatrix} 1 & 0 \\ 0 & 1 \end{pmatrix}\quad $ and
$\displaystyle A(u)=\frac{1}{d+u}\begin{pmatrix} a+u & b \\ c & a-u \end{pmatrix}.\quad $\\
We have that:
$$A(u,\theta)=\theta^{uh}A(u)\theta^{-uh}.$$
We have the following:
\begin{theor}
$$\lim_{s\rightarrow\infty}\Big(s^{-ah}h^s\big(\prod_{r=-s}^{s}A(u+r;\theta)\big)h^ss^{ah}\Big)=A^q(u),$$
where:
$$A^q(u)=\frac{1}{\sin(\pi(d+u))}
\begin{pmatrix} \sin(\pi(a+u))
	& \frac{\pi b(\theta+\theta^{-1})^{2a}}{\Gamma(1+a+\lambda)\Gamma(1+a-\lambda)} \\
\frac{\pi c(\theta+\theta^{-1})^{-2a}}{\Gamma(1-a+\lambda)\Gamma(1-a-\lambda)}
	& \sin(\pi(a-u))
\end{pmatrix}$$
\end{theor}
\end{examp}

\begin{reftc}
\cite{TV1,TV2}
\end{reftc}

\end{document}